\documentclass[A4paper, 12pt]{article}
\usepackage[utf8]{inputenc}
\usepackage[OT2, T1]{fontenc}
\usepackage{cite}
\usepackage{amsmath}
\usepackage{amsfonts}
\usepackage{amssymb}
\usepackage{amsxtra}
\usepackage{mathrsfs}
\usepackage{tensor}
\usepackage{calligra}
\usepackage{scalerel,stackengine}
\stackMath
\usepackage[OT2, T1]{fontenc}
\usepackage[thmmarks, amsmath, thref, amsthm]{ntheorem}
\usepackage[nottoc]{tocbibind}
\usepackage[pagebackref=true]{hyperref}

\usepackage{graphicx}
\usepackage[all]{xy}
\usepackage[portrait, top=2cm, bottom=2cm, left=2cm, right=2cm]{geometry}
\usepackage{enumerate}

\numberwithin{equation}{subsection}

\theoremstyle{plain}
\newtheorem{thm}{Theorem}[section]
\newtheorem{lm}[thm]{Lemma}
\newtheorem{prop}[thm]{Proposition}
\newtheorem{cor}[thm]{Corollary}

\newtheorem{customthm}{Theorem}

\theoremstyle{definition}

\newtheorem{rmk}[thm]{Remark}
\newtheorem{ex}[thm]{Example}
\newtheorem{question}{Question}

\newcommand{\pair}[1]{\left\langle #1 \right\rangle}
\newcommand{\ol}[1]{\overline{#1}}
\newcommand{\ps}[1]{(\!(#1)\!)}
\newcommand{\wh}[1]{\widehat{#1}}
\newcommand\wc[1]{
\savestack{\tmpbox}{\stretchto{
  \scaleto{
    \scalerel*[\widthof{\ensuremath{#1}}]{\kern-.3pt\bigwedge\kern-.3pt}
    {\rule[-\textheight/2]{1ex}{\textheight}}
  }{\textheight}
}{0.35ex}}
\stackon[1pt]{#1}{\scalebox{-1}{\tmpbox}}
}
\newcommand{\sqdot}{\mathbin{\vcenter{\hbox{\rule{.3ex}{.3ex}}}}}

\newcommand{\Abb}{\mathbb{A}}
\newcommand{\Cbb}{\mathbb{C}}
\newcommand{\Gbb}{\mathbb{G}}
\newcommand{\Hbb}{\mathbb{H}}
\newcommand{\Lbb}{\mathbb{L}}
\newcommand{\Pbb}{\mathbb{P}}
\newcommand{\Qbb}{\mathbb{Q}}
\newcommand{\Rbb}{\mathbb{R}}
\newcommand{\Zbb}{\mathbb{Z}}

\newcommand{\Ccal}{\mathcal{C}}
\newcommand{\Dcal}{\mathcal{D}}
\newcommand{\Fcal}{\mathcal{F}}
\newcommand{\Gcal}{\mathcal{G}}
\newcommand{\Tcal}{\mathcal{T}}
\newcommand{\Ocal}{\mathcal{O}}
\newcommand{\Xcal}{\mathcal{X}}
\newcommand{\Ycal}{\mathcal{Y}}

\renewcommand{\H}{\mathrm{H}}
\newcommand{\h}{\mathrm{h}}

\newcommand{\Hscr}{\mathscr{H}}
\newcommand{\Uscr}{\mathscr{U}}
\newcommand{\Vscr}{\mathscr{V}}

\newcommand{\ab}{\operatorname{ab}}
\newcommand{\Ab}{\mathbf{Ab}}
\newcommand{\Aut}{\operatorname{Aut}}
\newcommand{\AV}{\operatorname{AV}}
\newcommand{\Br}{\operatorname{Br}}
\newcommand{\can}{\operatorname{can}}
\newcommand{\cd}{\operatorname{cd}}
\newcommand{\cores}{\operatorname{cor}}
\newcommand{\cts}{\operatorname{cts}}
\renewcommand{\div}{\operatorname{div}}
\newcommand{\Div}{\operatorname{Div}}
\newcommand{\et}{\textnormal{\'et}}
\newcommand{\Ext}{\operatorname{Ext}}
\newcommand{\Gal}{\operatorname{Gal}}
\newcommand{\GL}{\operatorname{GL}}
\newcommand{\Hom}{\operatorname{Hom}}
\newcommand{\id}{\operatorname{id}}
\newcommand{\Img}{\operatorname{Im}}
\newcommand{\Ker}{\operatorname{Ker}}
\newcommand{\lien}{\mathbf{lien}}
\newcommand{\loc}{\operatorname{loc}}
\newcommand{\nr}{\operatorname{nr}}
\newcommand{\Pic}{\operatorname{Pic}}
\newcommand{\PT}{\operatorname{PT}}
\newcommand{\res}{\operatorname{res}}
\newcommand{\Res}{\operatorname{Res}}
\newcommand{\sconnect}{\operatorname{sc}}
\newcommand{\SL}{\operatorname{SL}}
\newcommand{\Spec}{\operatorname{Spec}}
\newcommand{\Sp}{\operatorname{Sp}}
\newcommand{\ssimple}{\operatorname{ss}}
\newcommand{\tors}{\operatorname{tors}}
\newcommand{\type}{\mathbf{type}}
\newcommand{\umult}{\mathbf{umult}}
\newcommand{\Zar}{\operatorname{Zar}}

\DeclareSymbolFont{cyrletters}{OT2}{wncyr}{m}{n}
\DeclareMathSymbol{\Sha}{\mathalpha}{cyrletters}{"58}

\DeclareMathOperator{\iExt}{\mathscr{E}\text{\kern -3pt {\calligra\large xt}\kern -.7pt}\ }
\DeclareMathOperator{\iHom}{\mathscr{H}\text{\kern -3pt {\calligra\large om}\kern -1pt}\ }

\title{Arithmetics of homogeneous spaces over $p$-adic function fields}
\author{Nguy$\tilde{\hat{\text{e}}}$n M$\d{\text{a}}$nh Linh}
\date{\today}

\begin{document}
\maketitle
\begin{abstract}
    Let $K$ be the function field of a smooth projective geometrically integral curve over a finite extension of $\mathbb{Q}_p$. Following the works of Harari, Scheiderer, Szamuely, Izquierdo, and Tian, we study the local-global and weak approximation problems for homogeneous spaces of $\textrm{SL}_{n,K}$ with geometric stabilizers extension of a group of multiplicative type by a unipotent group. The tools used are arithmetic (local and global) duality theorems in Galois cohomology, in combination with techniques similar to those used by Harari, Szamuely, Colliot-Th\'el\`ene, Sansuc, and Skorobogatov. As a consequence, we show that any finite abelian group is a Galois group over $K$, rediscovering the positive answer to the abelian case of the inverse Galois problem over $\mathbb{Q}_p(t)$. In the case where the curve is defined over a higher-dimensional local field instead of a finite extension of $\mathbb{Q}_p$, coarser results are also given.
\end{abstract}

\tableofcontents
\bigskip
	
\section{Introduction} \label{sec:Intro}
\subsection{Context of the problem} \label{subsec:Context}

Let $k$ be a field and let $\Omega$ be a set of places ({\em i.e.} equivalence classes of non-trivial absolute values) of $k$. For each $v \in \Omega$, we denote by $k_v$ the corresponding completion. We say that a class of smooth $k$-varieties {\em satisfy the local-global principle} (with respect to $\Omega$) if for every variety $X$ in this class, $\prod_{v \in \Omega} X(k_v) \neq \varnothing$ implies $X(k) \neq \varnothing$. 

Suppose that $X$ is a $k$-variety with $X(k) \neq \varnothing$. For each $v \in \Omega$, the set $X(k_v)$ is equipped with a local topology induced by the topology of $k_v$ (see for example \cite[Proposition 3.1]{Conrad}). We say that $X$ {\em has weak approximation in $\Omega$} if the diagonal embedding $X(k) \hookrightarrow \prod_{v \in \Omega} X(k_v)$ has dense image (where the product space is equipped with the product of local topologies, hence it is enough to work with finite subsets $S \subseteq \Omega$). For example, affine spaces have weak approximation (Artin--Whaples). 

These arithmetic properties (local-global principle and weak approximation) are $k$-stable birational invariants of smooth complete $k$-varieties, thanks to Serre's generalized version of the implicit function theorem \cite[Part II, Chapter III, \S 10.2]{SerreLie} combined with the theorem of Nishimura--Lang \cite{Nishimura}. The classical case is when $k$ is a number field (or its geometric counterpart, a global function field). In this case, Manin \cite{Manin} defined an obstruction to the existence of rational points using the Brauer--Grothendieck group $\Br X:= \H^2_{\et}(X,\Gbb_m)$ and the global reciprocity law (Albert--Brauer--Hasse--Noether). These were later used by Colliot-Th\'el\`ene and Sansuc to define an obstruction to weak approximation, which is now known as the {\em Brauer--Manin obstruction}.

The present article is motivated by the recent interest in studying these two arithmetic problems over fields of cohomological dimension greater than $2$. Of particular interest is the case of a function field $K$ of a (smooth, complete, geometrically integral) variety $V$ over a field $k$, and where the set of places is that of points $v \in V$ of codimension $1$. In this article, we focus on the case where $k$ is a finite extension of $\Qbb_p$ (except in section \ref{sec:Higher}, where $k$ is a {\em higher-dimensional local field}), and where $V$ is a curve (of course, the classical case of global function fields corresponds to the case where $k$ is finite). The corresponding function field $K$ is then a field of cohomological dimension $3$. 

The main tool for our approach will be Galois cohomology. Over a $p$-adic function field $K$ as above, the group $\H^3_{\et}(X,\Qbb/\Zbb(2))$ (and its variants) would play the role of $\Br X$ in the classical case. Using (local and global) duality theorems for tori, the works of Harari, Scheiderer and Szamuely gave us a certain understanding of the local-global principle \cite{HShasse} and weak approximation \cite{HSS} problems for tori. Tian \cite{TianThese,TianWeak} extended their ideas to study weak approximation for connected reductive groups $G$ with the property that the universal cover $G^{\sconnect}$ of its derived subgroup $G^{\ssimple}$ has weak approximation and contains a split maximal torus. Tian's work relies on the powerful machinery of Borovoi called {\em abelianization of non-abelian Galois cohomology} \cite{BorovoiSecond,BorovoiFirst,BorovoiReductive}. Izquierdo also provided some results over function fields of curves over higher-dimensional local fields \cite[Th\'eor\`eme 0.1]{IzquierdoII}. 

It should be noted that if we consider the completions coming from {\em all} the rank $1$ discrete valuations on $K$ (not just those coming from the closed points of the curve), then many other local-global principles have been established. For these results, the readers are invited to consult the recent survey of Wittenberg \cite[\S 4]{Wittenberg}.

In this article, the varieties under consideration will be homogeneous spaces of $\SL_n$ (or more generally, a {\em special}, simply connected semisimple algebraic group that has weak approximation). We shall see that in the heart of the proof of each of the main theorems, there lies the use of a Poitou--Tate style duality theorem and the Poitou--Tate sequence. Aside from them, the techniques used here are the {\em reinterpretation of the Brauer--Manin pairing} by Harari--Szamuely \cite{HSmotive,HShasse} for the existence of rational points, and by Colliot-Th\'el\`ene--Sansuc--Skorobogatov \cite{CTSdescent,SkorobogatovBeyond}
for the weak approximation property.

{\em Acknowledgements.} The author is funded by a ``Contrat doctoral sp\'eficique normalien'' from \'Ecole normale sup\'erieure de Paris. I am grateful to my PhD advisor, David Harari, for his guidance and support. I thank Jean-Louis Colliot-Th\'el\`ene, Yisheng Tian, Haowen Zhang, and Olivier Wittenberg for the relevant discussions. I also thank the Laboratoire de Math\'ematiques d'Orsay, Universit\'e Paris-Saclay for the excellent working condition.

\subsection{Statement of the main results} \label{subsec:Statement}

Here are the main results of the article. 

\begin{customthm} [Theorems \ref{thm:Hasse} and \ref{thm:HasseModified}] \label{customthm:Hasse} 
    Let $K$ be the function field of a smooth projective geometrically integral curve $\Omega$ over a $p$-adic field, $X$ a homogeneous space of $\SL_{n,K}$ whose geometric stabilizers are extensions of a group of multiplicative type by a unipotent group. Then the unramified first obstruction to the local-global principle for $X$ (that is, the obstruction relative to the subgroup of $\H^3(K(X),\Qbb/\Zbb(2))$ consisting of elements unramified over $K$ and whose restriction to $\H^3(K_v(X),\Qbb/\Zbb(2))$ comes from $\H^3(K_v,\Qbb/\Zbb(2))$ for all closed points $v \in \Omega$) is the only one.
\end{customthm}

\begin{customthm} [Theorems \ref{thm:Weak} and \ref{thm:WeakModified}] \label{customthm:Weak}
    Let $K$ be the function field of a smooth projective geometrically integral curve $\Omega$ over a $p$-adic field, $X$ a homogeneous space of $\SL_{n,K}$ whose geometric stabilizers are extensions of a group of multiplicative type by a unipotent group. Assume that $X(K) \neq \varnothing$. Then the reciprocity obstruction to weak approximation for $X$, relative to the subgroup of $\H^3(K(X),\Qbb/\Zbb(2))$ consisting of elements unramified over $K$ and whose restriction to $\H^3(K_v(X),\Qbb/\Zbb(2))$ comes from $\H^3(K_v,\Qbb/\Zbb(2))$ for all but finitely many closed points $v \in \Omega$, is the only one.
\end{customthm}

Using Theorem \ref{customthm:Weak}, we also give a positive answer to the abelian case of the inverse Galois problem over $p$-adic function fields (see Corollary \ref{cor:InverseGalois}).

Actually, for each of the two above theorems, we present two proofs. The first ones (Theorems \ref{thm:Hasse} and \ref{thm:Weak}) rely on the observation that the two arithmetic problems (the local-global principle and weak approximation) for homogeneous spaces of $\SL_n$ is closely related to the question of {\em universal torsors} over their smooth compactifications. This is part of a ``descent theory'' for torsors under tori, which (over number fields) originated from the formidable work of Colliot-Th\'el\`ene and Sansuc \cite[\S 3]{CTSdescent}, and to which Skorobogatov then made some complements \cite[\S 6]{SkorobogatovTorsors}. Hence, in the course of proving Theorems \ref{customthm:Hasse} and \ref{customthm:Weak}, we shall develop this descent theory for varieties over $p$-adic function fields. In particular, the following analogue of \cite[Th\'eor\`eme 3.8.1, Proposition 3.8.7]{CTSdescent} (see also \cite[Corollary 6.1.3]{SkorobogatovTorsors}) is proposed. 

\begin{customthm} [Theorem \ref{thm:Descent}] \label{customthm:Descent}
    Let $K$ be the function field of a smooth projective geometrically integral curve $\Omega$ over a $p$-adic field. Let $X$ be a smooth proper geometrically integral variety over $K$ such that $\Pic X_{\ol{K}}$ is a finitely generated free abelian group. If the universal torsors $Y \to X$ (see its definition at the beginning of section \ref{sec:Descent}) satisfy the local-global principle (resp. the local-global principle and weak approximation), then the reciprocity obstruction (that is, the obstruction relative to the subgroup of $\H^3(K(X),\Qbb/\Zbb(2))$ consisting of elements unramified over $K$) to the local-global principle (resp. the local-global principle and weak approximation) on $X$ is the only one.
\end{customthm}

The second proofs of Theorems \ref{customthm:Hasse} and \ref{customthm:Weak} (Theorems \ref{thm:HasseModified} and \ref{thm:WeakModified}, respectively) are inspired by the following remark, communicated to the author by Jean-Louis Colliot-Th\'el\`ene\footnote{The author is grateful to Colliot-Th\'el\`ene for allowing him to include these results in the present article.}: Every homogeneous space of $\SL_{n,K}$, whose geometric stabilizers are extensions of a group of multiplicative type by a unipotent group, is $K$-stably birational to a torus. This allows us to deduce Theorems \ref{thm:HasseModified} and \ref{thm:WeakModified} from \cite[Theorem 0.2]{HShasse} and \cite[Theorem 1.2]{HSS}, respectively.

For function fields of curves over higher-dimensional local fields, we propose weaker versions of Theorems \ref{customthm:Hasse} and \ref{customthm:Weak}. If $k$ is a $d$-dimensional local field $k$ (see paragraph \ref{subsec:Obstruction} below), consider the following condition.
\begin{equation*} \label{eq:Star}
    \text{\em $k = \Cbb\ps{t}$; or $d \ge 1$ and the $1$-local field associated with $k$ has characteristic $0$.}
    \tag{$\star$} 
\end{equation*}
In particular, $k$ has characteristic $0$. Note that the case where $k$ is $p$-adic corresponds to $d=1$. 

\begin{customthm} [Theorem \ref{thm:HigherHasse}] \label{customthm:HigherHasse} 
     Let $K$ be the function field of a smooth projective geometrically integral curve $\Omega$ over a $d$-dimensional local field $k$ satisfying \eqref{eq:Star}, $X$ a homogeneous space of $\SL_{n,K}$ with finite abelian geometric stabilizers. Then the adelic first obstruction to the local-global principle for $X$ (that is, the obstruction relative to the subgroup of $\H^{d+2}(X,\Qbb/\Zbb(d+1))$ consisting of elements whose restriction to $\H^{d+2}(X_{K_v}, \Qbb/\Zbb(d+1))$ comes from $\H^{d+2}(K_v,\Qbb/\Zbb(d+1))$ for all closed points $v \in \Omega$)
     is the only one.
\end{customthm}

\begin{customthm} [Theorem \ref{thm:HigherWeak}] \label{customthm:HigherWeak} 
     Let $K$ be the function field of a smooth projective geometrically integral curve $\Omega$ over a $d$-dimensional local field $k$ satisfying \eqref{eq:Star}, $X$ a homogeneous space of $\SL_{n,K}$ with finite abelian geometric stabilizers. Assume that $X(K) \neq \varnothing$. Then the generalized Brauer--Manin obstruction to weak approximation for $X$ (that is, the obstruction relative to the subgroup of $\H^{d+2}(X,\Qbb/\Zbb(d+1))$ consisting of elements whose restriction to $\H^{d+2}(X_{K_v}, \Qbb/\Zbb(d+1))$ comes from $\H^{d+2}(K_v,\Qbb/\Zbb(d+1))$ for all but finitely many closed points $v \in \Omega$) is the only one.
\end{customthm}

Theorems \ref{customthm:HigherHasse} and \ref{customthm:HigherWeak} do not hold for stabilizers of multiplicative type; counter-examples in the case where $d = 0$ and where the stabilizers are tori are given in Example \ref{ex:CounterExampleToricStabilizer}.

The condition \eqref{eq:Star} is crucial to establish duality theorems for finite modules. It can be slightly weakened by allowing $k$ (if $d = 0$) or the $1$-dimensional local associated with $k$ (if $d \ge 1)$ to have characteristic $p > 0$, provided that $p$ does not divide the order of the stabilizers.

The article is organized as follows. In paragraphs \ref{subsec:Notation}, \ref{subsec:Obstruction}, and \ref{subsec:Duality}, we recall the notations, conventions, and known results that will be used in the proofs of our main theorems. In section \ref{sec:Descent}, we develop a version of Colliot-Th\'el\`ene--Sansuc descent theory for varieties over $p$-adic function fields and prove Theorem \ref{customthm:Descent}. Theorems \ref{customthm:Hasse} and \ref{customthm:Weak} shall be proved in section \ref{sec:PAdic}, as consequences of the results already established in section \ref{sec:Descent}. Here we also prove that any finite abelian group is a Galois group over any $p$-adic function field. Finally, we prove Theorems \ref{customthm:HigherHasse} and \ref{customthm:HigherWeak} (which are results over function fields over higher-dimensional local fields) in section \ref{sec:Higher}.

\subsection{Notations and conventions} \label{subsec:Notation}

The following conventions shall be deployed throughout the article.

{\em Cohomology.} Unless stated otherwise, all (hyper-)cohomology groups will be {\em \'etale} or {\em Galois}. We abusively identify each object of an abelian category to the corresponding $1$-term complex concentrated in degree $0$. Over a field $K$, a {\em $K$-variety} is a separated scheme of finite type $X \to \Spec K$. We use $\ol{K}$ to denote a fixed {\em separable} closure of $K$, $\Gamma_K = \Gal(\ol{K}/K)$ to denote its absolute Galois group, and $\ol{X} = X \times_K \ol{K}$. We denote by $\Dcal^+(X)$ (resp. $\Dcal^+(K)$, resp. $\Dcal^+(\Ab)$) the bounded-below derived category of \'etale sheaves over $X$ (resp. of discrete $\Gamma_K$-modules\footnote{The category of $\Gamma_K$-modules is equivalent to that of sheaves over $(\Spec K)_{\et}$, we abusively identify these two.}, resp. of abelian groups). If $v$ is a place of $K$, $K_v$ denotes the corresponding completion, $X_v = X \times_K K_v$, and $\loc_v: \H^i(X,-) \to \H^i(X_v,-)$ denotes the localization map in cohomology. Whenever $K$ is the function field of a smooth projective geometrically integral curve $\Omega$ over a field $k$, $\Omega^{(1)}$ denotes the set of closed points of $\Omega$. For each $v \in \Omega^{(1)}$, $K_v$ (resp. $\Ocal_v$) denotes the $v$-adic completion of $K$ (resp. of the local ring $\Ocal_{\Omega,v}$), while $K_v^{\h}$ and $\Ocal_v^{\h}$ denote the corresponding henselizations, and $k(v)$ denotes the residue field of $v$.

{\em Motivic complexes.} Let $X$ be a smooth variety over a field $K$. Lichtenbaum defined arithmetic complexes $\Zbb(i)$ over $X_{\et}$ for $i = 0,1,2$ \cite{LichtenbaumConstruction,LichtenbaumNew} (we shall write $\Zbb_X(i)$ if we want to emphasize the variety $X$). We have quasi-isomorphisms $\Zbb(0) \cong \Zbb$ and $\Zbb(1) \cong \Gbb_m[-1]$. The complex $\Zbb(2)$ is concentrated in degrees $1$ and $2$. Also, there is a pairing
    \begin{equation} \label{eq:PairingMotivic}
        \Zbb(1) \otimes^{\Lbb} \Zbb(1) \to \Zbb(2).
    \end{equation}
If $\Fcal$ is a sheaf on $X_{\et}$, let $\Fcal(i) = \Fcal \otimes^{\Lbb} \Zbb(i)$. For $n$ invertible in $K$, there are quasi-isomorphisms $\Zbb/n(i) \cong \mu_n^{\otimes i}$, $i = 1,2$. It follows that $\Qbb/\Zbb(i) \cong \varinjlim_{n} \mu_n^{\otimes i}$ if $K$ has characteristic $0$ and $i \in \{1,2\}$. We also use this as the definition of the sheaf $\Qbb/\Zbb(i)$ for $i \notin \{1,2\}$.

{\em Abelian groups.} For a topological abelian group $A$ (the topology is understood to be discrete if not specified), $A^D = \Hom_{\cts}(A,\Qbb/\Zbb)$ denotes its Pontrjagin dual. For $n \ge 1$, we denote by $\tensor[_n]{A}{}$ the $n$-torsion subgroup of $A$, and $A_{\tors} = \varinjlim_n \tensor[_n]{A}{}$.

{\em Tori.} If $G$ is a smooth algebraic group over a field $K$, we denote by $\wh{G} = \iHom_K(G,\Gbb_m)$ its $\Gamma_K$-module of {\em geometric characters}, that is, $\wh{G} = \Hom_{\ol{K}}(\ol{G},\Gbb_m)$ equipped with the Galois action defined by the formula
    \begin{equation*}
        \forall \sigma \in \Gamma_K,\forall \chi \in \Hom_{\ol{K}}(\ol{G},\Gbb_m),\forall g \in G(\ol{K}), \quad (\tensor[^\sigma]{\chi}{})(g):= \tensor[^\sigma]{}{}(\chi(\tensor[^{\sigma^{-1}}]{g}{})).
    \end{equation*}
If $T$ is a torus over a field $K$, its {\em dual torus} is defined to be the torus $T'$ whose character module is the {\em cocharacter module} $\wc{T} = \iHom_K(\Gbb_m,T) = \Hom(\wh{T},\Zbb)$ of $T$. Since $T = \wc{T} \otimes \Gbb_m = \wc{T} \otimes^{\Lbb} \Gbb_m$ and $T' = \wh{T} \otimes \Gbb_m = \wh{T} \otimes^{\Lbb} \Gbb_m$, we obtain a pairing 
    \begin{equation*}
        T \otimes^{\Lbb} T' \to \Gbb_m \otimes^{\Lbb} \Gbb_m \cong \Zbb(1)[1] \otimes^{\Lbb} \Zbb(1)[1].
    \end{equation*}
Then \eqref{eq:PairingMotivic} induces a pairing
    \begin{equation}   \label{eq:PairingTorus}
        T \otimes^{\Lbb} T' \to \Zbb(2)[2]
    \end{equation}
in $\Dcal^+(K)$. We say that $T$ is {\em quasi-split} if it is isomorphic to $\Res_{A/K} \Gbb_{m,A}$ for some \'etale $K$-algebra $A$, where $\Res_{A/K}$ denotes the restriction of scalars {\em \`a la} Weil. This is equivalent to saying that $\wh{T}$ is a permutation module ({\em i.e.} it has a $\Gamma_K$-invariant $\Zbb$-basis). In this case, $\H^1(K, T) = 0$ by Shapiro's lemma and Hilbert's Theorem 90. Also, $T'$ is quasi-split.

{\em Tate--Shafarevich groups.} Let $K$ be the function field of a smooth projective geometrically integral curve $\Omega$ over a field $k$. Let $S \subseteq \Omega^{(1)}$ be a {\em finite} set of closed points, and $C$ a complex of $\Gamma_K$-modules. For $i \in \Zbb$, we define the groups
    \begin{align*}
        \Sha^i_S(K,C) = \Ker\left(\H^i(K,C) \to \prod_{v \notin S} \H^i(K_v,C)\right),\\
        \Sha^i(K,C) = \Sha^i_{\varnothing}(K,C), \quad \Sha^i_\omega(K,C) = \varinjlim_S \Sha^i_S(K,C).
    \end{align*}

{\em Unramified cohomology.} Let $X$ a smooth integral variety over a field $K$, $n \ge 1$ an integer invertible in $K$, $i \ge 0$, and $j \in \Zbb$. Let $\mu_n^{\otimes j} := \iHom_K(\mu_n^{\otimes (-j)},\Zbb/n)$ if $j < 0$ and $\mu_n^{\otimes 0} := \Zbb/n$. One defines the {\em unramified} part $\H^i_{\nr}(K(X)/K,\mu_n^{\otimes j})$ to be the subgroup of $\H^i(K(X),\mu_n^{\otimes j})$ consisting of elements $A$ that lift to $\H^i(\Ocal,\mu_n^{\otimes j})$ for every discrete valuation ring $\Ocal \supseteq K$ with field of fractions $K(X)$ (such a lifting is necessarily unique by the injectivity property, see \cite[Theorem 3.8.1]{CTgersten}). If $X$ is proper, this amounts to requiring that $A$ comes from $\H^i(\Ocal_{X,v},\mu_n^{\otimes j})$ for every point $v \in X$ of codimension 1 (see Theorem 4.1.1 in {\em loc. cit.}). The group $\H^i_{\nr}(K(X)/K,\mu_n^{\otimes j})$ is  a $K$-stable birational invariant \cite[Proposition 4.1.4]{CTgersten}). We define the ``evaluation'' pairing
    \begin{equation} \label{eq:EvaluationPairing}
        \H^i_{\nr}(K(X)/K, \mu_n^{\otimes j}) \times X(K') \to \H^i(K',\mu_n^{\otimes j}),
    \end{equation}
for any overfield $K'/K$, as follows. By Bloch--Ogus theorem (Gersten's conjecture for \'etale cohomology) \cite{BlochOgus} (see also \cite[Theorem 4.1.1]{CTgersten}), every class $A \in \H^i_{\nr}(K(X)/K, \mu_n^{\otimes j})$ comes from $\H^i(\Ocal_{X,P}, \mu_n^{\otimes j})$ for {\em any} point $P \in X$. Let $P': \Spec K' \to X$ be a $K'$-point with image $P \in X$. Lift $A$ to a unique element of $\H^i(\Ocal_{X,P}, \mu_n^{\otimes j})$, and then define $A(P')$ to be its image by the pullback $\H^i(\Ocal_{X,P}, \mu_n^{\otimes j}) \to \H^i(K', \mu_n^{\otimes j})$. 

Assume in addition that $K$ has characteristic $0$ and that $X$ is proper. For $i \ge 3$, there is a natural map (for the case where $i = 3$, see \cite[Proposition 2.9]{Kahn} and \cite[(19)]{HSS})
    \begin{equation} \label{eq:GerstenResolution}
        \H^{i+1}(X,\Zbb(2)) \to \H^i_{\nr}(K(X)/K,\Qbb/\Zbb(2)),
    \end{equation}
defined as follows. First, let $\Zbb(2)_{\Zar}$ be the complex concentrated in degrees $\le 2$ defined in a similar way to $\Zbb(2)$, but over the small {\em Zariski} site $X_{\Zar}$. Let $\alpha: X_{\et} \to X_{\Zar}$ denote the change-of-sites map. Then the adjunction $\Qbb(2)_{\Zar} \to \Rbb \alpha_\ast \Qbb(2)$ is an isomorphism \cite[Th\'eor\`eme 2.6(c)]{Kahn}. Since $\Qbb(2)_{\Zar}$ is concentrated in degrees $\le 2$, one has $\Rbb^i \alpha_\ast \Qbb(2) = \Rbb^{i+1} \alpha_\ast \Qbb(2) = 0$, hence $\Rbb^{i+1}\alpha_\ast \Zbb(2) \cong \Rbb^i \alpha_\ast \Qbb/\Zbb(2)$. The Leray spectral sequence for $\alpha$ yields an edge map
    \begin{equation*}
        \H^{i+1}(X,\Zbb(2)) \to \H^0_{\Zar}(X, \Rbb^{i+1}\alpha_\ast \Zbb(2)) \cong \H^0_{\Zar}(X, \Rbb^{i}\alpha_\ast \Qbb/\Zbb(2)).
    \end{equation*}
Since $X$ is proper, we have $\H^0_{\Zar}(X, \Rbb^{i}\alpha_\ast \Qbb/\Zbb(2)) \cong \H^i_{\nr}(K(X)/K,\Qbb/\Zbb(2))$ by the Gersten resolution \cite[Theorem 4.1.1]{CTgersten}. The map \eqref{eq:GerstenResolution} is defined. For $X = \Spec K$, this yields an isomorphism $\H^{i+1}(K,\Zbb(2)) \xrightarrow{\cong} \H^i(K,\Qbb/\Zbb(2))$. According to its construction, the map \eqref{eq:GerstenResolution} is functorial in $X$ and in $K$. Applying this functoriality to the natural morphism $\Spec (K(X)) \to X$ yields a commutative diagram

\begin{equation} \label{eq:GerstenResolutionDiagram}
    \xymatrix{
        \H^{i+1}(X,\Zbb(2)) \ar[r] \ar[d] & \H^i_{\nr}(K(X)/K,\Qbb/\Zbb(2)) \ar@{_{(}->}[d] \\
        \H^{i+1}(K(X),\Zbb(2)) \ar[r]^-{\cong} & \H^i(K(X),\Qbb/\Zbb(2)).
    }
\end{equation}

{\em Algebraic groups.} Let $K$ be a field, $X$ a smooth $K$-variety and $G$ a smooth $K$-group scheme. The (non-abelian) \'etale cohomology (pointed) set $\H^1(X,G)$ classifies $X$-torsors under $G$. Let $f: Y \to X$ be such a torsor, with class $[Y] \in \H^1(X,G)$. For each Galois cocycle $c: \Gamma_K \to G(\ol{K})$, we may {\em twist} $G$ to obtain an inner $K$-form $\tensor[_c]{}{} G$ of $G$, and a torsor $\tensor[_c]{f}{}: \tensor[_c]{}{} Y \to X$ under $\tensor[_c]{}{} G$. For each $P \in X(K)$, we have the following equivalence:
\begin{equation} \label{eq:Twisting}
    [Y](P) = [c] \in \H^1(K,G) \Leftrightarrow [\tensor[_c]{}{} Y](P) = 1 \in \H^1(K,\tensor[_c]{}{}G) \Leftrightarrow P \in \tensor[_c]{f}{}(\tensor[_c]{}{} Y(K))
\end{equation}
(see \cite[p. 22]{SkorobogatovTorsors}). If $G$ is commutative, then $\tensor[_c]{}{} G = G$ and $[\tensor[_c]{}{}Y] = [Y] - \pi^\ast[c] \in \H^1(X,G)$, where $\pi: X \to \Spec K$ denotes the structure morphism. The following continuity result is perhaps well-known. We present a proof here for the lack of a reference.

\begin{lm} \label{lm:Continuity}
     Let $K$ be a topological field ({\em i.e.} its addition, multiplication and inversion are continuous), for example, a field equipped with an absolute value. Let $X$ be a smooth integral variety over $K$ and consider the natural topology on $X(K)$ \cite[Proposition 3.1]{Conrad}.
     \begin{enumerate}
        \item \label{lm:Continuity1} Let $G$ be a smooth, quasi-projective, commutative group scheme over $K$. Then for any $i \ge 0$ and $A \in \H^i(X,G)$, the induced evaluation $X(K) \to \H^i(K,G)$ is locally constant.
        
        \item \label{lm:Continuity2} For any $n \ge 1$, $i \ge 0$, $j \in \Zbb$, and $A \in \H^i_{\nr}(K(X)/K,\mu_n^{\otimes j})$, the induced evaluation $X(K) \to \H^i(K,\mu_n^{\otimes j})$ is locally constant.

        \item \label{lm:Continuity3} Suppose that $K$ is complete with respect to an absolute value. Let $G$ be a (not necessarily commutative) smooth $K$-group scheme. Then for any torsor $Y \to X$ under $G$ with class $[Y] \in \H^1(X,G)$, the induced evaluation $X(K) \to \H^1(K, G)$ is locally constant.
     \end{enumerate}
\end{lm}
\begin{proof}
    We prove \ref{lm:Continuity1}, the proof of \ref{lm:Continuity2} being similar. Denote by $\pi: X \to \Spec K$ the structure morphism. Fix a $K$-point $P_0 \in X(K)$ (if there is any). Replacing $A$ by $A - \pi^\ast A(P_0)$, we may assume that $A(P_0) = 0$. Since $\H^i(K,G) \cong \H^i(\Ocal_{X,{P_0}}^{\h}, G)$ \cite[Chapter III, Remark 3.11]{MilneEtale}, there is an \'etale neighborhood $f: X' \to X$ of $P_0$ and a point $P'_0 \in X'(K)$ such that $f(P'_0) = P_0$ and $f^\ast A = 0$. Now, $f$ induces a homeomorphism between a neighborhood (for the $K$-topology) $\Uscr' \subseteq X'(K)$ of $P_0'$ and a neighborhood $\Uscr \subseteq X(K)$ of $P_0$. Each point $P \in \Uscr(K)$ then has a factorization $\Spec K \to X' \xrightarrow{f} X$, thus the evaluation-at-$P$ factorizes through the pullback map $f^\ast: \H^i(X,G) \to \H^i(X', G)$. Since $f^\ast A = 0$, it follows that $A(P) = 0$.

    Let us now show \ref{lm:Continuity3}. Fix a $K$-point $P_0 \in X(K)$ and let $c: \Gamma_K \to G(\ol{K})$ be a Galois cocycle representing $[Y](P_0) \in \H^1(K,G)$. Twisting by $c$ yields a torsor $\tensor[_c]{}{} Y$ under the inner form $\tensor[_c]{}{}G$, which contains a $K$-point lying over $P_0$ (in view of \eqref{eq:Twisting}). Since $K$ is complete with respect to an absolute value, we may apply Serre's generalized version of the implicit function theorem \cite[Part II, Chapter III, \S 10.2]{SerreLie}. This assures the existence of a neighborhood $\Uscr \subseteq X(K)$ of $P_0$, whose $K$-points $P$ has can be lifted to $K$-points of $\tensor[_{c}]{}{} Y$. In other words, $[Y](P) = [c] = [Y](P_0)$ (again, by \eqref{eq:Twisting}) for all $P \in \Uscr$.
\end{proof}

Following \cite[\S 4.2]{CTSrational}, we shall say that $G$ is {\em special} if any torsor under $G$ over a $K$-variety is locally trivial for the Zariski topology. This implies $\H^1(L,G) = 1$ for any overfield $L/K$. It follows that if $H \subseteq G$ is a Zariski closed subgroup and $X = H \backslash G$ (hence the projection $G \to X$ is a torsor under $H$), then the evaluation map
    \begin{equation*}
        X(L) \to \H^1(L,H), \quad P \mapsto [G](P)
    \end{equation*}
induces a bijection $X(L)/G(L) \cong \H^1(L,H)$. Serre proved that any special group is linear and connected \cite[\S 4.1, Th\'eor\`eme 1]{SerreChow}. Examples of special groups are the general linear group $\GL_n$ (or more generally, $\GL_A$ for a central simple algebra $A$ over $K$), the special linear group $\SL_n$, or the symplectic group $\Sp_{2n}$.

\subsection{Reciprocity obstruction} \label{subsec:Obstruction}

The generalized Weil reciprocity law is the key to defining the reciprocity obstructions to the local global principle and weak approximation for varieties over fields with arithmetico-geometric nature. In this paper, we work with function fields of curves over $p$-adic fields, and more generally, higher-dimensional local fields. We shall recall their definition and properties here. By definition, a {\em $0$-dimensional local field} is either a finite field or the field $\Cbb\ps{t}$ for some algebraically closed field $\Cbb$ of characteristic $0$. These fields have absolute Galois group $\wh{\Zbb}$ (for $\Cbb\ps{t}$, this is the celebrated Puiseux's theorem), hence have cohomological dimension $1$. For $d \ge 1$, a {\em $d$-dimensional local field} is a complete discretely valued field whose residue field is a $(d-1)$-dimensional local field.

\begin{prop} \label{prop:TraceIsomorphism}
    Let $k$ be a $d$-local field satisfying the condition \eqref{eq:Star} from page \pageref{eq:Star}.
    \begin{enumerate}
        \item \label{prop:TraceIsomorphism1} The field $k$ has cohomological dimension $\cd(k) = d+1$.
        
        \item \label{prop:TraceIsomorphism2} We have a canonical trace isomorphism $\H^{d+1}(k,\Qbb/\Zbb(d)) \cong \Qbb/\Zbb$. More precisely, if $k_{i-1}$ denotes the residue field of $k_i$ for $i \in \{1,\ldots,d\}$ and $k_d = k$, then this isomorphism is given by the composite of the ``residue'' maps
        \begin{equation} \label{eq:ResidueIsomorphism}
            \H^{d+1}(k,\mu_n^{\otimes d}) \overset{\partial}{\cong} \H^d(k_{d-1},\mu_n^{\otimes (d-1)}) \overset{\partial}{\cong} \cdots \overset{\partial}{\cong} \H^1(k_0,\Zbb/n) \cong \Zbb/n.
        \end{equation}
        
        \item \label{prop:TraceIsomorphism3} For any finite extension $\ell/k$, the trace isomorphisms fit into a commutative diagram
             \begin{equation} \label{eq:Restriction}
                \xymatrix{
                    \H^{d+1}(k,\Qbb/\Zbb(d)) \ar[r]^-{\cong} \ar[d]^{\res_{\ell/k}} & \Qbb/\Zbb \ar[d]^{[\ell : k] \cdot} \\
                    \H^{d+1}(\ell,\Qbb/\Zbb(d)) \ar[r]^-{\cong} & \Qbb/\Zbb,
                }
            \end{equation}
         where $\res_{\ell/k}: \H^{d+1}(k,\Qbb/\Zbb(d)) \to \H^{d+1}(\ell,\Qbb/\Zbb(d))$ denotes the restriction map. Consequently, we have a commutative diagram
        \begin{equation}   \label{eq:Corestriction}
            \xymatrix{
                \H^{d+1}(\ell,\Qbb/\Zbb(d)) \ar[r]^-{\cong} \ar[d]^{\cores_{\ell/k}} & \Qbb/\Zbb \ar@{=}[d] \\
                \H^{d+1}(k,\Qbb/\Zbb(d)) \ar[r]^-{\cong} & \Qbb/\Zbb,
                }
        \end{equation}
        where $\cores_{\ell/k}: \H^{d+1}(\ell,\Qbb/\Zbb(d)) \to \H^{d+1}(k,\Qbb/\Zbb(d))$ denotes the corestriction map.
    \end{enumerate}
\end{prop}
\begin{proof}
    Note that the commutativity of \eqref{eq:Corestriction} follows from that of \eqref{eq:Restriction}, since the composite $\cores_{\ell/k} \circ \res_{\ell/k}: \H^{d+1}(k,\Qbb/\Zbb(d)) \to \H^{d+1}(k,\Qbb/\Zbb(d))$ is the multiplication by $[\ell:k]$, and since $\Qbb/\Zbb$ is divisible. We prove \ref{prop:TraceIsomorphism1}, \ref{prop:TraceIsomorphism2}, and \ref{prop:TraceIsomorphism3} by induction on $d$. For $d = 0$, one has $\Gamma_k \cong \wh{\Zbb}$, so \ref{prop:TraceIsomorphism1} and \ref{prop:TraceIsomorphism2} are obvious. Furthermore, for any finite extension $\ell/k$, the inclusion $\Gamma_{\ell} \hookrightarrow \Gamma_k$ is the multiplication by $[\ell:k]$ on $\wh{\Zbb}$, which implies \ref{prop:TraceIsomorphism3}.

    Suppose that $d \ge 1$, and let $\kappa = k_{d-1}$ denote the residue field of $k$ (which is a $(d-1)$-local field). We recall from \cite[Chapitre II, Annexe, (2.2)]{SerreGalois} that for any $n \ge 1$, $i \ge 0$ and $j \in \Zbb$, the Hochschild--Serre spectral sequence yields an exact sequence
        \begin{equation} \label{eq:Residue}
            0 \to \H^{i+1}(\kappa,\mu_n^{\otimes j}) \to \H^{i+1}(k,\mu_n^{\otimes j}) \xrightarrow{\partial} \H^i(\kappa,\mu_n^{\otimes (j-1)}) \to 0
        \end{equation}
    where $\partial$ is the {\em residue} map. Since $\cd(\kappa) = d$, setting $i=j=d$ in \eqref{eq:Residue} yields
        \begin{equation*}
            \H^{d+1}(k,\mu_n^{\otimes d}) \overset{\partial}{\cong} \H^d(\kappa,\mu_n^{\otimes (d-1)}) \cong \Zbb/n
        \end{equation*}
    Furthermore, one has $\cd(k) = d+1$ by \cite[Chapitre II, \S 4.3, Proposition 12]{SerreGalois}. It remains to establish diagram \eqref{eq:Restriction}. For this, we make use of the following standard property of $\partial$. If $\ell$ is a finite extension of $k$ with residue field $\lambda$, then we have a commutative diagram
        \begin{equation*}
            \xymatrix{
                \H^{i+1}(k,\mu_n^{\otimes j}) \ar[r]^-{\partial} \ar[d]^{\res_{\ell / k}} & \H^i(\kappa,\mu_n^{\otimes (j-1)}) \ar[d]^{e \cdot \res_{\lambda/\kappa}} \\
                \H^{i+1}(\ell,\mu_n^{\otimes j}) \ar[r]^-{\partial} & \H^i(\lambda,\mu_n^{\otimes (j-1)}),
            }
        \end{equation*}
    where $e$ is the ramification index (this follows easily from the description of the residue map in \cite[\S 1, (1.3)(i)]{Kato}). Since $[\ell:k] = e \cdot [\lambda:\kappa]$, we deduce the commutativity of \eqref{eq:Restriction} from that of the same diagram for the extension $\lambda/\kappa$.
\end{proof}

From now on, let $k$ be a $d$-dimensional local field satisfying the condition \eqref{eq:Star} from page \pageref{eq:Star} and let $K$ be the function field of a smooth projective geometrically integral curve $\Omega$ over $k$. For each closed point $v \in \Omega^{(1)}$, the $v$-adic completion $K_v$ is a $(d+1)$-dimensional local field, a field of cohomological dimension $\cd(K_v) = d+2$. The field $K$ itself has cohomological dimension $\cd(K) \le d+2$. Indeed, by \cite[Proposition 5.10]{HarariCohomologie}, it is enough to show that $\cd(k(t)) \le d+2$. Let $A$ be any torsion $\Gamma_{k(t)}$-module. Since $\Gamma_k \cong \Gal(\ol{k}(t)/k(t))$, we have a spectral sequence  $\H^p(k,\H^q(\ol{k}(t),A)) \Rightarrow \H^{p+q}(k(t),A)$. But $\cd(\ol{k}(t)) \le 1$ by Tsen's theorem \cite[Proposition 6.2.3, Theorem 6.2.8]{GS} and $\cd(k) = d+1$, hence $\H^i(k(t),A) = 0$ for $i > d+2$, or $\cd(k(t)) \le d+2$.

\begin{prop} [Generalized Weil reciprocity law] \label{prop:WeilReciprocity}
    We have a complex
        \begin{equation} \label{eq:WeilReciprocity}
        \H^{d+2}(K,\Qbb/\Zbb(d+1)) \xrightarrow{(\loc_v)_{v \in \Omega^{(1)}}} \bigoplus_{v \in \Omega^{(1)}} \H^{d+2}(K_v,\Qbb/\Zbb(d+1)) \xrightarrow{\sigma} \Qbb/\Zbb,
    \end{equation}
    where $\sigma$ is the sum of the isomorphisms $\H^{d+2}(K_v,\Qbb/\Zbb(d+1)) \cong \Qbb/\Zbb$ from Proposition \ref{prop:TraceIsomorphism}\ref{prop:TraceIsomorphism2}.
 \end{prop}
 \begin{proof}
     We prove this by applying the reciprocity law for {\em cycle modules}. The family of Galois cohomology groups $\{\H^i(-,\mu_n^{\otimes j})\}_{i,j}$ form a cycle module, that is, a {\em cycle premodule} in the sense of \cite[Definition 1.1]{Rost} satisfying the axioms \textbf{(FD)} and \textbf{(C)} of Definition 2.1 in {\em loc. cit.} For the proofs, we refer to \cite[Remarks 1.11 and 2.5]{Rost} and \cite[Proposition 1.7]{Kato}. For each $n \ge 1$, applying the property \textbf{(RC)} of \cite[Proposition 2.2]{Rost}, we obtain a complex
    \begin{equation*}
        \H^{d+2}(K,\mu_n^{\otimes (d+1)}) \xrightarrow{(\partial_v)_{v \in \Omega^{(1)}}} \bigoplus_{v \in \Omega^{(1)}} \H^{d+1}(k(v),\mu_n^{\otimes d}) \xrightarrow{\sum_{v \in \Omega^{(1)}} \cores_{k(v)/k}} \H^{d+1}(k,\mu_n^{\otimes d}),
    \end{equation*}
    where $\partial_v$ is the composite of the residue map $\H^{d+2}(K_v, \mu_n^{\otimes (d+1)}) \xrightarrow{\cong} \H^{d+1}(k(v), \mu_n^{\otimes d}) \cong \Zbb/n$ from \eqref{eq:ResidueIsomorphism} and $\loc_v$. Taking \eqref{eq:Corestriction} into account, one obtains \eqref{eq:WeilReciprocity}.
\end{proof}

Let $X$ be a smooth geometrically integral $K$-variety such that $\prod_{v \in \Omega^{(1)}} X(K_v) \neq \varnothing$. The evaluation pairings defined in \eqref{eq:EvaluationPairing} for the extensions $K_v/K$ reassemble into a pairing
    \begin{equation*}
        \H^{d+2}_{\nr}(K(X)/K,\Qbb/\Zbb(d+1)) \times \prod_{v \in \Omega^{(1)}} X(K_v) \to \Qbb/\Zbb, \quad (A,(P_v)_{v \in \Omega^{(1)}}) \mapsto \sum_{v \in \Omega^{(1)}} A(P_v)
    \end{equation*}
(the above sum is finite by \cite[Proposition 2.5(i)]{CTPS}, see also \cite[Lemma 5.1]{HSS} for the case where $k$ is $p$-adic). Thanks to \eqref{eq:WeilReciprocity}, the above pairing vanishes on the image of $\H^{d+2}(K,\Qbb/\Zbb(d+1))$ in $\H^{d+2}_{\nr}(K(X)/K,\Qbb/\Zbb(d+1))$. Hence it induces a pairing
    \begin{equation} \label{eq:UnramifiedPairing}
        \tfrac{\H^{d+2}_{\nr}(K(X)/K,\Qbb/\Zbb(d+1))}{\Img\H^{d+2}(K,\Qbb/\Zbb(d+1))} \times \prod_{v \in \Omega^{(1)}} X(K_v) \to \Qbb/\Zbb,
    \end{equation}
which vanishes on the diagonal image of $X(K)$ on $\prod_{v \in \Omega^{(1)}} X(K_v)$ (by \eqref{eq:WeilReciprocity} again). The absence of a family $(P_v)_{v \in \Omega^{(1)}} \prod_{v \in \Omega^{(1)}} X(K_v)$ orthogonal to $\H^{d+2}_{\nr}(K(X)/K,\Qbb/\Zbb(d+1))$ is an obstruction to the existence of $K$-rational points on $X$. We refer to this as the {\em reciprocity obstruction to the local-global principle for $X$}. We are also interested in the following coarser obstruction. Restricting \eqref{eq:UnramifiedPairing} to locally constant classes yields a map
    \begin{equation} \label{eq:UnramifiedReciprocity}
        \rho_X: \Ker\left(\tfrac{\H^{d+2}_{\nr}(K(X)/K,\Qbb/\Zbb(d+1))}{\Img\H^{d+2}(K,\Qbb/\Zbb(d+1))} \to \prod_{v \in \Omega^{(1)}}  \tfrac{\H^{d+2}_{\nr}(K_v(X)/K_v,\Qbb/\Zbb(d+1))}{\Img\H^{d+2}(K_v,\Qbb/\Zbb(d+1))} \right) \to \Qbb/\Zbb.
    \end{equation}
If $X(K) \neq \varnothing$, then $\rho_X = 0$. We refer to the non-vanishing of $\rho_X$ as the {\em (unramified) first obstruction to the local-global principle for $X$}. 

Suppose that $X(K) \neq \varnothing$. By Lemma \ref{lm:Continuity}\ref{lm:Continuity2}, the pairing \eqref{eq:UnramifiedPairing} vanishes on the closure (for the product of $v$-adic topologies) of $X(K)$ in $\prod_{v \in \Omega^{(1)}} X(K_v)$. The non-orthogonality to $\tfrac{\H^{d+2}_{\nr}(K(X)/K,\Qbb/\Zbb(d+1))}{\Img\H^{d+2}(K,\Qbb/\Zbb(d+1))}$ of a family $(P_v)_{v \in \Omega^{(1)}} \in \prod_{v \in \Omega^{(1)}} X(K_v)$ is then an obstruction to approximating this family by $K$-rational points. We refer to this as the {\em reciprocity obstruction to weak approximation for $X$}. We are also interested in the following coarser obstruction to weak approximation. For any finite set $S \subseteq \Omega^{(1)}$, \eqref{eq:UnramifiedPairing} restricts to a pairing
    \begin{equation} \label{eq:UnramifiedPairingS}
        (-,-)_S: \Ker\left(\tfrac{\H^{d+2}_{\nr}(K(X)/K,\Qbb/\Zbb(d+1))}{\Img\H^{d+2}(K,\Qbb/\Zbb(d+1))} \to \prod_{v \notin S}  \tfrac{\H^{d+2}_{\nr}(K_v(X)/K_v,\Qbb/\Zbb(d+1))}{\Img\H^{d+2}(K_v,\Qbb/\Zbb(d+1))} \right) \times \prod_{v \in S} X(K_v) \to \Qbb/\Zbb,
    \end{equation}
which vanishes on the closure of $X(K)$. The non-orthogonality of a family $(P_v)_{v \in S}$ to the subgroup of ``constant-outside-$S$'' elements is then an obstruction to approximating this family by $K$-rational points.

We also note that in the case where $d=1$ (for example, when $k$ is $p$-adic), by the Gersten resolution \eqref{eq:GerstenResolution}, the pairing \eqref{eq:UnramifiedPairing} induces a pairing 
    \begin{equation} \label{eq:UnramifiedPairingGersten}
        \tfrac{\H^4(X,\Zbb(2))}{\Img\H^4(K,\Zbb(2))} \times \prod_{v \in \Omega^{(1)}}X(K_v) \to \Qbb/\Zbb.
    \end{equation}
Similarly, the map \eqref{eq:UnramifiedReciprocity} induces a map
    \begin{equation} \label{eq:UnramifiedReciprocityGersten}
        \rho_X: \Ker\left(\tfrac{\H^4(X,\Zbb(2))}{\Img\H^4(K,\Zbb(2))} \to \prod_{v \in \Omega^{(1)}}  \tfrac{\H^4(X_v,\Zbb(2))}{\Img\H^4(K_v,\Zbb(2))} \right) \to \Qbb/\Zbb.
    \end{equation}
    
Actually, the present article only deals with the (unramified) reciprocity obstruction in the case where $k$ is $p$-adic. For function fields of curves over higher-dimensional fields, we do not prove that the obstruction is unramified. Instead, we shall work with the following adapted version of the pairing \eqref{eq:UnramifiedPairing}. Again, let $K$ be the function field of a smooth projective geometrically integral curve $\Omega$ over a $d$-dimensional local field $k$ satisfying the condition \eqref{eq:Star} from page \pageref{eq:Star}. Let $X$ be a smooth geometrically integral $K$-variety. For a non-empty open subset $U \subseteq \Omega$, we shall call a smooth, separated, finitely presented scheme $\Xcal \to U$ such that $\Xcal \times_U K = X$ an {\em integral model of $X$ over $U$}. Such an integral model exists when $U$ is sufficiently small. For $v \in U^{(1)}$, we have $\Xcal(\Ocal_v) \subseteq X(K_v)$ by the valuative criterion for separatedness. We define the set $X(\Abb_K)$ of {\em adelic points} of $X$ as the subset of families $(P_v)_{v \in \Omega^{(1)}} \in \prod_{v \in \Omega^{(1)}} X(K_v)$ such that $P_v \in \Xcal(\Ocal_v)$ for all but finitely many $v \in U^{(1)}$. If $X(\Abb_K) \neq \varnothing$, we consider the pairing
    \begin{equation*}
        \H^{d+2}(X,\Qbb/\Zbb(d+1)) \times X(\Abb_K) \to \Qbb/\Zbb, \qquad (A,(P_v)) \mapsto \sum_{v \in \Omega^{(1)}} A(P_v),  
    \end{equation*}
which is well-defined. Indeed, $A$ comes from $\H^{d+2}(\Xcal_V, \Qbb/\Zbb(d+1))$ for some non-empty open subset $V \subseteq U$ (where $\Xcal_V = \Xcal \times_U V$), and $P_v$ comes from $\Xcal(\Ocal_v)$ for all but finitely many $v \in V^{(1)}$. It follows that for these $v$, $A(P_v)$ comes from $\H^{d+2}(\Ocal_v,\Qbb/\Zbb(d+1))$. By \cite[Chapter III, Remark 3.11]{MilneEtale}, one has $\H^{d+2}(\Ocal_v,\Qbb/\Zbb(d+1)) \cong \H^{d+2}(k(v),\Qbb/\Zbb(d+1)) = 0$ since $\cd(k(v)) = d+1$. Thus, $A(P_v) = 0$ for all but finitely many $v \in \Omega^{(1)}$ (note that $\Omega \setminus V$ is finite). The above pairing vanishes on the image of $\H^{d+2}(K,\Qbb/\Zbb(d+1)) \to \H^{d+2}(X,\Qbb/\Zbb(d+1))$ by virtue of the generalized Weil reciprocity law \eqref{eq:WeilReciprocity}. Thus, we obtain a pairing
    \begin{equation} \label{eq:AdelicPairing}
        \tfrac{\H^{d+2}(X,\Qbb/\Zbb(d+1))}{\Img \H^{d+2}(K,\Qbb/\Zbb(d+1))} \times X(\Abb_K) \to \Qbb/\Zbb, \qquad (A,(P_v)) \mapsto \sum_{v \in \Omega^{(1)}} A(P_v), 
    \end{equation}
which vanishes on $X(K)$ (again, by \eqref{eq:WeilReciprocity}). Restricting \eqref{eq:AdelicPairing} to locally constant classes yields a map
    \begin{equation} \label{eq:AdelicReciprocity}
        \rho_X: \Ker\left( \tfrac{\H^{d+2}(X,\Qbb/\Zbb(d+1))}{\Img \H^{d+2}(K,\Qbb/\Zbb(d+1))} \to \prod_{v \in \Omega^{(1)}}\tfrac{\H^{d+2}(X_v,\Qbb/\Zbb(d+1))}{\Img\H^{d+2}(K_v,\Qbb/\Zbb(d+1))}\right) \to \Qbb/\Zbb, 
    \end{equation}
which vanishes whenever $X(K) \neq \varnothing$. We refer to the non-vanishing of $\rho_X$ as the {\em (adelic) first obstruction to the local-global principle for $X$}.

Now, suppose that $X(K) \neq \varnothing$. For any finite set $S \subseteq \Omega^{(1)}$, \eqref{eq:AdelicPairing} restricts to a pairing
    \begin{equation} \label{eq:AdelicPairingS}
        (-,-)_S: \Ker\left(\tfrac{\H^{d+2}(X,\Qbb/\Zbb(d+1))}{\Img\H^{d+2}(K,\Qbb/\Zbb(d+1))} \to \prod_{v \notin S}\tfrac{\H^{d+2}(X_v,\Qbb/\Zbb(d+1))}{\Img\H^{d+2}(K_v,\Qbb/\Zbb(d+1))}\right) \times \prod_{v \in S}X(K_v) \to \Qbb/\Zbb,
    \end{equation}
which vanishes on the closure (for the product of $v$-adic topologies) of $X(K)$ in $\prod_{v \in S} X(K_v)$ (this uses Lemma \ref{lm:Continuity}\ref{lm:Continuity1}). The non-vanishing of a family $(P_v)_{v \in S}$ to the group on the left-hand side of \eqref{eq:AdelicPairingS} is an obstruction to approximating this family by $K$-rational points. We refer to this as the {\em generalized Brauer--Manin obstruction to weak approximation for $X$ in $S$}.

\subsection{Arithmetic duality theorems} \label{subsec:Duality}

This paragraph reassembles the tools for proving our main theorems, namely the local and global duality theorems for Galois cohomology of curves over higher-dimensional local fields. 

We start with the case where $k$ is a $d$-dimensional local field satisfying the condition \eqref{eq:Star} from page \pageref{eq:Star}. Recall that for $n \ge 1$, we have an isomorphism $\H^{d+1}(k,\mu_n^{\otimes d}) \cong \Zbb/n$ from \eqref{eq:ResidueIsomorphism}. By \cite[Chapter I, Theorem 2.17]{MilneDuality}, if $F$ is a finite $n$-torsion $\Gamma_k$-module and $F' = \Hom_k(F,\mu_n^{\otimes d})$, the cup-product pairing
    \begin{equation*}
        \H^i(k,F) \times \H^{d+1-i}(k,F') \to \H^{d+1}(k,\mu_n^{\otimes d}) \cong \Zbb/n
    \end{equation*}
is a perfect duality of finite groups for $0 \le i \le d+1$. Suppose that $d \ge 1$ and that $F$ extends to a group scheme $\Fcal$ over the ring of integers $\Ocal$ of $k$ (thus $F'$ also extends to $\Fcal' = \iHom_{\Ocal}(\Fcal,\mu_n^{\otimes d})$). Then $\H^i(\Ocal,\Fcal)$ is a subgroup of $\H^i(k, F)$, $\H^{d+1-i}(\Ocal,\Fcal')$ is a subgroup of $\H^{d+i-i}(k, F')$, and these subgroups are exact annihilators of each other under the above duality pairing \cite[Proposition 2.5]{IzquierdoDualite}.

Let $K$ be the function field of a smooth projective geometrically integral curve $\Omega$ over $k$. If $j: U \hookrightarrow \Omega$ is a non-empty open subset and $\Fcal$ is a complex of sheaves on $U_{\et}$, we define the {\em hypercohomology groups with compact support} $\H^i_c(U,\Fcal) := \H^i(\Omega, j_{!} \Fcal)$, where $j_{!}$ denotes the extension by zero \cite[p. 93]{MilneEtale}. From the proof of \cite[Lemme 1.3]{IzquierdoArxiv}, one has a canonical isomorphism $\H^{d+3}_c(U,\Qbb/\Zbb(d+1)) \cong \Qbb/\Zbb$.

Let $F$ be a finite $\Gamma_K$-module. For a sufficiently small non-empty open subset $U \subseteq \Omega$, $F$ extends to a finite \'etale group scheme $\Fcal \to U$. Let $\Fcal' = \iHom_U(\Fcal,\Qbb/\Zbb(d+1))$, which extends the finite $\Gamma_K$-module $F' = \Hom_K(F,\Qbb/\Zbb(d+1))$. For $0 \le i \le d+3$, one has the Yoneda product pairing for cohomology with compact support (see \cite[p. 168]{MilneEtale})
    \begin{equation*}
        \sqdot: \Ext^i_U(\Fcal', \Qbb/\Zbb(d+1)) \times \H^{d+3-i}_c(U,\Fcal') \to \H^{d+3}_c(U, \Qbb/\Zbb(d+1)) \cong \Qbb/\Zbb.
    \end{equation*}
On the other hand, since $\Fcal \cong \iHom_U(\Fcal',\Qbb/\Zbb(d+1))$, the spectral sequence 
    \begin{equation*}
        \H^p(U,\iExt^q_U(\Fcal',\Qbb/\Zbb(d+1))) \Rightarrow \Ext^{p+q}_U(\Fcal',\Qbb/\Zbb(d+1))
    \end{equation*}
yields an edge map $\H^i(U,\Fcal) \to \Ext^i_U(\Fcal', \Qbb/\Zbb(d+1))$. Hence, we obtain a pairing
    \begin{equation} \label{eq:ArtinVerdierFinite}
        \pair{-,-}_{\AV}: \H^i(U,\Fcal) \times \H^{d+3-i}_c(U,\Fcal') \to \H^{d+3}_c(U,\Qbb/\Zbb(2)) \cong \Qbb/\Zbb,
    \end{equation}
which is in fact a perfect duality of finite groups \cite[Proposition 2.1]{IzquierdoDualite}. We refer to it as the {\em Artin--Verdier duality pairing}. This pairing induces a perfect duality pairing
    \begin{equation} \label{eq:PoitouTateFinite}
        \pair{-,-}_{\PT}: \Sha^i(K,F) \times \Sha^{d+3-i}(K,F') \to \Qbb/\Zbb,
    \end{equation}
of finite groups \cite[Th\'eor\`eme 2.4]{IzquierdoDualite}, as follows. Let $\eta \in \Sha^i(K,F)$ and $\alpha \in \Sha^{d+3-i}(K,F')$. If $U$ is sufficiently small, we may lift $\eta$ to an element $\eta_U \in \H^i(U,\Fcal)$ and $\alpha$ to an element $\alpha_U \in \H^{d+3-i}(U,\Fcal')$. By the localization exact sequence for cohomology with compact support
    \begin{equation} \label{eq:LocalizationFinite}
        \H^{d+3-i}_c(U,\Fcal') \to \H^{d+3-i}(U,\Fcal') \to \bigoplus_{v \notin U} \H^{d+3-i}(K_v,F')
    \end{equation}
(which can be proved by exactly the same argument as in \cite[Chapter II, Lemma 2.4]{MilneDuality}), $\alpha_U$ comes from an element $\alpha_U^c \in \H^{d+3-i}_c(U,\Fcal')$. Then $\pair{\eta,\alpha}_{\PT} = \pair{\eta_U,\alpha_U^c}_{\AV}$. The non-degeneracy of \eqref{eq:PoitouTateFinite} shall serve in the proof of Theorem \ref{customthm:HigherHasse}.

In order the prove Theorem \ref{customthm:HigherWeak}, one needs the exact sequence \cite[Lemme 1.2]{IzquierdoII}
    \begin{equation} \label{eq:ExactSequenceFinite}
        \H^i(K,F) \to \prod_{v \in S} \H^i(K_v,F) \xrightarrow{\theta} \Sha^{d+2-i}_S(K,F')^D \to \Sha^{d+2-i}(K,F')^D \to 0,
    \end{equation}
for any finite subset $S \subseteq \Omega^{(1)}$ and $1 \le i \le d+1$, which is established in the course of establishing the Poitou--Tate sequence for finite modules. Here, the map $\theta$ is defined by
    \begin{equation*}
        \forall (f_v)_{v \in S} \in \prod_{v \in S}\H^i(K_v,F), \forall \alpha \in \Sha^{d+2-i}_S(K,F'), \quad \theta((f_v)_{v \in S})(\alpha) = \sum_{v \in S} f_v \cup \loc_v(\alpha),
    \end{equation*}
where $f_v \cup \loc_v(\alpha) \in \H^{d+2}(K_v,\Qbb/\Zbb(d+1)) \cong \Qbb/\Zbb$ for each $v \in S$.

Let us now focus on the case where $k$ is a $p$-adic field, {\em i.e.} a finite extension of $\Qbb_p$ (hence $d = 1$). Let $T$ be a $K$-torus. Recall that the dual torus $T'$ is the torus with character module $\wh{T'} = \wc{T}$. For each $v \in \Omega^{(1)}$, local duality between the tori $T$ and $T'$ over the $2$-dimensional local field $K_v$ asserts that the pairing \eqref{eq:PairingTorus} induces a cup-product pairing
    \begin{equation} \label{eq:LocalDualityTorus}
        \H^1(K_v,T) \times \H^1(K_v,T') \to \H^4(K_v,\Zbb(2)) \cong \H^3(K_v,\Qbb/\Zbb(2)) \cong \Qbb/\Zbb,
    \end{equation}
which is a perfect duality of finite groups. For a sufficiently small non-empty open subset $U \subseteq \Omega$, $T$ (resp. $T'$) extends to a $U$-torus $\Tcal$ (resp. $\Tcal'$). Then $\H^1(\Ocal_v,\Tcal)$ is a subgroup of $\H^1(K_v,T)$, $\H^1(\Ocal,\Tcal')$ is a subgroup of $\H^1(K_v,T')$, and these subgroups are exact annihilators of each other under the above duality pairing (see \cite[Lemma 2.1, Proposition 2.2]{HShasse} and \cite[Proposition 2.3(a)]{HSS}).

The pairing \eqref{eq:PairingTorus} extends to a pairing $\Tcal \otimes^{\Lbb} \Tcal' \to \Gbb_m \otimes^{\Lbb} \Gbb_m \to \Zbb(2)[2]$ in $\Dcal^+(U)$. For $0 \le i \le 3$, we have a pairing 
    \begin{equation} \label{eq:ArtinVerdierTorus1}
        \Ext^i_U(\Tcal', \Tcal \otimes^{\Lbb} \Tcal') \times \H^{3-i}_c(U,\Tcal') \xrightarrow{\sqdot} \H^3_c(U, \Tcal \otimes^{\Lbb} \Tcal') \to \H^5_c(U,\Zbb(2)) \cong \Qbb/\Zbb,
    \end{equation}
where the first arrow is the Yoneda product pairing for cohomology with compact support \cite[p. 168]{MilneEtale}, the second arrow is induced by the above pairing, and the isomorphism $\H^5_c(U,\Zbb(2)) \cong \Qbb/\Zbb(2)$ is  \cite[Lemma 1.1]{HShasse}. Furthermore, we have a composite map
    \begin{equation} \label{eq:ArtinVerdierTorus2}
        \H^i(U,\Tcal) \to \H^i(U, \iHom_U(\Tcal', \Tcal \otimes^{\Lbb} \Tcal')) \to \Ext^i_U(\Tcal', \Tcal \otimes^{\Lbb} \Tcal'),
    \end{equation}
where the first arrow is induced by the natural morphism $\Tcal \to \iHom_U(\Tcal', \Tcal \otimes^{\Lbb} \Tcal')$ in $\Dcal^+(U)$, and the second arrow is the edge map from the spectral sequence
    \begin{equation*}
        \H^p(U, \iExt_U^q(\Tcal', \Tcal \otimes^{\Lbb} \Tcal')) \Rightarrow \Ext^{p+q}_U(\Tcal', \Tcal \otimes^{\Lbb} \Tcal').
    \end{equation*}
From \eqref{eq:ArtinVerdierTorus1} and \eqref{eq:ArtinVerdierTorus2}, we obtain a pairing
    \begin{equation} \label{eq:ArtinVerdierTorus}
        \pair{-,-}_{\AV}: \H^i(U,\Tcal) \times \H^{3-i}_c(U,\Tcal') \to \H^5_c(U,\Zbb(2)) \cong \Qbb/\Zbb
    \end{equation}
We refer to it as the {\em Artin--Verdier pairing (for tori)}. Unlike in the finite case, this is {\em not} a perfect duality (nor are the concerning cohomology groups finite). Nevertheless, for $i = 1$ (and by exchanging $\Tcal$ and $\Tcal'$), this induces a perfect duality pairing
    \begin{equation} \label{eq:PoitouTateTorus}
        \pair{-,-}_{\PT}: \Sha^2(K,T) \times \Sha^1(K,T') \to \Qbb/\Zbb
    \end{equation}
of finite groups, as follows (see Theorem 1.3 and the proof of Theorem 4.1 in \cite{HShasse}). Let $\eta \in \Sha^2(K,T)$ and $\alpha \in \Sha^1(K,T')$. If $U$ is sufficiently small, we may lift $\eta$ to an element $\eta_U \in \H^2(U,\Tcal)$ and $\alpha$ to an element $\alpha_U \in \H^1(U,\Tcal')$. By the localization exact sequence
    \begin{equation} \label{eq:LocalizationTorus}
        \cdots \to \H^i_c(U,\Tcal) \to \H^i (U,\Tcal) \to \bigoplus_{v \notin U} \H^i (K_v,T) \to \H^{i+1}_c(U,\Tcal) \to \cdots
    \end{equation}
starting from degree $1$ \cite[Corollary 3.2]{HShasse}, $\eta_U$ comes from an element $\eta_U^c \in \H^2_c(U,\Tcal)$. Then $\pair{\eta,\alpha}_{\PT} = \pair{\eta_U^c,\alpha_U}_{\AV}$. In this article, we shall need the fact that \eqref{eq:PoitouTateTorus} is also induced by \eqref{eq:ArtinVerdierTorus} for $i = 2$ (without exchanging $\Tcal$ and $\Tcal'$). To see this, it suffices to show the following

\begin{lm}
    We have a commutative diagram of pairings
        \begin{equation*}
            \xymatrix@C-2pc{
                \H^2_c(U,\Tcal) \ar[d] & \times & \H^1(U,\Tcal') \ar[rrrrrrrr]^-{\pair{-,-}_{\AV}} &&&&&&&& \H^5_c(U,\Zbb(2)) \ar@{=}[d] \\
                \H^2(U,\Tcal) & \times & \H^1_c(U,\Tcal') \ar[rrrrrrrr]^-{\pair{-,-}_{\AV}} \ar[u] &&&&&&&& \H^5_c(U,\Zbb(2)), 
            }
        \end{equation*}  
    that is, for all $\eta_U \in \H^2(U,\Tcal)$ coming from $\eta_U^c \in \H^2_c(U,\Tcal)$ and $\alpha_U \in \H^1(U,\Tcal')$ coming from $\alpha_U^c \in \H^1_c(U,\Tcal')$, one has $\pair{\eta_U^c,\alpha_U}_{\AV} = \pair{\eta_U,\alpha_U^c}_{\AV}$.
\end{lm}
\begin{proof}
    As in \cite[(26)]{HShasse}, for any $n \ge 1$, one has two commutative diagrams of pairings
        \begin{equation*}
            \xymatrix@C-2pc{
                \H^2_c(U,\tensor[_n]{}{}\Tcal) \ar[d]^{\iota_n} & \times & \H^2(U,\tensor[_n]{}{}\Tcal') \ar[rrrrrrrr]^-{\pair{-,-}_{\AV}} &&&&&&&& \H^4_c(U,\mu_n^{\otimes 2}) \ar[d]^{\partial} \\
                \H^2_c(U,\Tcal) & \times & \H^1(U,\Tcal') \ar[rrrrrrrr]^-{\pair{-,-}_{\AV}} \ar[u]_{\partial_n} &&&&&&&& \H^5_c(U,\Zbb(2)) 
            }
        \end{equation*}  
    and
        \begin{equation*}
            \xymatrix@C-2pc{
                \H^2(U,\tensor[_n]{}{}\Tcal) \ar[d]^{\iota_n} & \times & \H^2_c(U,\tensor[_n]{}{}\Tcal') \ar[rrrrrrrr]^-{\pair{-,-}_{\AV}} &&&&&&&& \H^4_c(U,\mu_n^{\otimes 2}) \ar[d]^{\partial} \\
                \H^2(U,\Tcal) & \times & \H^1_c(U,\Tcal') \ar[rrrrrrrr]^-{\pair{-,-}_{\AV}} \ar[u]_{\partial_n} &&&&&&&& \H^5_c(U,\Zbb(2)), 
            }
        \end{equation*}  
    where the top rows are the Artin--Verdier pairing for finite modules \eqref{eq:ArtinVerdierFinite}, and the maps $\iota_n$, $\partial_n$, $\partial$ are induced by the respective morphisms $\Zbb/n \to \Zbb(1)[1]$, $\Zbb(1)[1] \to \Zbb/n(1)[1]$, and $\Zbb/n(2) \to \Zbb(2)[1]$ in the derived category. On the other hand, by the same argument as in Lemma 4.2 and the proof of Lemma 4.7(3) in \cite{DH}, we have a commutative diagram of pairings
        \begin{equation*}
            \xymatrix@C-2pc{
                \H^2_c(U,\tensor[_n]{}{}\Tcal) \ar[d] & \times & \H^2(U,\tensor[_n]{}{}\Tcal') \ar[rrrrrrrr]^-{\pair{-,-}_{\AV}} &&&&&&&& \H^4_c(U,\mu_n^{\otimes 2}) \ar@{=}[d] \\
                \H^2(U,\tensor[_n]{}{}\Tcal) & \times & \H^2_c(U,\tensor[_n]{}{}\Tcal') \ar[rrrrrrrr]^-{\pair{-,-}_{\AV}} \ar[u] &&&&&&&& \H^4_c(U,\mu_n^{\otimes 2}).
            }
        \end{equation*}
    Let $\eta^c_U \in \H^2_c(U,\Tcal)$ (with image $\eta_U \in \H^2(U,\Tcal)$) and $\alpha^c_U \in \H^1_c(U,\Tcal)$ (with image $\alpha_U \in \H^1(U,\Tcal)$). Since the group $\H^2_c(U,\Tcal)$ is torsion by \cite[Corollary 3.3]{HShasse}, one has $\eta_U^c = \iota_n(\eta_{U,n}^c)$ for some $n \ge 1$ and some $\eta_{U,n}^c \in \H^2_c(U, \tensor[_n]{}{}\Tcal)$. Denote by $\eta_{U,n} \in \H^2(U, \tensor[_n]{}{}\Tcal)$ the image of $\eta_{U,n}^c$. From the commutativity of the above diagrams, one has
        \begin{equation*}
            \pair{\eta^c_U, \alpha_U}_{\AV} = \partial(\pair{\eta_{U,n}^c, \partial_n(\alpha_U)}_{\AV}) = \partial(\pair{\eta_{U,n}, \partial_n(\alpha_U^c)}_{\AV}) = \pair{\eta_U, \alpha_U^c}_{\AV}.
        \end{equation*}
    The Lemma is now proved.
\end{proof}

It follows that the pairing \eqref{eq:PoitouTateTorus} has the following alternative description. Let $\eta \in \Sha^2(K,T)$ and  $\alpha \in \Sha^1(K,T')$. Lift $\eta$ to an element $\eta_U \in \H^2(U,\Tcal)$ and $\alpha$ to an element $\alpha_U \in \H^1(U,\Tcal')$ (shrinking $U$ if necessary). By the localization sequence \eqref{eq:LocalizationTorus}, $\alpha_U$ comes from an element $\alpha_U^c \in \H^1_c(U,\Tcal')$. Then $\pair{\eta,\alpha}_{\PT} = \pair{\eta_U, \alpha_U^c}_{\AV}$.

The non-degeneracy of \eqref{eq:PoitouTateTorus} shall serve in the proof of Theorem \ref{customthm:Hasse}. As for Theorem \ref{customthm:Descent}, one needs the following part of the Poitou--Tate exact sequence for tori \cite[Proposition 3.5]{HSS}:
    \begin{equation} \label{eq:PoitouTateSequenceTorus}
        \H^1(K,T) \to \Pbb^1(K,T) \xrightarrow{\theta} \H^1(K,T')^D.
    \end{equation}
This is an exact sequence of {\em topological abelian groups}. The groups $\H^1(K,T)$ and $\H^1(K,T')$ are discrete. The group $\Pbb^1(K, T)$ is, by definition, the topological restricted product of the finite groups $\H^1(K_v, T)$ relative to the subgroups $\H^1(\Ocal_v,\Tcal)$, $v \in U^{(1)}$ (recall that $\Tcal$ is a $U$-torus extending $T$). The map $\theta$ is defined by
    \begin{equation*}
        \forall (t_v)_{v \in \Omega^{(1)}} \in \Pbb^1(K,T), \forall \alpha \in \H^1(K,T'), \quad \theta((t_v)_{v \in \Omega^{(1)}})(\alpha) = \sum_{v \in \Omega^{(1)}} t_v \cup \loc_v(\alpha).
    \end{equation*}
Here, $t_v \cup \loc_v(\alpha) \in \H^4(K_v,\Zbb(2)) \cong \Qbb/\Zbb$ {\em via} the local cup-product \eqref{eq:LocalDualityTorus}.

In order to prove Theorem \ref{customthm:Weak}, we require the following exact sequence. Let $S \subseteq \Omega^{(1)}$ be any finite set. Since any element of $\prod_{v \in S} \H^1(K_v,T)$ can completed by $0$ into an element of $\Pbb^1(K,T)$, \eqref{eq:PoitouTateSequenceTorus} restricts to the three last terms of the exact sequence
    \begin{equation*} 
        0 \to \Sha^1(K,T) \to \Sha^1_S(K,T) \to \prod_{v \in S} \H^1(K_v,T) \to \H^1(K,T')^D
    \end{equation*}
of discrete abelian groups. Dualizing this sequence and exchanging $T$ and $T'$, one obtains an exact sequence
    \begin{equation} \label{eq:ExactSequenceTorus}
           \H^1(K,T) \to \prod_{v \in S} \H^1(K_v,T) \xrightarrow{\theta} \Sha^1_S(K,T')^D \to \Sha^1(K,T')^D \to 0,
    \end{equation}
where the map $\theta$ is defined by
    \begin{equation*}
        \forall (t_v)_{v \in S} \in \prod_{v \in S} \H^1(K_v,T), \forall \alpha \in \Sha^1_S(K,T'), \quad \theta((t_v)_{v \in S})(\alpha) = \sum_{v \in S} t_v \cup \loc_v(\alpha).
    \end{equation*}

We conclude this section with the following lemma, which is \cite[Lemma 4.2(a)]{HSS}, whose proof relies on Tate--Lichtenbaum duality for $p$-adic curves. It is crucial for the constructions used in the proofs of Theorems \ref{customthm:Hasse}, \ref{customthm:Weak}, and \ref{customthm:Descent}.    

\begin{lm} \label{lm:QuasiTrivialTorusPAdic}
    If $Q$ is a quasi-split torus over a $p$-adic function field $K$, then  $\Sha^2_\omega(K,Q) = 0$.
\end{lm}

This lemma is used not only to establish the unramified nature of the obstructions but also the Poitou--Tate sequence for tori. For $K$ a function field of a curve over a $d$-dimensional local field (where $d \ge 2$), one would need the vanishing of $\Sha_\omega^{d+2}(K,\Zbb(d)) = \Sha^{d+2}(K,\Zbb(d))$, where $\Zbb(d)$ is the shifted weight $d$ cycle complex $z(\bullet,-)[-2d]$ defined by Bloch \cite{Bloch}, or equivalently, the vanishing of $\Sha^2(K,\Gbb_m)$ \cite[Lemme 3.15]{IzquierdoArxiv}. Unfortunately, this is not always the case (this problem was studied by Izquierdo in \cite[\S 4 and \S 5]{IzquierdoI}).

\begin{rmk}
    An independent but interesting consequence of Lemma \ref{lm:QuasiTrivialTorusPAdic} is the following local-global principle. Let $X$ be a Severi--Brauer variety over the function field $K$ of a $p$-adic curve $\Omega$. If $X(K_v) \neq \varnothing$ for all but finitely many $v \in \Omega^{(1)}$, then $X(K) \neq \varnothing$.
\end{rmk}

\section{Descent theory} \label{sec:Descent}

This section is devoted to the proof of Theorem \ref{customthm:Descent}. 

\subsection{Preliminary remarks} \label{subsec:DescentRemark}

We recall some facts. Let $X$ be a smooth geometrically integral variety over a field $K$ of characteristic $0$. Following Colliot-Th\'el\`ene and Sansuc ({\em cf.} \cite[Theorem 2.3.4, Definition 2.3.5]{SkorobogatovTorsors}), we define the {\em elementary obstruction} $e_X \in \Ext^2_K(\Pic\ol{X},\ol{K}[X]^\times)$ to be the inverse class of the $2$-fold extension
    \begin{equation} \label{eq:ElementaryObstruction}
        1 \to \ol{K}[X]^\times \to \ol{K}(X)^\times \xrightarrow{\div} \Div\ol{X} \to \Pic \ol{X} \to 0
    \end{equation}
of $\Gamma_K$-modules. If $M$ is a $K$-group of multiplicative type, recall that $\wh{M} = \iHom_K(M,\Gbb_m)$. The {\em type} of a torsor $Y \to X$ under $M$ is by definition the $\Gamma_K$-equivariant homomorphism 
    \begin{equation*}
        \wh{M} \to \Pic\ol{X} = \H^1(\ol{X},\Gbb_m), \qquad \chi \mapsto \chi_\ast [\ol{Y}].
    \end{equation*}
Conversely, for any given $\Gamma_K$-equivariant homomorphism $\lambda: \wh{M} \to \Pic\ol{X}$, the existence of $X$-torsors under $M$ of type $\lambda$ is equivalent to the vanishing of $\lambda^\ast e_X \in \Ext^2_K(\wh{M},\ol{K}[X]^\times)$. In the case where $X(K) \neq \varnothing$, we have $e_X = 0$ and torsors of any type exist. When $\ol{K}[X]^\times = \ol{K}^\times$, the spectral sequence $\H^p(K,\iExt^q_K(\wh{M},\Gbb_m)) \Rightarrow \Ext^{p+q}_K(\wh{M},\Gbb_m)$ yield the edge maps
    \begin{equation} \label{eq:ExtOfGroupOfMultiplicativeType}
       \H^p(K,M) = \H^p(K,\iHom_K(\wh{M},\Gbb_m)) \xrightarrow{\cong} \Ext^p_K(\wh{M},\Gbb_m),
    \end{equation}
which are isomorphisms for $p > 0$ \cite[Lemma 2.3.7]{SkorobogatovTorsors}. In particular, we may regard $\lambda^\ast e_X$ as an element of $\H^2(K,M)$. The equivalence between the vanishing of this element and the existence of $X$-torsors of type $\lambda$ is part of the ``fundamental exact sequence'' of Colliot-Th\'el\`ene and Sansuc (see Theorem 2.3.6 and Corollary 2.3.9 in {\em loc. cit.}), which reads 
    \begin{equation} \label{eq:FundamentalExactSequence}
        \H^1(K,M) \to \H^1(X,M) \xrightarrow{\type} \iHom_K(\wh{M},\Gbb_m) \to \H^2(K,M) \to \H^2(X,M).
    \end{equation}
If $\Pic\ol{X}$ is finitely generated as an abelian group and $M$ is the $K$-group of multiplicative type such that $\wh{M} = \Pic\ol{X}$ (that is, an isomorphism $\wh{M} \cong \Pic\ol{X}$ is {\em fixed}), we call a torsor $Y \to X$ under $M$ {\em universal} if its type is the identity morphism of $\wh{M}$\footnote{This differs slightly from the usual convention, where a torsor is defined to be universal if its type is {\em any} isomorphism of Galois module (not just a fixed one from the beginning).}. Indeed, the existence of such a torsor is equivalent to $e_X = 0$. 

Assume furthermore that $\Pic\ol{X}$ is free,  then the {\it N\'eron--Severi torus} of $X$ is by definition the $K$-torus $T$ such that $\wh{T} = \Pic\ol{X}$. For example, this is case when $X$ is projective and {\it rationally connected} (combine \cite[\S 8.4, Theorem 1]{BLR}, \cite[Corollary 4.18]{Debarre}, and \cite[Th\'eor\`eme 5.1]{Kleiman}), that is, for any algebraically closed overfield $K'/K$ and two general points $P_0,P_1 \in X(K')$, there exists a morphism $\gamma: \Pbb^1_{K'} \to X_{K'}$ such that $\gamma(0) = P_0$ and $\gamma(1) = P_1$. Examples of such varieties are smooth compactifications of geometrically unirational varieties (such as homogeneous spaces of connected linear algebraic groups; indeed, a celebrated theorem of Chevalley asserts that connected linear algebraic groups are geometrically rational, even $K$-unirational \cite{Chevalley}). 

Before starting, let us restate the main result of this section ({\em i.e.} Theorem \ref{customthm:Descent}).

\begin{thm} \label{thm:Descent}
    Let $K$ be the function field of a smooth proper geometrically integral curve $\Omega$ over a $p$-adic field $k$, and let $X$ be a smooth proper geometrically integral variety over $K$ such that the abelian group $\Pic\ol{X}$ is finitely generated and free (for example, $X$ is projective and rationally connected). Let $T$ be the N\'eron--Severi torus of $X$. There exists a homomorphism
    \begin{equation*}
        u: \H^1(K,T') \to \tfrac{\H^4(X,\Zbb(2))}{\Img \H^4(K,\Zbb(2))}
    \end{equation*}
    with the following properties. Suppose that $\prod_{v \in \Omega^{(1)}} X(K_v) \neq \varnothing$, then
    \begin{enumerate}
        \item \label{thm:Descent1} universal $X$-torsors exist (that is, $e_X = 0$) if and only if there exists a family of $\prod_{v \in \Omega^{(1)}} X(K_v)$ which is orthogonal to $u(\Sha^1(K,T'))$ relative to the pairing \eqref{eq:UnramifiedPairingGersten};

        \item \label{thm:Descent2} if $\left(\prod_{v \in \Omega^{(1)}} X(K_v) \right)^{\Img(u)}$ denotes the subset of $\prod_{v \in \Omega^{(1)}} X(K_v)$ consisting of the families orthogonal to $\Img(u)$ relative to the pairing \eqref{eq:UnramifiedPairingGersten}, then
            \begin{equation*}
                \left(\prod_{v \in \Omega^{(1)}} X(K_v) \right)^{\Img(u)} = \bigcup_{\substack{f: Y \to X \\ \type(Y) = \id}} f\left(\prod_{v \in \Omega^{(1)}} Y(K_v)\right);
            \end{equation*}

        \item \label{thm:Descent3} there are only finitely many isomorphism classes of universal torsors $Y \to X$ such that $\prod_{v \in \Omega^{(1)}} Y(K_v) \neq \varnothing$;

        \item \label{thm:Descent4} if the universal torsors $Y \to X$ satisfy the local-global principle (resp. the local-global principle and weak approximation), then the reciprocity obstruction \eqref{eq:UnramifiedPairingGersten} to the local-global principle (resp. the local-global principle and weak approximation) on $X$ attached to $\Img(u)$ is the only one.
    \end{enumerate}
\end{thm}

Let $K$ be a field of cohomological dimension $\cd(K) \le 3$ and let $\pi: X \to \Spec K$ be a smooth proper geometrically integral variety such that the abelian group $\Pic\ol{X}$ is finitely generated and free. Let $T$ be the N\'eron--Severi torus of $X$ (that is, $\wh{T} = \Pic\ol{X}$). We construct a map 
    \begin{equation} \label{eq:MapU}
        u: \H^1(K,T') \to \tfrac{\H^4(X,\Zbb(2))}{\Img \H^4(K,\Zbb(2))}
    \end{equation}
as in the statement of Theorem \ref{thm:Descent}, as follows. First, since $\ol{K}[X]^\times = \ol{K}^\times$, we have the following distinguished triangle in the category $\Dcal^+(K)$:
    \begin{equation} \label{eq:TriangleHatT}
        \Gbb_m \to \tau_{\le 1} \Rbb \pi_\ast \Gbb_{m,X} \to \wh{T}[-1] \to \Gbb_m[1].
    \end{equation}
Applying the exact functor $- \otimes^{\Lbb} \Gbb_m$ yields a distinguished triangle
    \begin{equation} \label{eq:MapU1}
        \Gbb_m \otimes^{\Lbb} \Gbb_m \to (\tau_{\le 1} \Rbb \pi_\ast \Gbb_{m,X}) \otimes^{\Lbb} \Gbb_m \to T'[-1] \to \Gbb_m \otimes^{\Lbb} \Gbb_m[1].
    \end{equation}
Let $\theta_1: \Rbb \pi_\ast (\Gbb_{m,X} \otimes^{\Lbb} \Gbb_{m,X}) \to \Rbb \pi_\ast \Zbb_X(2)[2]$ be the map induced by the pairing \eqref{eq:PairingMotivic}. Next, denote by $\theta_2$ the composite
    \begin{align*}
        (\tau_{\le 1} \Rbb \pi_\ast \Gbb_{m,X}) \otimes^{\Lbb} \Gbb_m & \to (\tau_{\le 1} \Rbb \pi_\ast \Gbb_{m,X}) \otimes^{\Lbb} \Rbb \pi_\ast \Gbb_{m,X} \to \Rbb \pi_\ast ((\pi^\ast \tau_{\le 1} \Rbb \pi_\ast \Gbb_{m,X}) \otimes^{\Lbb} \Gbb_{m,X}) \\
         & \to \Rbb \pi_\ast ((\pi^\ast \Rbb\pi _\ast \Gbb_{m,X}) \otimes^{\Lbb} \Gbb_{m,X}) \to \Rbb \pi_\ast (\Gbb_{m,X} \otimes^{\Lbb} \Gbb_{m,X}),
    \end{align*}
where 
    \begin{itemize}
        \item the first arrow is induced by the natural map $\Gbb_m \to \Rbb\pi_\ast \Gbb_{m,X}$,
        \item the second arrow is the canonical ``base change'' morphism constructed in \cite[p. 306]{Fu},
        \item the third arrow is induced by the natural map $\tau_{\le 1} \Rbb \pi_\ast \Gbb_{m,X} \to \Rbb \pi_\ast \Gbb_{m,X}$,
        \item and the last arrow is induced by the adjunction $\pi^\ast \Rbb \pi_\ast \Gbb_{m,X} \to \Gbb_{m,X}$.
    \end{itemize}
Finally, let
    \begin{equation} \label{eq:MapTheta}
        \theta = \theta_1 \circ \theta_2: (\tau_{\le 1} \Rbb \pi_\ast \Gbb_{m,X}) \otimes^{\Lbb} \Gbb_m \to \Rbb \pi_\ast \Zbb_X(2)[2].
    \end{equation}
By the functoriality of the pairing \eqref{eq:PairingMotivic}, $\theta$ fits into a commutative diagram
    \begin{equation*}
        \xymatrix{
        \Gbb_m \otimes^{\Lbb} \Gbb_m \ar[r] \ar[d] & (\tau_{\le 1} \Rbb \pi_\ast \Gbb_{m,X}) \otimes^{\Lbb} \Gbb_m \ar[d]^{\theta} \\
        \Zbb(2)[2] \ar[r] & \Rbb \pi_\ast \Zbb_X(2)[2].
    }\end{equation*}
Let $\Zbb_{X/K}(2)$ denote the cone of $\Zbb_X(2) \to \Rbb \pi_\ast \Zbb_X(2)$. It follows from the axioms of triangulated categories that there exists a dashed arrow making the diagram
    \begin{equation} \label{eq:MapU2}
        \xymatrix{
            \Gbb_m \otimes^{\Lbb} \Gbb_m \ar[r] \ar[d] & (\tau_{\le 1} \Rbb \pi_\ast \Gbb_{m,X}) \otimes^{\Lbb} \Gbb_m \ar[d]^{\theta} \ar[r] & T'[-1] \ar@{-->}[d]^{\lambda} \ar[r] & \Gbb_m \otimes^{\Lbb} \Gbb_m[1] \ar[d] \\
            \Zbb(2)[2] \ar[r] & \Rbb \pi_\ast \Zbb_X(2)[2] \ar[r] & \Zbb_{X/K}(2)[2] \ar[r] & \Zbb(2)[3]
        }
    \end{equation}
commute (the top row being \eqref{eq:MapU1}). Since  $\H^5(K,\Zbb(2)) \cong \H^4(K,\Qbb/\Zbb(2)) = 0$ under the assumption $\cd(K) \le 3$, taking cohomology of the bottom row of \eqref{eq:MapU2} yields an identification 
    \begin{equation} \label{eq:MapU3}
        \H^4(K,\Zbb_{X/K}(2)) \cong \tfrac{\H^4(X,\Zbb(2))}{\Img \H^4(K,\Zbb(2))}.
    \end{equation} 
We take the map $u$ in $\eqref{eq:MapU}$ to be the composite
    \begin{equation*}
        \H^1(K,T') \xrightarrow{\lambda_\ast} \H^4(K,\Zbb_{X/K}(2)) \cong \tfrac{\H^4(X,\Zbb(2))}{\Img \H^4(K,\Zbb(2))}.
    \end{equation*}

\subsection{Existence of universal torsors} \label{sub:Existence}

In this paragraph, we prove Theorem \ref{thm:Descent}\ref{thm:Descent1}. Let $K$ be the function field of a smooth projective geometrically integral curve $\Omega$ over a $p$-adic field $k$. The point is to relate the first obstruction \eqref{eq:UnramifiedReciprocityGersten} and the global Poitou--Tate duality pairing \eqref{eq:PoitouTateTorus}, as in the following analogue of \cite[Lemme 3.3.3]{CTSdescent} (see also \cite[(6.4)]{SkorobogatovTorsors}).

\begin{prop} \label{prop:Descent1}
    Let $\pi: X \to \Spec K$ be a smooth proper geometrically integral variety such that the abelian group $\Pic\ol{X}$ is finitely generated and free. Let $T$ be the N\'eron--Severi torus of $X$. Assume in addition that $\prod_{v \in \Omega^{(1)}}X(K_v) \neq \varnothing$. In particular, the class $\eta \in \H^2(K,T)$ corresponding to the elementary obstruction $e_X \in \Ext^2_K(\wh{T},\Gbb_m)$ (under the identification \eqref{eq:ExtOfGroupOfMultiplicativeType}) belongs to $\Sha^2(K,T)$. Then, for all $\alpha \in \Sha^{1}(K,T')$, one has the equality
        \begin{equation*}
            \rho_X(u(\alpha)) = -\pair{\eta,\alpha}_{\PT}.
        \end{equation*}
    Here, the map $\rho_X$ was defined in \eqref{eq:UnramifiedReciprocityGersten}, and $\pair{-,-}_{\PT}$ is the pairing \eqref{eq:PoitouTateTorus}.
\end{prop}
\begin{proof}
    We follow the argument in \cite[\S 3]{HSmotive} and \cite[Proposition 5.3]{HShasse}. First, we inspect the pairing $\pair{\eta,-}_{\PT}$. By \cite[Lemma 2.3]{BvHcrelle}, the object $\tau_{\le 1} \Rbb \pi_\ast \Gbb_{m,X}$ in \eqref{eq:TriangleHatT} is represented by the complex $[\ol{K}(X)^\times \xrightarrow{\div} \Div \ol{X}]$ concentrated in degree $-1$ and $0$. It follows that the class $-e_X \in \Ext^2_K(\wh{T},\Gbb_m)$ of the $2$-fold extension \eqref{eq:ElementaryObstruction} (note that $\ol{K}[X]^\times = \ol{K}^\times$) is also represented by a morphism $\wh{T} \to \Gbb_m[2]$ in $\Dcal^+(K)$ associated with  triangle \eqref{eq:TriangleHatT}.

    Let $\pi^U: \Xcal \to U$ be an integral model of $X$ over some non-empty open subset $U \subseteq \Omega$. We may assume that $T$ extends to a $U$-torus $\Tcal$. Then $\wh{T}$ (resp. $T'$) extends to the finitely presented locally constant group scheme $\wh{\Tcal} = \iHom_U(\Tcal,\Gbb_m)$ (resp. the $U$-torus $\Tcal' = \wh{\Tcal} \otimes \Gbb_m$). For $U$ sufficiently small, $\pi^U_\ast \Gbb_{m,\Xcal} = \Gbb_m$ and $\Rbb^1\pi^U_\ast \Gbb_{m,\Xcal} = \wh{\Tcal}$, hence one has a distinguished triangle
        \begin{equation} \label{eq:Descent11}
            \Gbb_m \to \tau_{\le 1} \Rbb \pi^U_\ast \Gbb_{m,\Xcal} \to \wh{\Tcal}[-1] \to \Gbb_m[1].
        \end{equation}
    in $\Dcal^+(U)$, which extends \eqref{eq:TriangleHatT}. The inverse class $e_U \in \Ext^2_U(\wh{\Tcal},\Gbb_m)$ of the morphism $\wh{\Tcal} \to \Gbb_m[2]$ associated with \eqref{eq:Descent11} is a lifting of $e_X \in \Ext^2_K(\wh{T},\Gbb_m)$.  

    We claim that there is a commutative diagram
        \begin{equation} \label{eq:Descent12}
            \xymatrix@C-1pc{
                & \H^2(U,\iHom_U(\wh{\Tcal},\Gbb_m)) \ar[r]^-{r_1} & \Ext^2_U(\wh{\Tcal},\Gbb_m) \ar@{=}[rd] \\
                \H^2(U,\Tcal) \ar[r] \ar[ru]^{\cong} \ar[rd] & \H^2(U,\iHom_U(\wh{\Tcal},\Tcal \otimes^{\Lbb} \wh{\Tcal})) \ar[u] \ar[d] \ar[r]^-{r_2} & \Ext^2_U(\wh{\Tcal},\Tcal \otimes^{\Lbb} \wh{\Tcal}) \ar[d]^{-\otimes^{\Lbb} \Gbb_m} \ar[u]^{\gamma_1} \ar[r]^-{\gamma_1} & \Ext^2_U(\wh{\Tcal},\Gbb_m) \ar[d]^{-\otimes^{\Lbb} \Gbb_m} \\
                & \H^2(U,\iHom_U(\Tcal',\Tcal \otimes^{\Lbb} \Tcal')) \ar[r]^-{r_3} & \Ext^2_U(\Tcal', \Tcal \otimes^{\Lbb} \Tcal') \ar[r]^-{\gamma_2} & \Ext^2_U(\Tcal',\Gbb_m \otimes^{\Lbb} \Gbb_m) \ar[d]^{\gamma_3} \\
                &&& \Ext^4_U(\Tcal',\Zbb(2))
            }
        \end{equation}    
    The two left triangles of \eqref{eq:Descent12} are obtained by applying $\H^2(U,-)$ to the commutative diagram
        \begin{equation*} 
            \xymatrix{
                && \Tcal \ar[lld]_-{\cong} \ar[d] \ar[rrd] \\
                \iHom_U(\wh{\Tcal},\Gbb_m) && \iHom_U(\wh{\Tcal}, \Tcal \otimes^{\Lbb} \wh{\Tcal}) \ar[ll] \ar[rr]^-{- \otimes^{\Lbb} \Gbb_m} && \iHom_U(\Tcal', \Tcal \otimes^{\Lbb} \Tcal')
            }
        \end{equation*}
    in $\Dcal^+(U)$. The map $r_1, r_2, r_3$ in \eqref{eq:Descent12} are the edge maps from the spectral sequences 
        \begin{equation*}
            \H^p(U,\iExt^q_U(\wh{\Tcal},\Fcal)) \Rightarrow \Ext^{p+q}_U(\wh{\Tcal},\Fcal),
        \end{equation*}
    for $\Fcal = \Gbb_m, \Tcal \otimes^{\Lbb} \wh{\Tcal}, \Tcal \otimes^{\Lbb} \Tcal'$ respectively. The two middle squares of \eqref{eq:Descent12} commute by the functoriality of these spectral sequences. The map $\gamma_1, \gamma_2, \gamma_3$ in \eqref{eq:Descent12} are induced by the respective pairings $\Tcal \otimes^{\Lbb} \wh{\Tcal} \to \Gbb_m$, $\Tcal \otimes^{\Lbb} \Tcal' \to \Gbb_m \otimes^{\Lbb} \Gbb_m$, and \eqref{eq:PairingMotivic}. The other triangle and square of \eqref{eq:Descent12} obviously commute. Since $r_1$ is an isomorphism by \cite[Lemma 2.3.7]{SkorobogatovTorsors}, the class $e_U \in \Ext^2_U(\wh{\Tcal},\Gbb_m)$ comes from an element $\eta_U \in \H^2(U,\Tcal)$ lifting $\eta \in \H^2(K,T)$. Let $\varepsilon_U, \varepsilon_U'$ denote its respective images in $\Ext^2_U(\wh{\Tcal},\Tcal \otimes^{\Lbb} \wh{\Tcal})$ and $\Ext^2_U(\Tcal',\Tcal \otimes^{\Lbb} \Tcal')$ by \eqref{eq:Descent12}. Then $-\gamma_2(\varepsilon_U')$ is represented by the morphism $\Tcal' \to \Gbb_m \otimes^{\Lbb} \Gbb_m[2]$ associated with the distinguished triangle obtained by applying $-\otimes^{\Lbb} \Gbb_m$ to \eqref{eq:Descent11}, {\em i.e.}
        \begin{equation} \label{eq:Descent13}
            \Gbb_m \otimes^{\Lbb} \Gbb_m \to (\tau_{\le 1} \Rbb \pi^U_\ast \Gbb_{m,\Xcal}) \otimes^{\Lbb} \Gbb_m \to \Tcal'[-1] \to \Gbb_m \otimes^{\Lbb} \Gbb_m [1],
        \end{equation}
    and $-\gamma_3(\gamma_2(\varepsilon'_U))$ is represented by the composite $\Tcal' \to \Zbb(2)[4]$ of \eqref{eq:PairingMotivic} with this morphism. On the other hand, we have a commutative diagram of pairings
        \begin{equation} \label{eq:Descent14}
            \xymatrix@C-2pc{
                \H^2(U,\Tcal) \ar[d] & \times & \H^1_c(U,\Tcal') \ar@{=}[d] \ar[rrrr] &&&& \H^3_c(U,\Tcal \otimes^{\Lbb} \Tcal') \ar@{=}[d]  \\
                \Ext^2_U(\Tcal',\Tcal \otimes^{\Lbb} \Tcal') \ar[d]^{\gamma_2} & \times & \H^1_c(U,\Tcal') \ar@{=}[d] \ar[rrrr]^-{\sqdot} &&&& \H^3_c(U,\Tcal \otimes^{\Lbb} \Tcal') \ar[d] \\
                \Ext^2_U(\Tcal',\Gbb_m \otimes^{\Lbb} \Gbb_m) \ar[d]^{\gamma_3} & \times & \H^1_c(U,\Tcal') \ar@{=}[d] \ar[rrrr]^-{\sqdot} &&&& \H^3_c(U,\Gbb_m \otimes^{\Lbb} \Gbb_m) \ar[d] \\
                \Ext^4_U(\Tcal',\Zbb(2)) & \times & \H^1_c(U,\Tcal') \ar[rrrr]^-{\sqdot} &&&& \H^5_c(U,\Zbb(2)) \cong \Qbb/\Zbb 
            }
        \end{equation}    
    where $\sqdot$ means the Yoneda product, and where the top square commutes thanks to the construction of the cup-product (Artin--Verdier) pairing for cohomology with compact support (see \eqref{eq:ArtinVerdierTorus1} and \eqref{eq:ArtinVerdierTorus2}). Hence, we have the following equality for all $\alpha_U^c \in \H^1_c(U,\Tcal')$:
        \begin{equation} \label{eq:Descent15}
            \gamma_3(\gamma_2(\varepsilon'_U)) \sqdot \alpha_U^c = \pair{\eta_U, \alpha_U^c}_{\AV} \in \H^5_c(U,\Zbb(2)) \cong \Qbb/\Zbb.
        \end{equation}
   Let $\alpha \in \Sha^1(K,T')$. By the localization sequence \eqref{eq:LocalizationTorus}, when $U$ is sufficiently small, $\alpha$ lifts to an element $\alpha_U^c \in \H^1_c(U,\Tcal)$. The right hand side of \eqref{eq:Descent15} is $\pair{\eta,\alpha}_{\PT}$ by the discussion following the construction of \eqref{eq:PoitouTateTorus} in paragraph \ref{subsec:Duality}.
   
   The next step is to inspect the element $\rho_X(u(\alpha))$. Consider the commutative diagram
    \begin{equation} \label{eq:Descent16}
        \xymatrix@C-1pc{
            & \H^4(K,\Zbb(2)) \ar[r] \ar[d] & \H^4(X,\Zbb(2)) \ar[r] \ar[d] & \frac{\H^4(X,\Zbb(2))}{\Img \H^4(K,\Zbb(2))} \ar[r] \ar[d] & 0 \\
            0 \ar[r] & \prod\limits_{v \in \Omega^{(1)}} \H^4(K_v,\Zbb(2)) \ar[r] & \prod\limits_{v \in \Omega^{(1)}} \H^4(X_v,\Zbb(2)) \ar[r] & \prod\limits_{v \in \Omega^{(1)}} \frac{\H^4(X_v,\Zbb(2))}{\Img \H^4(K_v,\Zbb(2))}, 
        }
    \end{equation}
    with exact rows (each map $\H^4(K_v,\Zbb(2)) \to \H^4(X_v,\Zbb(2))$ is injective since $X(K_v) \neq \varnothing$). Since $\alpha \in \Sha^1(K,T')$, $u(\alpha)$ lies in the kernel of the right vertical arrow in \eqref{eq:Descent16}. Let $\beta \in \frac{\prod_{v \in \Omega^{(1)}} \H^4(K_v,\Zbb(2))}{\Img \H^4(K,\Zbb(2))}$ be its image by the snake lemma construction.

    \begin{lm} \label{lm:Descent1}
        We have $\beta \in \frac{\bigoplus_{v \in \Omega^{(1)}} \H^4(K_v,\Zbb(2))}{\Img \H^4(K,\Zbb(2))}$, and its image by the sum map
            \begin{equation*}
                \sigma: \tfrac{\bigoplus_{v \in \Omega^{(1)}} \H^4(K_v,\Zbb(2))}{\Img \H^4(K,\Zbb(2))} \to \Qbb/\Zbb.
            \end{equation*}
        is precisely $\rho_X(u(\alpha))$.
    \end{lm}
    \begin{proof}
        Let $A \in \H^4(X,\Zbb(2))$ be a lifting of $u(\alpha)$. For each $v \in \Omega^{(1)}$, choose any point $P_v \in X(K_v)$. Then the constant element $\loc_v(A) \in \H^4(X_v,\Zbb(2))$ comes from $A(P_v) \in \H^4(K_v,\Zbb(2))$. By definition of the snake lemma construction, the family $(A(P_v))_{v \in \Omega^{(1)}}$ is a lifting of $\beta$. Thanks to \cite[Lemma 5.1]{HSS}, we have $A(P_v) = 0$ for all but finitely many $v \in \Omega^{(1)}$, so that $\beta \in \frac{\bigoplus_{v \in \Omega^{(1)}} \H^4(X_v,\Zbb(2))}{\H^4(K_v,\Zbb(2))}$. Finally, $\sigma(\beta) = \sum_{v \in \Omega^{(1)}} A(P_v) = \rho_X(u(\alpha))$.
    \end{proof} 
   Return to the proof of Proposition \ref{prop:Descent1}. Now we study the element $\beta$ by repeating the argument in \cite[Lemma 5.4]{HShasse}. Diagram \eqref{eq:MapU2} extends to a commutative diagram
    \begin{equation} \label{eq:Descent17}
        \xymatrix{
            \Gbb_m \otimes^{\Lbb} \Gbb_m \ar[r] \ar[d] & (\tau_{\le 1} \Rbb \pi_\ast^U \Gbb_{m,\Xcal}) \otimes^{\Lbb} \Gbb_m \ar[d] \ar[r] & \Tcal'[-1] \ar[d]^{\lambda_U} \ar[r] & \Gbb_m \otimes^{\Lbb} \Gbb_m[1] \ar[d] \\
            \Zbb(2)[2] \ar[r] & \Rbb \pi_\ast^U \Zbb_{\Xcal}(2)[2] \ar[r] & \Zbb_{\Xcal/U}(2)[2] \ar[r] & \Zbb(2)[3]
        }
    \end{equation}
    in $\Dcal^+(U)$, with distinguished rows (the top row being \eqref{eq:Descent13}; we recall that $\Tcal$ is a $U$-torus extending $T$, and the arrow $\Tcal' \to \Gbb_m \otimes^{\Lbb} \Gbb_m[2]$ represents $-\gamma_2(\varepsilon_U')$). Remember that $\alpha_U^c \in \H^1_c(U,\Tcal')$ is a lifting of $\alpha \in \Sha^1(K,T')$. For $v \in \Omega^{(1)}$, denote by $j_v:\Spec K^{\h}_v  \to U$ be the natural map, and consider the commutative diagram
        \begin{equation} \label{eq:Descent18}
        \xymatrix{
            \Zbb(2)[2] \ar[r] \ar[d] & \Rbb \pi_\ast^U \Zbb_{\Xcal}(2)[2] \ar[r] \ar[d] & \Zbb_{\Xcal/U}(2)[2] \ar[d] \\
            \bigoplus\limits_{v \notin U} j_{v\ast} j_v^{\ast} \Zbb(2)[2] \ar[r] & \bigoplus\limits_{v \notin U} j_{v\ast} j_v^{\ast} \Rbb \pi_\ast^U \Zbb_{\Xcal}(2)[2] \ar[r] & \bigoplus\limits_{v \notin U} j_{v\ast} j_v^{\ast} \Zbb_{\Xcal/U}(2)[2],
        }
        \end{equation}
    whose rows are parts of exact triangles. Let $\Ccal_l$ (resp. $\Ccal_r$) denote the cone of the left (resp. right) vertical arrow of \eqref{eq:Descent18}. Using the localization sequence
        \begin{equation*}
            \cdots \to \H^4_c(U,\Zbb_{\Xcal/U}(2)) \to \H^4(U,\Zbb_{\Xcal/U}(2)) \to \bigoplus_{v \notin U} \H^4(K_v^{\h}, j_v^\ast \Zbb_{\Xcal/U}(2)) \to \cdots 
        \end{equation*}
    (see for example \cite[Proposition 3.1]{HShasse}), we identify $\lambda_{U\ast} \alpha_U^c \in \H^4_c(U,\Zbb_{\Xcal/U}(2))$ to an element of $\Hscr^1(\Ccal_r)$. Taking cohomology of \eqref{eq:Descent17} yields a commutative diagram
        \begin{equation} \label{eq:Descent19}
        \xymatrix{
            \H^1_c(U,\Tcal') \ar[r] \ar[d]^{\lambda_{U\ast}} & \H^3_c(U,\Gbb_m \otimes^{\Lbb} \Gbb_m) \ar[d] \\
            \H^4_c(U,\Zbb_{\Xcal/U}(2)) \ar[r] & \H^5_c(U,\Zbb(2)) \cong \Qbb/\Zbb.
        }
    \end{equation}
    In particular, the bottom arrow maps $\lambda_{U \ast} \alpha_U^c$ to $\beta_U^c = -\gamma_3(\gamma_2(\varepsilon'_U)) \sqdot \alpha_U^c$, which is equal to $-\pair{\eta_U,\alpha_U^c}_{\AV}$ by \eqref{eq:Descent15}. On the other hand, $\beta_U^c$ can be identified to an element of $\Hscr^2(\Ccal_l)$. Passing to the direct limit over $U$ smaller and smaller, $\lambda_{U \ast} \alpha_U^c$ becomes $\lambda_\ast \alpha = u(\alpha)$, and $\beta_U^c$ becomes the image of $u(\alpha)$ by the snake lemma construction applied to
        \begin{equation*}
        \xymatrix@C-1pc{
            & \H^4(K,\Zbb(2)) \ar[r] \ar[d] & \H^4(X,\Zbb(2)) \ar[r] \ar[d] & \frac{\H^4(X,\Zbb(2))}{\Img \H^4(K,\Zbb(2))} \ar[r] \ar[d] & 0 \\
            0 \ar[r] & \prod\limits_{v \in \Omega^{(1)}} \H^4(K_v^{\h},\Zbb(2)) \ar[r] & \prod\limits_{v \in \Omega^{(1)}} \H^4(X_v^{\h},\Zbb(2)) \ar[r] & \prod\limits_{v \in \Omega^{(1)}} \frac{\H^4(X_v^{\h},\Zbb(2))}{\Img \H^4(K_v^{\h},\Zbb(2))}, 
        }
    \end{equation*}
    with exact rows (where $X_v^{\h}:=X \times_K K_v^{\h}$). Let us show that we may replace henselizations by completions in the above construction. Indeed, by \cite[Chapter III, Remark 3.11]{MilneEtale}, $\H^3(\Ocal_{\Omega,v}^{\h},\Qbb/\Zbb(2)) \cong \H^3(\Ocal_v,\Qbb/\Zbb(2)) \cong \H^3(k(v),\Qbb/\Zbb(2)) = 0$ since $\cd (k(v)) = 2$. It follows, by the localization sequence in \'etale cohomology, that we have a chain of isomorphisms
        \begin{equation*}
            \H^4(K_v^{\h},\Zbb(2)) \cong \H^3(K_v^{\h},\Qbb/\Zbb(2)) \cong \H^2(k(v),\Qbb/\Zbb(1)) \cong \H^3(K_v,\Qbb/\Zbb(2)) \cong \H^4(K_v,\Zbb(2)).
        \end{equation*}
    It follows that $\beta_U^c$ becomes $\sigma(\beta)$ by taking limit, where $\beta \in \tfrac{\prod_{v \in \Omega^{(1)}} \H^4(K_v,\Zbb(2))}{\Img \H^4(K,\Zbb(2))}$ is the image of $u(\alpha)$ by the snake lemma construction. By Lemma \ref{lm:Descent1}, one has
        \begin{equation*}
            \rho_X(u(\alpha)) = \sigma(\beta) = -\pair{\eta,\alpha}_{\PT},
        \end{equation*}
    which concludes the proof of Proposition \ref{prop:Descent1}.
\end{proof}
\begin{proof}[Proof of Theorem \ref{thm:Descent}\ref{thm:Descent1}]
     If there exists a family $(P_v)_{v \in \Omega^{(1)}}$ orthogonal to $u(\Sha^1(K,T'))$ relative to the pairing \eqref{eq:UnramifiedPairingGersten}, then the $\rho_X(u(\alpha)) = 0$ for all $\alpha \in \Sha^1(K,T')$. By Proposition \ref{prop:Descent1} and the non-degeneracy of $\pair{-,-}_{\PT}$, one has $e_X = 0$. The converse is obvious.
\end{proof}

\subsection{Description of the obstruction using universal torsors} \label{subsec:DescentDescription}

In this paragraph, we prove Theorem \ref{thm:Descent}\ref{thm:Descent2}. When universal torsors exist, they give an explicit description of the map $u$ in \eqref{eq:MapU} as in the following analogue of \cite[Lemme 3.5.2]{CTSdescent} (see also \cite[Lemma 3]{SkorobogatovBeyond}).

\begin{prop} \label{prop:Descent2}
    Let $\pi: X \to \Spec K$ be a smooth proper geometrically integral variety such that the abelian group $\Pic\ol{X}$ is finitely generated and free. Let $T$ be the N\'eron--Severi torus of $X$. Suppose that $Y \to X$ is a universal torsor\footnote{By our convention, this means its type is the identity of $\wh{T}$.} under $T$. Then the map $u$ constructed in \eqref{eq:MapU} is equal to the composite
        \begin{equation*}
            \H^1(K,T') \xrightarrow{[Y] \cup \pi^\ast(-)} \H^4(X,\Zbb(2)) \to \tfrac{\H^4(X,\Zbb(2))}{\Img \H^4(K,\Zbb(2))},
        \end{equation*}
    where the cup product $\H^1(X,T) \times \H^1(X,T') \xrightarrow{\cup} \H^4(X,\Zbb(2))$ is induced by the pairing \eqref{eq:PairingTorus}.
\end{prop}
\begin{proof}
    Let $\alpha \in \H^1(K,T') = \Ext^1_K(\Zbb,T')$, which can be represented by a morphism $\Zbb \to T'[1]$ in $\Dcal^+(K)$. This morphism gives rise to the vertical arrows in the following commutative diagram in $\Dcal^+(\Ab)$:
        \begin{equation} \label{eq:Descent21}
            \xymatrix{
               \Rbb\Hom_K(T',\Rbb\pi_\ast\Zbb_X(2)) \ar[r] \ar[d] & \Rbb\Hom_K(T',\Zbb_{X/K}(2)) \ar[d] \\
                \Hbb(K,\Rbb\pi_\ast\Zbb_X(2))[1] \ar[r] & \Hbb(K,\Zbb_{X/K}(2))[1].
            }
        \end{equation}
    The horizontal arrows in \eqref{eq:Descent21} are induced by the map $\Rbb\pi_\ast\Zbb(2) \to \Zbb_{X/K}(2)$. In what follows, we shall make use of the fact that $\Rbb\Hom_X(T', -) = \Rbb\Hom_K(T', -) \circ \Rbb \pi_\ast$ and $\Hbb(X,-) = \Hbb(K,-) \circ \Rbb \pi_\ast$ (see \cite[Corollary 10.8.3]{Weibel}).
    We claim that there is a commutative diagram
        \begin{equation} \label{eq:Descent22}
            \xymatrix@C-1pc{
                \H^1(X,T) \ar@{=}[r] \ar[d]^{\cong} & \H^1(X,T) \ar[r]^-{\type} \ar[d]^{\cong} & \Hom_K(\wh{T},\Pic\ol{X}) \ar@{=}[d] \\
                \Ext^1_X(\wh{T},\Gbb_m) \ar[d]^{-\otimes^{\Lbb} \Gbb_m} & \Ext^1_K(\wh{T},\tau_{\le 1} \Rbb \pi_\ast \Gbb_{m,X}) \ar[r] \ar[l]_-{\cong} \ar[d]^{-\otimes^{\Lbb} \Gbb_m} & \Hom_K(\wh{T},\wh{T}) \ar[d]^{-\otimes^{\Lbb} \Gbb_m} \\
                \Ext^1_X(T',\Gbb_m \otimes^{\Lbb} \Gbb_m) \ar[d]^{\theta_{1\ast}} & \Ext^1_K(T',(\tau_{\le 1} \Rbb \pi_\ast \Gbb_{m,X}) \otimes^{\Lbb} \Gbb_m) \ar[l]_-{\theta_{2\ast}} \ar[r] \ar[d]^{\theta_\ast} & \Hom_K(T',T') \ar[d]^{\lambda_\ast} \\
                \Ext_X^3(T',\Zbb(2)) \ar@{=}[r] \ar[d]^{\sqdot \, \pi^\ast \alpha} & \Ext_K^3(T',\Rbb\pi_\ast\Zbb(2)) \ar[r] \ar[d]^{\sqdot \, \alpha} 
                 & \Ext_K^3(T',\Zbb_{X/K}(2)) \ar[d]^{\sqdot \, \alpha} \\
                \H^4(X,\Zbb(2)) \ar@{=}[r] & \H^4(K,\Rbb \pi_\ast \Zbb(2)) \ar[r] & \H^4(K,\Zbb_{X/K}(2)),
            }
        \end{equation}
         where $\sqdot$ means the Yoneda product, and where the maps $\theta_1,\theta_2$ were defined in the course of constructing the map $\theta$ from \eqref{eq:MapTheta}. To see this, let us consider the four rectangles of \eqref{eq:Descent22} from the top to the bottom. As for the first rectangle, the left bottom horizontal arrow is induced by the natural map $\tau_{\le 1} \Rbb \pi_\ast \Gbb_{m,X} \to \Rbb \pi_\ast \Gbb_{m,X}$ (keeping in mind that $\Ext^1_K(\wh{T},\Rbb \pi_\ast \Gbb_{m,X}) = \Ext^1_X(\wh{T},\Gbb_m)$), and the right bottom horizontal arrow is induced by the map from triangle \eqref{eq:TriangleHatT}. The commutativity of this rectangle and the established isomorphisms are well-known; see for example \cite[Proof of Proposition 8.1, Appendix B]{HSdescent}. The second rectangle commutes by the functoriality of $- \otimes^{\Lbb} -$ (bearing in mind that $\Ext^1_K(T',\Rbb \pi_\ast(\Gbb_{m,X} \otimes^{\Lbb} \Gbb_{m,X})) = \Ext^1_X(T',\Gbb_m \otimes^{\Lbb} \Gbb_m)$). As for the third rectangle, the left square obviously commutes (noting that, of course, $\Ext^3_K(T',\Rbb\pi_\ast \Zbb(2)) = \Ext^3_X(T',\Zbb(2))$), and the right square is induced by diagram \eqref{eq:MapU2}. As for the fourth rectangle, the left square obviously commutes, and the right square is obtained by taking cohomology of \eqref{eq:Descent21}.

         Let $Y \to X$ be a universal torsor. By our convention, its type is the identity of $\wh{T}$. Hence, the image of $[Y] \in \H^1(X,T)$ in $\H^4(K,\Zbb_{X/K}(2))$ by \eqref{eq:Descent22} is precisely $\lambda_\ast \alpha$. Under the identification \eqref{eq:MapU3}, this is the same as $u(\alpha)$. Thus, in order to prove Proposition \ref{prop:Descent2}, it remains to show that the image of $[Y]$ in $\H^4(X,\Zbb(2))$ by \eqref{eq:Descent22} is precisely $[Y] \cup \pi^\ast \alpha$. To this end, we argue as in the proof of Proposition \ref{prop:Descent1} to obtain a commutative diagram
        \begin{equation} \label{eq:Descent23}
            \xymatrix@C-1pc{
                & \H^1(X,\iHom_X(\wh{T},\Gbb_m)) \ar[r]^-{\cong} & \Ext^1_X(\wh{T},\Gbb_m) \ar@{=}[rd] \\
                \H^1(X,T) \ar[r] \ar[ru]^{\cong} \ar[rd] & \H^1(X,\iHom_X(\wh{T},T \otimes^{\Lbb} \wh{T})) \ar[u] \ar[d] \ar[r] & \Ext^1_X(\wh{T},T \otimes^{\Lbb} \wh{T}) \ar[d]^{-\otimes^{\Lbb} \Gbb_m} \ar[u]^{\gamma_1} \ar[r]^-{\gamma_1} & \Ext^1_X(\wh{T},\Gbb_m) \ar[d]^{-\otimes^{\Lbb} \Gbb_m} \\
                & \H^1(X,\iHom_X(T',T \otimes^{\Lbb} T')) \ar[r] & \Ext^1_X(T', T \otimes^{\Lbb} T') \ar[r]^-{\gamma_2} & \Ext^1_X(T',\Gbb_m \otimes^{\Lbb} \Gbb_m) \ar[d]^{\theta_{1\ast}} \\
                &&& \Ext^3_X(T',\Zbb(2)),
            }
        \end{equation}    
         similar to \eqref{eq:Descent12}. Denote by $\varepsilon$ the image of $[Y]$ in $\Ext^1_X(T',T \otimes^{\Lbb} T')$ by \eqref{eq:Descent23}. Then the image of $[Y]$ in $\Ext^1_X(T',\Gbb_m \otimes^{\Lbb} \Gbb_m)$ by \eqref{eq:Descent22} is precisely $\gamma_2(\varepsilon)$. Now, we have a commutative diagram of pairings
        \begin{equation} \label{eq:Descent24}
            \xymatrix@C-2pc{
                \H^1(X,T) \ar[d]  & \times & \H^1(X,T') \ar@{=}[d] \ar[rrrr]^-{\cup} &&&& \H^2(X,T \otimes^{\Lbb} T') \ar@{=}[d]  \\
                \Ext^1_X(T',T \otimes^{\Lbb} T') \ar[d]^{\gamma_2} & \times & \H^1(X,T') \ar@{=}[d] \ar[rrrr]^-{\sqdot} &&&& \H^2(X,T \otimes^{\Lbb} T') \ar[d] \\
                \Ext^1_X(T',\Gbb_m \otimes^{\Lbb} \Gbb_m) \ar[d]^{\theta_{1\ast}} & \times & \H^1(X,T') \ar@{=}[d] \ar[rrrr]^-{\sqdot} &&&& \H^2(X,\Gbb_m \otimes^{\Lbb} \Gbb_m) \ar[d] \\
                \Ext^3_X(T',\Zbb(2)) & \times & \H^1(X,T') \ar[rrrr]^-{\sqdot} &&&& \H^4(X,\Zbb(2)),
            }
        \end{equation}    
        similar to \eqref{eq:Descent14} (for the commutativity of the top square, see \cite[Chapter V, Proposition 1.20]{MilneEtale}). This yields the identity $\theta_{1\ast} \gamma_2(\varepsilon) \sqdot \pi^\ast \alpha = [Y] \cup \pi^\ast \alpha \in \H^4(X,\Zbb(2))$, which is exactly what we need. Proposition \ref{prop:Descent2} is hence proved. 
\end{proof}

\begin{proof} [Proof of Theorem \ref{thm:Descent}\ref{thm:Descent2}]
    We start with the inclusion ``$\subseteq$''. Suppose that there exists a family $(P_v)_{v \in \Omega^{(1)}}$ orthogonal to $\Img(u)$ relative to the pairing \eqref{eq:UnramifiedPairingGersten}. By \ref{thm:Descent1}, we know that there exists a universal torsor $f: Y \to X$. In the light of Proposition \ref{prop:Descent2}, we have
        \begin{equation} \label{eq:Descent25}
            \sum_{v \in \Omega^{(1)}} [Y](P_v) \cup \loc_v(\alpha) = \sum_{v \in \Omega^{(1)}}([Y] \cup \pi^\ast \alpha)(P_v) = 0 \in \Qbb/\Zbb
        \end{equation}
    for all $\alpha \in \H^1(K,T')$. Note that if $\Xcal \to U$ is a proper integral model of $X$ over some non-empty open subset $U \subseteq \Omega$, then $X(K_v) = \Xcal(\Ocal_v)$ for all $v \in U^{(1)}$ by the valuative criterion for properness. Furthermore, shrinking $U$ if necessary, we may assume that $Y$ extends to a torsor $\Ycal \to \Xcal$ under a $U$-torus $\Tcal$ extending $T$. Thus, $[Y](P_v)$ comes from $[\Ycal](P_v) \in \H^1(\Ocal_v,\Tcal)$ for all $v \in U^{(1)}$, or $([Y](P_v))_{v \in \Omega^{(1)}} \in \Pbb^1(K,T)$. By virtue of \eqref{eq:Descent25} and the exact sequence \eqref{eq:PoitouTateSequenceTorus}, there exists $t \in \H^1(K,T)$ such that $\loc_v(t) = [Y](P_v)$ for all $v \in \Omega^{(1)}$. Twisting by a Galois cocycle representing $t$ yields a torsor $\tensor[_t]{f}{}: \tensor[_t]{}{} Y \to X$ (that is, $[\tensor[_t]{}{} Y] = [Y] - \pi^\ast t \in \H^1(X,T)$) such that $P_v \in \tensor[_t]{f}{}(\tensor[_t]{}{} Y (K_v))$ (see \eqref{eq:Twisting}) for all $v \in \Omega^{(1)}$. The torsor $\tensor[_t]{}{} Y$ is again universal by the fundamental exact sequence \eqref{eq:FundamentalExactSequence}. This proves the inclusion ``$\subseteq$''.

    Conversely, if $f: Y \to X$ is a universal torsor such that $(P_v)_{v \in \Omega^{(1)}} \in f \left(\prod_{v \in \Omega^{(1)}} Y(K_v)\right)$, then $[Y](P_v) = 0 \in \H^1(K_v,T)$ for all $v \in \Omega^{(1)}$. This obviously implies the identity \eqref{eq:Descent25}, which means $(P_v)_{v \in \Omega^{(1)}}$ orthogonal to $\Img(u)$ by Proposition \ref{prop:Descent2}. This proves the inclusion ``$\supseteq$''.
\end{proof}

\subsection{End of the proof of Theorem \ref{customthm:Descent}} \label{subsec:DescentEnd}

In this paragraph, we finish the proof of Theorem \ref{thm:Descent} ({\em i.e.} Theorem \ref{customthm:Descent}).

\begin{proof} [Proof of Theorem \ref{thm:Descent}\ref{thm:Descent3}]
    Suppose that there is a universal torsor $Y \to X$ (otherwise, there would be nothing to prove). In view of the fundamental exact sequence \eqref{eq:FundamentalExactSequence}\footnote{Recall that a torsor is universal if its type is the identity of $\wh{T}$.}, we have to show that there are only finitely many classes $t \in \H^1(K,T)$ for which $\prod_{v \in \Omega^{(1)}}\tensor[_t]{}{}Y(K_v) \neq \varnothing$. Equivalently, by \eqref{eq:Twisting}, we have to show that the property
        \begin{equation} \label{eq:Descent31}
            \text{``there exists $(P_v)_{v \in \Omega^{(1)}} \in \prod_{v \in \Omega^{(1)}} X(K_v)$ with $[Y](P_v) = \loc_v(t)$ for all $v \in \Omega^{(1)}$''}
        \end{equation}
    holds for only finitely many classes $t \in \H^1(K,T)$. Let $\Xcal \to U$ be a proper integral model of $X$ over some non-empty open subset $U \subseteq \Omega$. Shrinking $U$ if necessary, we may assume that $T$ extends to a $U$-torus $\Tcal$ and $Y$ extends to a torsor $\Ycal \to \Xcal$ under $\Tcal$. Suppose that $t \in \H^1(K,T)$ satisfies \eqref{eq:Descent31}. For all $v \in U^{(1)}$, since $X(K_v) = \Xcal(\Ocal_v)$ by the valuative criterion for properness, the class $\loc_v(t)$ comes from $[\Ycal](P_v) \in \H^1(\Ocal_v,\Tcal)$. On the other hand, we have an exact sequence
        \begin{equation*}
            \H^1(U,\Tcal) \to \prod_{v \notin U} \H^1(K_v,T) \times \prod_{v \in  U^{(1)}} \H^1(\Ocal_v,\Tcal) \to \H^1(K,T')^D,
        \end{equation*}
    obtained in the course of establishing the exact sequence \eqref{eq:PoitouTateSequenceTorus} (see \cite[Proof of Proposition 3.5]{HSS}). Since $\H^1(K,T)$ is orthogonal to $\H^1(K,T')^D$ by the generalized Weil reciprocity law \eqref{eq:WeilReciprocity}, we see that $t$ comes from $\H^1(U,\Tcal)$. Hence, it suffices to show that image of the map $\H^1(U,\Tcal) \to \H^1(K,T)$ is finite. Indeed, if $n$ is the degree of a field extension splitting $T$, then $n\H^1(K,T) = 0$ by Hilbert's Theorem 90 and restriction-corestriction. Thus, the map $\H^1(U,\Tcal) \to \H^1(K,T)$ factors through $\H^1(U,\Tcal)/n$. But this latter injects into the group $\H^2(U,\tensor[_n]{}{} \Tcal)$, which is finite (combine \cite[Chapter VI, Corollary 2.8]{MilneEtale} with the localization sequence \eqref{eq:LocalizationFinite}, nothing that each group $\H^2(K_v,\tensor[_n]{}{} T)$ is finite, and that $\Omega \setminus U$ is itself finite).
\end{proof}

\begin{proof} [Proof of Theorem \ref{thm:Descent}\ref{thm:Descent4}]
    Suppose that the universal $X$-torsors under $T$ satisfy the local-global principle. If there exists a family $(P_v)_{v \in \Omega^{(1)}} \in \prod_{v \in \Omega^{(1)}} X(K_v)$ orthogonal to $\Img(u)$ relative to the pairing \eqref{eq:UnramifiedPairingGersten}, then it follows from \ref{thm:Descent2} that there exists a universal torsor $f: Y \to X$ such that $\prod_{v \in \Omega^{(1)}} Y(K_v) \neq \varnothing$. By our assumption that $Y$ satisfies the local-global principle, one has $Y(K) \neq \varnothing$, {\em a fortiori} $X(K) \neq \varnothing$. 

    Suppose that the universal $X$-torsors under $T$ satisfy the local-global principle and weak approximation. Let $(P_v)_{v \in \Omega^{(1)}} \in \prod_{v \in \Omega^{(1)}} X(K_v)$ be a family orthogonal to $\Img(u)$, $S \subseteq \Omega^{(1)}$ a finite set of closed points, and $\Uscr_v \subseteq X(K_v)$ a $v$-adic neighborhood of $P_v$ for each $v \in S$. By \ref{thm:Descent2}, there exists a universal torsor $f: Y \to X$ and a family $(Q_v)_{v \in \Omega^{(1)}} \in \prod_{v \in \Omega^{(1)}} Y(K_v)$ such that $f(Q_v) = P_v$ for all $v \in \Omega^{(1)}$. By our assumption on $Y$, there exists a point $Q \in Y(K)$ which belongs to $\prod_{v \in S} f^{-1}(\Uscr_v) \times \prod_{v \notin S} Y(K_v)$. Then the point $f(Q) \in X(K)$ belongs to $\prod_{v \in S} \Uscr_v \times \prod_{v \notin S} X(K_v)$. 
\end{proof}

We conclude this section with the following interesting

\begin{rmk} \label{rmk:ImageOfU}
    In fact, the image of the map $u$ from \eqref{eq:MapU} is contained in 
        \begin{equation*}
            \Img\left(\tfrac{\H^3(X,\Qbb/\Zbb(2))}{\Img\H^3(K,\Qbb/\Zbb(2))} \to \tfrac{\H^4(X,\Zbb(2))}{\Img\H^4(K,\Zbb(2))}\right).
        \end{equation*}
    Indeed, since the group $\H^1(K,T')$ has finite exponent by Hilbert's Theorem 90 and restriction-corestriction, it would be enough to show that the torsion elements of $\frac{\H^4(X,\Zbb(2))}{\Img\H^4(K,\Zbb(2))}$ comes from $\frac{\H^3(X,\Qbb/\Zbb(2))}{\Img\H^3(K,\Qbb/\Zbb(2))}$. To see this, let $A \in \H^4(X,\Zbb(2))$ and $n \ge 1$ such that $nA = \pi^\ast c$ for some $c \in \H^4(K,\Zbb(2))$, where $\pi: X \to \Spec K$ denotes the structure morphism. One has $\H^4(K,\Zbb/n(2)) = \H^4(K,\mu_n^{\otimes 2}) = 0$ because $\cd(K) \le 3$, hence $c = nc'$ for some $c' \in \H^4(K,\Zbb(2))$. It follows that the element $A - \pi^\ast c' \in \H^4(X,\Zbb(2))$ is $n$-torsion, hence it comes from $\H^3(X,\Zbb/n(2))$, {\em a fortiori} from $\H^3(X,\Qbb/\Zbb(2))$. Since $\H^4(K,\Zbb(2)) \cong \H^3(K,\Qbb/\Zbb(2))$, we conclude that $A$ itself comes from $\H^3(X,\Qbb/\Zbb(2))$.

    Consequently, the statements \ref{thm:Descent1}, \ref{thm:Descent2}, and \ref{thm:Descent4} of Theorem \ref{thm:Descent} can be refined as follows
        \begin{enumerate}
            \item Universal $X$-torsors exist if and only if the map $\rho_X$ from \eqref{eq:AdelicReciprocity} is the zero map.

            \item A family $(P_v)_{v \in \Omega^{(1)}} \in \prod_{v \in \Omega^{(1)}} X(K_v)$ can be lifted to a universal torsor $f: Y \to X$ if it is orthogonal to $\H^3(X,\Qbb/\Zbb(2))$.

            \setcounter{enumi}{3}

            \item If the universal torsors $Y \to X$ satisfy the local-global principle (resp. the local-global principle and weak approximation), then the obstruction to the local-global principle (resp. the local-global principle and weak approximation) on $X$, defined by the pairing \eqref{eq:AdelicPairing}, is the only one.
        \end{enumerate}
\end{rmk}
\section{Local-global principle and weak approximation} \label{sec:PAdic}

This section is devoted to the proof of Theorems \ref{customthm:Hasse} and \ref{customthm:Weak}. For each of these results, we offer two proofs. The first proofs invoke the results from section \ref{sec:Descent} (Theorem \ref{thm:Descent}\ref{thm:Descent1} for the local-global principle and Proposition \ref{prop:Descent2} for the weak approximation). They are presented in paragraphs \ref{subsec:Hasse} and \ref{subsec:Weak} respectively. The second proofs, which use the fibration methods, rely on an observation communicated to the author by Jean-Louis Colliot-Th\'el\`ene. They shall be presented in paragraph \ref{subsec:Modified}. We also discuss some questions related to weak approximation in  paragraph \ref{subsec:Example}. In particular, we show that any finite abelian group is a Galois group over any $p$-adic function field, rediscovering the positive answer to the abelian case of the inverse Galois problem over $\Qbb_p(t)$.

\subsection{Local-global principle for stabilizers of type $\umult$} \label{subsec:Hasse}

We establish Theorem \ref{customthm:Hasse} in this paragraph. First, recall some facts.

As a rule, the problem of the existence of rational points on homogeneous spaces is harder than that of weak approximation. It requires the general machinery of {\em liens} (or bands, kernels) and non-abelian Galois cohomology, which has been systematically studied in the last 30 years \cite{BorovoiSecond,DLA}. We refer to \cite[\S 1]{FSS} for a complete exposition. Let $X$ be a homogeneous space of a smooth algebraic group $G$ over a field $K$. Let $\ol{H}$ denote the stabilizer of a $\ol{K}$-point of $X$, which is supposed to be smooth. If $\ol{H}$ is commutative, it has natural $K$-form $H$. Otherwise, $\ol{H}$ need not be defined over $K$. Nevertheless, one can always define the associated {\em Springer $K$-lien} $L_X$ ({\em grosso modo}, it is the $\ol{K}$-group $\ol{H}$ equipped with a natural {\em outer Galois action}, {\em i.e.} a Galois action modulo conjugation), the set $\H^2(K,L_X)$ of non-abelian Galois $2$-cohomology, and the {\em Springer class} $\eta_X \in \H^2(K,L_X)$. The class $\eta_X$ is {\em neutral} if and only if $X$ is dominated by a principal homogeneous space of $G$ (if $\H^1(K,G) = 1$, for example, when $G$ is special, this is equivalent to $X(K) \neq \varnothing$). 

Since the derived subgroup $[\ol{H},\ol{H}]$ is characteristic in $\ol{H}$, the canonical outer Galois action induces an {\em action} on $\ol{H}^{\ab} = \ol{H}/[\ol{H},\ol{H}]$. Thus we obtain a $K$-form $H^{\ab}$ of $\ol{H}^{\ab}$. Since every character of $\ol{H}$ factors through $\ol{H}^{\ab}$, the group $\Hom_{\ol{K}}(\ol{H},\Gbb_m) = \Hom_{\ol{K}}(\ol{H}^{\ab},\Gbb_m)$ is equipped with a structure of $\Gamma_K$-module {\em via} this $K$-form of $\ol{H}^{\ab}$.

When $H$ is a $K$-group, we denote by $\lien(H)$ the canonical $K$-lien associated with $H$. If $H$ is abelian, $\H^2(K,\lien(H))$ is just the usual Galois cohomology group $\H^2(K,H)$, and its only neutral class is $0$. Finally, a morphism $L \to L'$ of algebraic $K$-liens induces a {\em relation} $\H^2(K,L) \multimap \H^2(K,L')$. This turns out to be a {\em map} if either the underlying $\ol{K}$-group of $L'$ is commutative or the underlying morphism between $\ol{K}$-groups is surjective.

The following description of the Picard groups of homogeneous spaces is due to Popov \cite[Corollary to Theorem 4]{Popov}, see also  \cite{BDH} and \cite[Theorem 5.8]{BvHtrans}.

\begin{lm} \label{lm:PicardGroupOfHomogeneousSpace}
    Let $X$ be a homogeneous space a smooth, simply connected semisimple linear algebraic group $G$ over a field $K$, with smooth geometric stabilizers $\ol{H}$. Then, as a $\Gamma_K$-module, $\Pic\ol{X}$ is isomorphic to the character group of $\ol{H}$ (the Galois action on this group was defined above). The isomorphism is given by pushing forward the class $[\ol{G}] \in \H^0(K,\H^1(\ol{X},\ol{H}))$ of the torsor $\ol{G} \to \ol{X}$ under $\ol{H}$. 
\end{lm}
\begin{proof}
    Here we used the fact that $\ol{K}[X]^\times = \ol{K}^\times$ because $\ol{K}[G]^\times = \ol{K}^\times$ by Rosenlicht's lemma \cite[Proposition 3]{Rosenlicht}, and $\Pic \ol{G} = 0$ since $G$ is simply connected semisimple \cite[\S 4.3]{Voskresenskii}.
\end{proof}

Let $K$ be the function field of a smooth projective geometrically integral curve $\Omega$ over a $p$-adic field $k$, and let $X$ be a homogeneous space of a simply connected semisimple linear algebraic group $G$ over $K$, with geometric stabilizers $\ol{H}$ of type $\umult$, hence an extension of a group $\ol{M}$ of multiplicative type by a unipotent group $\ol{U}$. Since $\ol{U}$ (the unipotent radical of $\ol{H}$) is characteristic in $\ol{H}$, we have a natural Galois action on $\ol{M} = \ol{H}/\ol{U}$ (hence a $K$-form $M$ of $\ol{M}$). Since $\ol{U}$ does not have any non-trivial characters, the character module of $\ol{H}$ is just $\wh{M}$, hence $\Pic \ol{X} = \wh{M}$ by Lemma \ref{lm:PicardGroupOfHomogeneousSpace}. Let $X^c$ be a smooth projective compactification of $X$. Since $X^c$ is smooth, projective and geometrically unirational, the abelian group $\Pic\ol{X}^c$ is finitely generated and free (see the discussion preceding Theorem \ref{thm:Descent}). There is an exact sequence
    \begin{equation*}
        0 \to \Div_{\infty} \ol{X}^c \to \Pic \ol{X}^c \to \Pic \ol{X} \to 0,
    \end{equation*}
where $\Div_{\infty} \ol{X}^c$ denotes the group of Weil divisors on $\ol{X}^c$ supported in $\ol{X}^c \setminus \ol{X}$ (it is a permutation $\Gamma_K$-module). Note that the injectivity of $\Div_{\infty} \ol{X}^c \to \Pic \ol{X}^c$ follows from the fact that $\ol{K}[X]^\times = \ol{K}[G]^\times = \ol{K}^\times$ \cite[Proposition 3]{Rosenlicht}. Let $T$ (resp. $Q$) be the $K$-torus with character module $\Pic \ol{X}^c$ (resp. $\Div_{\infty} \ol{X}^c$). Then $Q$ is quasi-split. We have exact sequences
    \begin{equation} \label{eq:QTM}
        0 \to \wh{Q} \to \wh{T} \to \wh{M} \to 0
    \end{equation}
and
    \begin{equation} \label{eq:MTQ}
        1 \to M \to T \to Q \to 1.
    \end{equation}
Let $M' = \wh{M} \otimes^{\Lbb} \Zbb(1)$. Applying the functor $- \otimes^{\Lbb} \Zbb(1)$ to \eqref{eq:QTM} yields a distinguished triangle
    \begin{equation} \label{eq:MQT}
        M' \to Q' \to T' \to M'[1]
    \end{equation}
in $\Dcal^+(K)$. In particular, $M$ (resp. $M'$) is quasi-isomorphic to the complex $[T \to Q]$ (resp. $[Q' \to T']$) concentrated in degrees $0$ and $1$.

\begin{rmk} \label{rmk:FiniteCase}
    If $M = F$ is finite abelian, the map $\wc{T} \to \wc{Q}$ on cocharacter modules is injective. Hence there is an exact sequence
        \begin{equation*}
            1 \to F' \to Q' \to T' \to 1,
        \end{equation*}
    where $F' = \iHom_K(F,\Qbb/\Zbb(2))$. In this case, we have a quasi-isomorphism $M' \cong F'$.
\end{rmk}

We construct a map
    \begin{equation} \label{eq:MapTau}
        \tau: \Sha^2_\omega(K,M') \to 
        \tfrac{\H^3_{\nr}(K(X)/K,\Qbb/\Zbb(2))}{\Img\H^3(K,\Qbb/\Zbb(2))}
    \end{equation}
as follows. First, by the Gersten resolution \eqref{eq:GerstenResolution}, we have a map
    \begin{equation} \label{eq:MapTau1}
        \tfrac{\H^4(X^c,\Zbb(2))}{\Img \H^4(K,\Zbb(2))} \to \tfrac{\H^3_{\nr}(K(X)/K,\Qbb/\Zbb(2))}{{\Img \H^3(K,\Qbb/\Zbb(2))}}.
    \end{equation}
Furthermore, since the torus $Q'$ is quasi-split, one has $\H^1(L,Q') = 0$ for any overfield $L/K$ and $\Sha^2_\omega(K,Q') = 0$ by Lemma \ref{lm:QuasiTrivialTorusPAdic}. The long exact sequence associated with \eqref{eq:MQT} gives 
    \begin{equation} \label{eq:MapTau2}
        \Sha^2_\omega(K,M') \cong \Sha^1_\omega(K,T') \subseteq \H^1(K,T')
    \end{equation}
More generally, one has
    \begin{equation} \label{eq:Sha2M}
        \Sha^2_S(K,M') \cong \Sha^1_S(K,T')
    \end{equation}
for any finite set $S \subseteq \Omega^{(1)}$. We define the map $\tau$ in \eqref{eq:MapTau} as the composite of \eqref{eq:MapTau1}, the map $u: \H^1(K,T') \to \tfrac{\H^4(X^c,\Zbb(2))}{\Img \H^4(K,\Zbb(2))}$ constructed in \eqref{eq:MapU}, and \eqref{eq:MapTau2}. This map $\tau$ shall serve as an obstruction to the local-global principle and weak approximation for $X$.

\begin{rmk}
    According to Remark \ref{rmk:ImageOfU}, the image of the map $\tau$ from \eqref{eq:MapTau} is contained in $\Img\left(\tfrac{\H^3(X^c,\Qbb/\Zbb(2))}{\Img \H^3(K,\Qbb/\Zbb(2))} \to \tfrac{\H^3_{\nr}(K(X)/K,\Qbb/\Zbb(2))}{{\Img \H^3(K,\Qbb/\Zbb(2))}}\right)$.
\end{rmk}

We can now state

\begin{thm} [Theorem \ref{customthm:Hasse}] \label{thm:Hasse}
    Let $K$ be the function field of a smooth projective geometrically integral curve $\Omega$ over a $p$-adic field $k$, and $X$ a homogeneous space of a special, simply connected semisimple algebraic group $G$ over $K$, with geometric stabilizers $\ol{H}$ of type $\umult$. We keep the above notations; in particular, there is a map $\tau$ as in \eqref{eq:MapTau}. If there exists a family $(P_v)_{v \in \Omega^{(1)}} \in \prod_{v \in \Omega^{(1)}} X(K_v)$ orthogonal to $\tau(\Sha^2(K,M'))$ relative to the pairing \eqref{eq:UnramifiedPairing}, then $X(K) \neq \varnothing$. In particular, the unramified first obstruction \eqref{eq:UnramifiedReciprocity} to the local-global principle for $X$ is the only one.
\end{thm}

First, we deal with unipotent stabilizers using the following well-known result. It shall also serve in the proof of Theorem \ref{customthm:Weak}.

\begin{lm} \label{lm:Unipotent}
    Let $K$ be a field of characteristic $0$ and $G$ a special algebraic group over $K$. Then homogeneous spaces of $G$ with unipotent geometric stabilizers have $K$-rational points. They have weak approximation if $G$ does.
\end{lm} 
\begin{proof}
    Let $X$ be such a homogeneous space. The Springer class $\eta_X$ is neutral by \cite[Chapitre IV, Th\'eor\`eme 1.3]{Douai} (see also \cite[Corollary 4.2]{BorovoiSecond}), hence $X$ is dominated by a principal homogeneous space of $G$. Since $G$ is special, this means there exists a $G$-equivariant morphism $\phi: G \to X$. In particular, $X(K) \neq \varnothing$. 
    
    Let $S$ be a finite set of places of $K$, $(P_v)_{v \in S} \in \prod_{v \in S}X(K_v)$, and $\Uscr_v \subseteq X(K_v)$ a neighborhood (for the local topology) of $P_v$, $v \in S$. Each fibre $\phi^{-1}(P_v)$ is a torsor under a unipotent $K_v$-group, hence has a $K_v$-point $Q_v$ by \cite[Lemme 1.13]{Sansuc}. If $G$ has weak approximation, we find a point $Q \in G(K)$ which belongs to $\prod_{v \in S} \phi^{-1}(\Uscr_v)$. Then $\phi(Q) \in X(K)$ belongs to $\prod_{v \in S} \Uscr_v$.
\end{proof}

\begin{proof} [Proof of Theorem \ref{thm:Hasse}]
    Keeps the notations as above. Denote by $\eta \in \H^2(K,M)$ the element corresponding to the elementary obstruction $e_X \in \Ext^2_K(\Pic\ol{X},\Gbb_m) = \Ext^2_K(\wh{M},\Gbb_m)$  under the identification \eqref{eq:ExtOfGroupOfMultiplicativeType} (recall that $\ol{K}[X]^\times = \ol{K}[G]^\times = \ol{K}^\times$ by \cite[Proposition 3]{Rosenlicht}). The projection $\ol{H} \to \ol{M}$ induces a surjective morphism $L_X \to \lien(M)$ of algebraic $K$-liens (recall that $L_X$ denotes the Springer lien of $X$), which in turn induces a {\em map} $\H^2(K,L_X) \to \H^2(K,M)$. 

\begin{lm} \label{lm:Hasse}
    The map $\H^2(K,L_X) \to \H^2(K,M)$ sends $\eta_X$ to $\eta$.
\end{lm}
\begin{proof}
    Actually, this result is valid for any ambient group $G$. It suffices to follow the proof of \cite[Theorem 9.5.1]{SkorobogatovTorsors}. This requires the description of the relation $\H^2(K,L_X) \multimap \H^2(K,M)$ in terms of {\em gerbes} (see for example \cite[\S 2.2]{DLA}). According to \cite[Chapitre IV, \S 3.2, Chapitre V, Propositions 3.1.6 and 3.2.1]{Giraud}, $\eta$ is represented by the gerbe $\Gcal$ of universal $X$-torsors under $M$, that is, for every finite extension $L/K$, the fibre category $\Gcal(L)$ is the groupoid of universal $X_L$-torsors under $M$. On the other hand, $\eta_X$ is represented by the gerbe $\Gcal_X$, whose fibre category $\Gcal_X(L)$ is for every finite extension $L/K$ the groupoid of $L$-torsors under $G$ dominating $X_L$ \cite[Chapitre IV, \S 5.1]{Giraud}. If $Y$ is such an $L$-torsor (equipped with a $G_L$-equivariant dominant morphism $Y \to X_L)$,  let $H_L = \Aut_{G_L}(Y/X_L)$. Then $H_L$ is an algebraic subgroup of $G_L$, and $Y \to X_L$ is a torsor under $H_L$ (see \cite[\S 9.2]{SkorobogatovTorsors}). In particular, $H_L$ is an $L$-form of $\ol{H}$. The contracted product $Z:=Y \times^{H_L}_L M_L$ is an $X_L$-torsor under $M$, and the map $\H^1(\ol{X},\ol{H}) \to \H^1(\ol{X},\ol{M})$ sends $[\ol{G}] = [\ol{Y}]$ to $[\ol{Z}]$. Combining with Lemma \ref{lm:PicardGroupOfHomogeneousSpace}, we see that the identification $\Pic\ol{X} = \wh{M}$ is given by pushing forward the class $[\ol{Z}] \in \H^0(K,\H^1(\ol{X},\ol{M}))$, {\em i.e.} the torsor $Z \to X_L$ has type $\id$. The construction $Y \mapsto Z$ defines a morphism $\Gcal_X \to \Gcal$ of algebraic $K$-gerbes, thus $\H^2(K,L_X) \to \H^2(K,M)$ maps $\eta_X$ to $\eta$.
\end{proof}
    
    Return to the proof of Theorem \ref{thm:Hasse}. Since \eqref{eq:ElementaryObstruction} is functorial in $X$, the map $\Ext^2_K(\wh{F},\Gbb_m) \to \Ext^2_K(\wh{T},\Gbb_m)$ sends $e_X$ to $e_{X^c}$. It follows that the map $\H^2(K,M) \to \H^2(K,T)$ sends $\eta$ to an element $\eta^c$ corresponding to $e_{X^c}$ (under the identification $\H^2(K,T) \cong \Ext^2_K(\wh{T},\Gbb_m)$ of \eqref{eq:ExtOfGroupOfMultiplicativeType}). If $(P_v)_{v \in \Omega^{(1)}} \in \prod_{v \in \Omega^{(1)}} X(K_v)$ is a family orthogonal to $\tau(\Sha^2(K,M'))$, then, as a family in $\prod_{v \in \Omega^{(1)}} X^c(K_v)$, it is orthogonal to $u(\Sha^1(K,T'))$ by the construction of $\tau$ (where $u$ is the map constructed in \eqref{eq:MapU}). Theorem \ref{thm:Descent}\ref{thm:Descent1} then implies that $e_{X^c} = 0$, or $\eta^c = 0$. On the other hand, since $\H^1(K,Q) = 0$ (the torus $Q$ being quasi-split), the long exact sequence associated with \eqref{eq:MTQ} assures that $\H^2(K,M) \to \H^2(K,T)$ is injective. It follows that $\eta = 0$. By Lemma \ref{lm:Hasse}, the map $\H^2(K,L_X) \to \H^2(K,M)$ sends $\eta_X$ to the neutral class $0$. Now, \cite[Theorem 3.4]{DLA} provides a diagram
        \begin{equation*}
            \xymatrix{
                & X_1 \ar[ld]_{\phi} \ar[rd]^{\psi} \\
                X && X_2,
            }
        \end{equation*}
    where
    \begin{itemize}
        \item $X_1$ is a homogeneous space of $G \times_K \SL_n$ with Springer lien $L_{X_1} \cong L_X$ and Springer class $\eta_{X_1} = \eta_X$,
        \item $X_2 = M \backslash \SL_n$ for some $K$-embedding $M \hookrightarrow \SL_n$ and some $n$,  
        \item $\phi$ is a torsor under $\SL_n$,
        \item the fibres of $\psi$ are homogeneous spaces of $G$ with geometric stabilizers $\Ker(\ol{H} \to \ol{M}) = \ol{U}$.
    \end{itemize}
    Since $X_2(K) \neq \varnothing$, we have $X_1(K) \neq \varnothing$ by Lemma \ref{lm:Unipotent}, hence $X(K) \neq \varnothing$.
\end{proof}

\subsection{Weak approximation for stabilizers of type $\umult$} \label{subsec:Weak}

In this paragraph, we establish Theorem \ref{customthm:Weak}. We start by recalling the following well-known result, which already appeared in \cite{Chernousov} (see also \cite{HarariQuelques,LA} for finite subgroups of $\SL_n$). We give a proof here for the sake of reference.

\begin{lm} \label{lm:WAForClassifyingSpaces}
    Let $G$ be a smooth algebraic group over a field $K$, $H$ a smooth Zariski closed subgroup of $G$, and $X = H \backslash G$. The projection $G \to X$ is then a torsor under $H$. For any finite set $S$ of places of $K$, if a family $(P_v)_{v \in S} \in \prod_{v \in S} X(K_v)$ lies in the closure (for the product of local topologies) of the diagonal image of $X(K)$, then $([G](P_v))_{v \in S}$ belongs to the image of the localization $\H^1(K,H) \to \prod_{v \in S} \H^1(K_v,H)$. The converse holds if $G$ is special and has weak approximation.
\end{lm} 
\begin{proof}
    Suppose that $(P_v)_{v \in S}$ lies in the closure of $X(K)$. By Lemma \ref{lm:Continuity}\ref{lm:Continuity3}, for each $v \in S$, there exists a neighborhood $\Uscr_v \subseteq X(K_v)$ of $P_v$ such that $[G](P'_v) = [G](P_v)$ for all $P'_v \in \Uscr_v$. Let $P \in \prod_{v \in S} \Uscr_v$ be a $K$-point. Then the element $[G](P) \in \H^1(K,H)$ satisfies $[G](P_v) = \loc_v([G](P))$ for all $v \in S$.
    
    Conversely, suppose that there exists $h \in \H^1(K,H)$ such that $\loc_v(h) = [G](P_v)$ for all $v \in S$. Since $G$ is special, the evaluation-at-$[G]$ map $X(L) \to \H^1(L,H)$ induces a bijection $X(L) / G(L) \cong \H^1(L,H)$ for any overfield $L/K$. In particular, we may write $h = [G](P)$ for some $P \in X(K)$. For each $v \in S$, let $\Uscr_v \subseteq X(K_v)$ be a neighborhood of $P_v$, and let $\Vscr_v$ denote its preimage by the continuous map
        \begin{equation*}
            G(K_v) \to X(K_v), \quad g_v \mapsto P \cdot g_v
        \end{equation*}
    Since $\loc_v([G](P)) = \loc_v(h) = [G](P_v)$, there exists $g_v \in G(K_v)$ such that $P \cdot g_v = P_v$, hence $\Vscr_v \neq \varnothing$. Under the hypothesis that $G$ has weak approximation, there exists $g \in G(K)$ which belongs to $\prod_{v \in S} \Vscr_v$. Then $P \cdot g \in X(K)$ belongs to $\prod_{v \in S} \Uscr_v$.
\end{proof}

\begin{rmk}
    Taking $G = \SL_n$ in \ref{lm:WAForClassifyingSpaces}, we see that weak approximation for the quotient $H \backslash \SL_n$ is an intrinsic property of the algebraic $K$-group $H$ (independent of the embedding $H \hookrightarrow \SL_n$). Indeed, if $H \hookrightarrow \SL_n$ and $H \hookrightarrow \SL_m$ are two embeddings, the quotient varieties $H \backslash \SL_n$ and $H \backslash \SL_m$ are $K$-stably birational by the ``no-name lemma'' \cite[\S 3.2]{CTSrational}.
\end{rmk}

Let $K$ be the function field of a smooth projective geometrically integral curve $\Omega$ over a $p$-adic field $k$. Let $H$ be a $K$-group of type $\umult$, hence an extension of a $K$-group of multiplicative type $M$ by a unipotent group $U$. Let $X = H \backslash G$ for some embedding $H \hookrightarrow G$ into a simply connected semisimple linear algebraic group $G$ over $K$. By Lemma \ref{lm:PicardGroupOfHomogeneousSpace},  $\Pic \ol{X} = \wh{H} = \wh{M}$ as $\Gamma_K$-modules. As in paragraph \ref{subsec:Hasse}, let $X^c$ be a smooth projective compactification of $X$, $T$ the $K$-torus with $\wh{T} = \Pic\ol{X}^c$, and $Q = T/M$ (it is a quasi-split $K$-torus). Finally, let $M' = M \otimes^{\Lbb} \Zbb(1)$ (it is represented by a $2$-term complex of tori). 

Since $X^c(K) \neq \varnothing$, by the fundamental exact sequence \eqref{eq:FundamentalExactSequence}, there exists a universal torsor $Y^c \to X^c$ under $T$. By Proposition \ref{prop:Descent2}, the map $\tau$ from \eqref{eq:MapTau} is the composite
    \begin{equation} \label{eq:MapTauAlternative}
        \Sha^2_\omega(K,M') \cong \Sha^1_\omega(K,T') \xrightarrow{[Y^c] \cup (-)_{X^c}} \H^4(X^c,\Zbb(2)) \to \tfrac{\H^4(X^c,\Zbb(2))}{\Img \H^4(K,\Zbb(2))} \to \tfrac{\H^3_{\nr}(K(X)/K,\Qbb/\Zbb(2))}{\Img \H^3(K,\Qbb/\Zbb(2))},
    \end{equation}
where the subscript ${}_{X^c}$ denotes pullback along the structure morphism $X^c \to \Spec K$, and the last arrow is induced by \eqref{eq:GerstenResolution}. Using the alternative description \eqref{eq:MapTauAlternative}, we shall prove the following ``functoriality'' property of $\tau$, which is required when we apply the fibration method.

\begin{lm} \label{lm:FunctorialityOfTau}
    The map $\tau$ enjoys the following properties.
    \begin{enumerate}
        \item \label{lm:FunctorialityOfTau1} $\tau$ does not depend on the choice of the universal torsor $Y^c$.
        
        \item \label{lm:FunctorialityOfTau2} $\tau$ does not depend on the choice of the smooth projective compactification $X^c$.
        
        \item \label{lm:FunctorialityOfTau3} Let $G, G_1$ be simply connected semisimple linear algebraic groups over $K$. Let $H, H_1$ be $K$-groups of type $\umult$, equipped with respective embeddings into $G, G_1$, and let $X = H \backslash G$, $X_1 = H_1 \backslash G_1$. Let $M, M_1$ be the respective parts of multiplicative type of $H, H_1$, and let $\tau, \tau_1$ be the respective maps constructed in \eqref{eq:MapTau}. Assume that there exists a morphism $\varphi: G \to G_1$ and a $\varphi$-equivariant dominant morphism $\phi: X \to X_1$. Then there is a commutative diagram
            \begin{equation*}
                \xymatrix{
                \Sha^2_\omega(K,M_1') \ar[r]^-{\tau_1} \ar[d]^{\phi^\ast} & \frac{\H^3_{\nr}(K(X_1)/K,\Qbb/\Zbb(2))}{\Img\H^3(K,\Qbb/\Zbb(2))} \ar[d]^{\phi^\ast} \\
                \Sha^2_\omega(K,M') \ar[r]^-{\tau} & \frac{\H^3_{\nr}(K(X)/K,\Qbb/\Zbb(2))}{\Img\H^3(K,\Qbb/\Zbb(2))}.
                }
            \end{equation*}
    \end{enumerate}
\end{lm}
\begin{proof}
    We prove \ref{lm:FunctorialityOfTau1}. If $Y_1^c \to X^c$ is a second universal torsor, then $[Y_1^c] = [Y^c] + t_{X^c}$ in $\H^1(X^c,T)$ for some $t \in \H^1(K, T)$ by virtue of the fundamental exact sequence \eqref{eq:FundamentalExactSequence}\footnote{Recall our convention: a torsor under $T$ is universal if its type is exactly the identity of $\wh{T}$ (not just an isomorphism).}. This yields a map $\H^1(K,T') \to \H^4(X^c,\Zbb(2))$ given by $\alpha \mapsto [Y^c] \cup \alpha_{X^c} + (t \cup \alpha)_{X^c}$. Since $t \cup \alpha$ is an element of $\H^4(K,\Zbb(2)) \cong \H^3(K,\Qbb/\Zbb(2))$, by the description \eqref{eq:MapTauAlternative}, we see that $\tau$ remains unchanged when $Y$ is replaced by $Y_1$.

    We prove \ref{lm:FunctorialityOfTau2}. Let $X_1^c$ be a second smooth projective compactification of $X$. Let $T_1$ be the $K$-torus with $\wh{T}_1 = \Pic\ol{X}_1^c$, and $Y^c_1 \to X_1^c$ a universal torsor. Then there is an exact sequence 
        \begin{equation*}
            0 \to \wh{Q}_1 \to \wh{T}_1 \to \wh{M} \to 0
        \end{equation*}
    similar to \eqref{eq:QTM}, where $Q_1$ is a quasi-split torus. Thus there are $K$-tori $R,R_1,T_2$ such that $R,R_2$ are quasi-split and $T \times_K R \cong T_1 \times R_1 \cong T_2$, and the diagram
        \begin{equation} \label{eq:FunctorialityOfTau1}
            \xymatrix{
                & \wh{T}_2 \ar[ld] \ar[rd] \\
                \wh{T} \ar[rd] && \wh{T}_1 \ar[ld] \\
                & \wh{M} = \Pic\ol{X}
            }
        \end{equation} 
    commutes.
    Let $Y \to X, Y_1 \to X$ denote the respective restrictions of $Y^c \to X^c$ and $Y_1^c \to X_1^c$, which are torsors whose type are respectively given by the bottom arrows of  \eqref{eq:FunctorialityOfTau1}. Since $X(K) \neq \varnothing$, by the fundamental exact sequence \eqref{eq:FundamentalExactSequence}, there exists a torsor $Y_2 \to X$ under $T_2$ whose type $\wh{T}_2 \to \Pic\ol{X}$ is given by either of the two composites in \eqref{eq:FunctorialityOfTau1}. By the very same sequence, the image of $[Y_2]$ in $\H^1(X,T)$ (resp. in $\H^1(X,T_1)$) has the form $[Y] + t_X$ (resp. $[Y_1] + t_{1,X}$) for some $t \in \H^1(K,T)$ (resp. $t_1 \in \H^1(K,T_1)$). Twisting the $X^c$-torsor $Y^c$ (resp. $Y^c_1$) by $t$ (resp. by $t_1$) yields another universal torsor; this is something allowed by \ref{lm:FunctorialityOfTau1}. Hence, we may assume that the image of $[Y_2]$ in $\H^1(X,T)$ (resp. in $\H^1(X,T_1)$) is precisely $[Y]$ (resp. $[Y_1]$). One then obtains a commutative diagram    
        \begin{equation*}
            \xymatrix{
            \H^1(K,T') \ar[rrr]^-{[Y^c] \cup (-)_{X^c}} \ar[rrd]^-{[Y] \cup (-)_X} &&& \H^4(X^c,\Zbb(2)) \ar[r] \ar[d] & \H^3_{\nr}(K(X)/K,\Qbb/\Zbb(2)) \ar@{_{(}->}[d] \\
            \H^1(K,T_2') \ar[u]_{\cong} \ar[d]^{\cong} \ar[rr]^-{[Y_2] \cup (-)_X} && \H^4(X,\Zbb(2)) \ar[r] & \H^4(K(X),\Zbb(2)) \ar[r]^-{\cong} & \H^3(K(X),\Qbb/\Zbb(2)) \\
            \H^1(K,T_1') \ar[rru]^-{[Y_1] \cup (-)_X} \ar[rrr]^-{[Y_1^c] \cup (-)_{X^c}} &&& \H^4(X_1^c,\Zbb(2)) \ar[r] \ar[u] & \H^3_{\nr}(K(X)/K,\Qbb/\Zbb(2)) \ar@{^{(}->}[u]
            }
        \end{equation*}
    (the two squares on the right-hand side are \eqref{eq:GerstenResolutionDiagram}). By the description \eqref{eq:MapTauAlternative}, this shows that the compactifications $X^c$ and $X^c_1$ yield the same map $\tau$.

     We prove \ref{lm:FunctorialityOfTau3}.  By Nagata's Theorem \cite{Nagata}, there exists a compactification $\phi^c: X^c \to X_1^c$ of $\phi$. Let $T$ (resp. $T_1$) be the $K$-torus with $\wh{T} = \Pic\ol{X}^c$ (resp. $\wh{T}_1 = \Pic \ol{X}_1^c$). Then $\phi^c$ induces a morphism $\psi: T \to T_1$ of $K$-tori. Let $Y^c \to X^c$ (resp. $Y_1^c \to X_1^c$) be a universal torsor under $T$ (resp. under $T_1$). Then both the contracted product $Y^c \times_K^T T_1$ and the pullback $Y_1^c \times_{X_1^c} X^c$ are $X^c$-torsors under $T_1$ of type $\psi^\ast: \wh{T}_1 \to \wh{T}$. Again, by the fundamental exact sequence \eqref{eq:FundamentalExactSequence}, the images in $\H^1(X^c,T_1)$ of $[Y^c] \in \H^1(X^c,T)$ and $[Y_1^c] \in \H^1(X_1^c,T_1)$ differ by a class $t_{1,X^c}$, where $t_1 \in \H^1(K,T_1)$. Twisting the $X^c_1$-torsor $Y^c_1$ by $t_1$ (which does not modify the map $\tau$, thanks to \ref{lm:FunctorialityOfTau1}), we may assume that these images coincide. Cup-products with this common value then yield the oblique arrow in the diagram
        \begin{equation*}
            \xymatrix{
            \H^1(K,T_1') \ar[rr]^-{[Y_1^c] \cup (-)_{X_1^c}} \ar[d] \ar[rrd] && \H^4(X_1^c,\Zbb(2)) \ar[d]^{\phi^\ast} \ar[r] & \H^3(K(X_1)/K,\Qbb/\Zbb(2)) \ar[d]^{\phi^\ast} \\
            \H^1(K,T') \ar[rr]^-{[Y^c] \cup (-)_{X^c}} && \H^4(X^c, \Zbb(2)) \ar[r] & \H^3(K(X)/K,\Qbb/\Zbb(2)).
            }
        \end{equation*}
    Since the two triangles commute, the square on the left-hand side also commutes. The square on the right-hand side commutes thanks to the functoriality of the map \eqref{eq:GerstenResolution}. By the description \eqref{eq:MapTauAlternative}, the commutative diagram in the statement of \ref{lm:FunctorialityOfTau3} is established.
\end{proof}

In the course of establishing duality theorems between $M$ and $M'$, Izquierdo defined a pairing
    \begin{equation} \label{eq:PairingComplexTori}
        M \otimes^{\Lbb} M' \to \Zbb(2)[1]
    \end{equation}
(see \cite[Lemme 4.3]{IzquierdoArxiv}). By its very construction, the pairing \eqref{eq:PairingComplexTori} is functorial and generalizes \eqref{eq:PairingTorus}. In particular, these pairings (for $Q, T$ and $M$) are compatible with respect to the exact sequence \eqref{eq:MTQ} and triangle \eqref{eq:MQT}. 

To establish Theorem \ref{customthm:Weak} in the case of stabilizers of multiplicative type ({\em i.e.} when $H = M$), the first step is to reinterpret the alternative description \eqref{eq:MapTauAlternative} of the map $\tau$.

\begin{lm} \label{lm:MultiplicativeType}
   Suppose that $H = M$ is a $K$-group of multiplicative type, so we have a class $[G] \in \H^1(X, M)$ of the torsor $G \to X$ under $M$. Consider the cup-product
        \begin{equation*}
            \H^1(X,M) \times \H^2(X,M') \to \H^4(X,\Zbb(2))
        \end{equation*}
    induced by the pairing \eqref{eq:PairingComplexTori}. Then (up to a sign) for all $\alpha \in \Sha^2_\omega(K,M')$, the image of the class $[G] \cup \alpha_X \in \H^4(X,\Zbb(2))$ in $\tfrac{\H^4(K(X),\Zbb(2))}{\Img \H^4(K,\Zbb(2))} \cong \tfrac{\H^3(K(X),\Qbb/\Zbb(2))}{\Img \H^3(K,\Qbb/\Zbb(2))}$ coincides with $\tau(\alpha)$. 
\end{lm}
\begin{proof}
    Let $Y^c \to X^c$ be a universal torsor\footnote{By our convention, it means its type is the identity of $\wh{T}$.} and let $Y \to X$ be its restriction to $X$, which is a torsor under $T$ whose type is the map $\wh{T} \to \wh{M} = \Pic\ol{X}$ from \eqref{eq:QTM}. Since the (universal) torsor $G \to X$ under $M$ has type $\id$ by Lemma \ref{lm:PicardGroupOfHomogeneousSpace}, the fundamental exact sequence \eqref{eq:FundamentalExactSequence} assures that the map $\H^1(X,M) \to \H^1(X,T)$ sends $[G]$ to $[Y] + t_X$ for some $t \in \H^1(K,T)$. Twisting $Y^c$ by $t$ (which does not change the map $\tau$, by Lemma \ref{lm:FunctorialityOfTau}\ref{lm:FunctorialityOfTau1}), we may assume that the image of $[G]$ in $\H^1(X,T)$ is $[Y]$. Since the pairing \eqref{eq:PairingComplexTori} is functorial, the diagram
        \begin{equation*}
            \xymatrix{
            \H^1(K,T') \ar[rrr]^-{[Y^c] \cup (-)_{X^c}} \ar[rrd]^-{[Y] \cup (-)_X} \ar[d] &&& \H^4(X^c,\Zbb(2)) \ar[ld] \ar[r] \ar[d] & \H^3_{\nr}(K(X)/K,\Qbb/\Zbb(2)) \ar@{_{(}->}[d] \\
            \H^2(K,M')  \ar[rr]^-{[G] \cup (-)_X} && \H^4(X,\Zbb(2)) \ar[r] & \H^4(K(X),\Zbb(2)) \ar[r]^-{\cong} & \H^3(K(X),\Qbb/\Zbb(2)),
            }
        \end{equation*}
    commutes (except the left triangle, which commutes up to a sign). Indeed, the square on the right-hand side is \eqref{eq:GerstenResolutionDiagram}. By the description \eqref{eq:MapTauAlternative}  of $\tau$, the Lemma is proved.
\end{proof}

The next step is to establish following analogue of the exact sequences \eqref{eq:ExactSequenceFinite} and \eqref{eq:ExactSequenceTorus}.

\begin{lm} 
    For each finite set $S \subseteq \Omega^{(1)}$ of closed points, there is an exact sequence \begin{equation} \label{eq:ExactSequenceComplexTori}
        \H^1(K,M) \to \prod_{v \in S} \H^1(K_v,M) \xrightarrow{\theta}  \Sha^2_S(K,M')^D \to \Sha^2(K,M')^D \to 0,
    \end{equation}
    where the map $\theta$ is defined by $\theta((m_v)_{v \in S})(\alpha) = \sum_{v \in S} (m_v \cup \loc_v(\alpha))$. Here the local cup-products $\H^1(K_v,M) \times \H^2(K_v,M') \to \H^4(K_v,\Zbb(2)) \cong \Qbb/\Zbb$ are induced by the pairing \eqref{eq:PairingComplexTori}. 
\end{lm}
\begin{proof}
    Consider the commutative diagram
    \begin{equation} \label{eq:ExactSequenceComplexTori1}
        \xymatrix{
           Q(K) \ar[r] \ar[d] & \H^1(K,M) \ar[r] \ar[d] & \H^1(K,T) \ar[r] \ar[d] & 0 \\
           \prod_{v \in S} Q(K_v) \ar[r] & \prod_{v \in S} \H^1(K_v,M) \ar[r] \ar[d] &  \prod_{v \in S}\H^1(K_v,T) \ar[r] \ar[d] &  0 \\
           & \Sha^2_S(K,M')^D \ar[d] \ar[r]^{\cong} & \Sha^1_S(K,T')^D \ar[d] & \\
           & \Sha^2(K,M')^D \ar[d] \ar[r]^{\cong} & \Sha^1(K,T')^D \ar[d] & \\
           & 0 & 0 & 
           }
    \end{equation}
    with exact rows (the two top rows are the exact sequences associated with \eqref{eq:MTQ}, noting that $\H^1(L,Q) = 0$ for any overfield $L/K$ since $Q$ is quasi-split). The two bottom horizontal arrows are isomorphisms since $\H^1(L,Q') = 0$ for any overfield $L/K$ and $\Sha^2(K,Q') = \Sha^2_S(K,Q') = 0$ by Lemma \ref{lm:QuasiTrivialTorusPAdic}. That \eqref{eq:ExactSequenceComplexTori1} commutes follows from the functoriality of \eqref{eq:PairingComplexTori}. The right column is exact because it is \eqref{eq:ExactSequenceTorus}. The left vertical arrow has dense image (the torus $Q$ has weak approximation since it is $K$-rational). For each $v \in S$, the map $Q(K_v) \to \H^1(K_v,M)$ is the evaluation induced by the class $[T] \in \H^1(Q,M)$ of the torsor $T \to Q$ under $M$, hence it is locally constant by Lemma \ref{lm:Continuity}\ref{lm:Continuity3}. A diagram chasing then shows that the left column ({\em i.e.} the sequence \eqref{eq:ExactSequenceComplexTori}) is exact. 
\end{proof}

\begin{thm} \label{thm:MultiplicativeType}
    Let $K$ be the function field of a smooth projective geometrically integral curve $\Omega$ over a $p$-adic field $k$, $M$ a $K$-group of multiplicative type, and let $X = M \backslash G$, where $M \hookrightarrow G$ is some embedding into a special, simply connected semisimple algebraic group $G$ over $K$ that has weak approximation. Let $M' = \wh{M} \otimes^{\Lbb} \Zbb(1)$, $\tau$ the map constructed in \eqref{eq:MapTau}, and $S \subseteq \Omega^{(1)}$ a finite set of closed points. Then any family $(P_v)_{v \in S} \in \prod_{v \in S} X(K_v)$ orthogonal to $\tau(\Sha^2_S(K,M'))$ relative to the pairing \eqref{eq:UnramifiedPairingS} lies in the closure of the diagonal image of $X(K)$. Moreover, $X$ has weak approximation in $S$ if and only if $\Sha^2_S(K,M') = \Sha^2(K,M')$.
\end{thm}
\begin{proof}
    Let $(P_v)_{v \in S} \in \prod_{v \in S} X(K_v)$ be a family orthogonal to $\tau(\Sha^2_S(K,M'))$ relative to the pairing \eqref{eq:UnramifiedPairingS}. This means $\sum_{v \in S} A(P_v) = 0$ for any $\alpha \in \Sha^2_S(K,M')$ and any lifting $A \in \H^3_{\nr}(K(X)/K,\Qbb/\Zbb(2))$ of $\tau(\alpha)$. By virtue of Lemma \ref{lm:MultiplicativeType}, this is equivalent to
        \begin{equation*}
            \sum_{v \in S} [G](P_v) \cup \loc_v(\alpha) = \sum_{v \in S} ([G] \cup \alpha_X)(P_v) = 0
        \end{equation*}
    for all $\alpha \in \Sha^2_S(K,M')$. Using the exact sequence \eqref{eq:ExactSequenceComplexTori}, we see that $([G](P_v))_{v \in S}$ lies in the image of the localization $\H^1(K,M) \to \prod_{v \in S} \H^1(K_v,M)$. But then $(P_v)_{v \in S}$ lies in the closure of the diagonal image of $X(K)$ by Lemma \ref{lm:WAForClassifyingSpaces}.
    
    Any family $(P_v)_{v \in S} \in \prod_{v \in S}X(K_v)$ is orthogonal to $\tau(\Sha^2(K,M'))$ since this group consists of everywhere locally constant classes and $X(K) \neq \varnothing$. It follows that $X$ has weak approximation in $S$ whenever $\Sha^2(K,M') = \Sha^2_S(K,M')$. Conversely, suppose that $\Sha^2(K,M') \subsetneq \Sha^2_S(K,M')$. Then the map $\Sha^2_S(K,M')^D \to \Sha^2(K,M')^D$ has non-trivial kernel. Again, exactness of \eqref{eq:ExactSequenceComplexTori} implies the existence of a family $(m_v)_{v \in S} \in \prod_{v \in S}\H^1(K_v,M)$ whose image in $\Sha^2_S(K,M')^D$ is non-zero, {\em i.e.} $(m_v)_{v \in S}$ does not come from $\H^1(K,M)$. But since $G$ is special, for each $v \in S$, there exists $P_v \in X(K_v)$ such that $m_v = [G](P_v)$. By Lemma \ref{lm:WAForClassifyingSpaces}, the family $(P_v)_{v \in S}$ does not lie in the closure of $X(K)$, thus $X$ fails approximation in $S$. This concludes the proof of the theorem.
\end{proof}

Finally, we extend Theorem \ref{thm:MultiplicativeType} to the main theorem of this section by allowing the stabilizers to have a unipotent part. The proof uses the fibration method.

\begin{thm} [Theorem \ref{customthm:Weak}] \label{thm:Weak}
    Let $K$ be the function field of a smooth projective geometrically integral curve $\Omega$ over a $p$-adic field $k$, $H$ a linear algebraic $K$-group extension of a group multiplicative type $M$ by a unipotent group $U$, and $X = H \backslash G$ for some $K$-embedding $H \hookrightarrow G$ into a special, simply connected semisimple linear algebraic group $G$ over $K$ that has weak approximation. Let $M' = \wh{M} \otimes^{\Lbb} \Zbb(1)$, $\tau$ the map constructed in \eqref{eq:MapTau}, and $S \subseteq \Omega^{(1)}$ a finite set of closed points. Then any family $(P_v)_{v \in S} \in \prod_{v \in S} X(K_v)$ orthogonal to $\tau(\Sha^2_S(K,M'))$ relative to the pairing \eqref{eq:UnramifiedPairingS} lies the closure of the diagonal image of $X(K)$. In particular, the reciprocity obstruction to weak approximation for $X$ is the only one. Moreover, $X$ has weak approximation in $S$ if and only if $\Sha^2_S(K,M') = \Sha^2(K,M')$.
\end{thm}
\begin{proof}
    Choose an embedding $M \hookrightarrow \SL_n$ for some $n$. The embedding $H \hookrightarrow G$ and the composite $H \to M \hookrightarrow \SL_n$ yield a diagonal embedding $H \hookrightarrow G \times_K \SL_n$. We have a diagram
        \begin{equation} \label{eq:Weak1}
            \xymatrix{
            & X_1 \ar[ld]_{\phi} \ar[rd]^{\psi} \\
            X && X_2
            }
        \end{equation}
    where $X_1 = H \backslash (G \times_K \SL_n)$, $X_2 = M \backslash \SL_n$, $\phi$ is a torsor under $\SL_n$, and the fibres of $\psi$ are homogeneous spaces of $G$ with geometric stabilizers $\Ker(\ol{H} \to \ol{M}) = \ol{U}$. Let $\tau,\tau_1,\tau_2$ be the respective maps associated with $X, X_1, X_2$ {\em via} the construction \eqref{eq:MapTau}. By Lemma \ref{lm:FunctorialityOfTau}\ref{lm:FunctorialityOfTau3}, diagram \eqref{eq:Weak1} yields a commutative diagram
        \begin{equation} \label{eq:Weak2}
            \xymatrix{
            & \Sha^2_\omega(K,M') \ar[d]^{\tau_1} \ar[ld]_{\tau} \ar[rd]^{\tau_2}
            \\
            \frac{\H^3_{\nr}(K(X)/K,\Qbb/\Zbb(2))}{\Img\H^3(K,\Qbb/\Zbb(2))} \ar[r]^{\phi^\ast} & \frac{\H^3_{\nr}(K(X_1)/K,\Qbb/\Zbb(2))}{\Img\H^3(K,\Qbb/\Zbb(2))} 
             & \frac{\H^3_{\nr}(K(X_2)/K,\Qbb/\Zbb(2))}{\Img\H^3(K,\Qbb/\Zbb(2))} \ar[l]_{\psi^\ast}. 
            }
        \end{equation}
    Since $\H^1(K(X),\SL_n) = 1$ by a variant of Hilbert's Theorem 90, the generic fibre of $\phi$ has a section. The extension $K(X_1)/K(X)$ is thus purely transcendental, hence $\phi^\ast$ is an isomorphism.
    
    Let $(P_v)_{v \in S} \in \prod_{v \in S} X(K_v)$ be orthogonal to $\tau(\Sha^2_S(K,M'))$, and let $\Uscr_v \subseteq X(K_v)$ be a $v$-adic neighborhood of $P_v$ for each $v \in S$. Since $\H^1(K_v,\SL_n) = 1$, each fibre $\phi^{-1}(P_v)$ has a $K_v$-point $Q_v$. In view of \eqref{eq:Weak2} (note that $\phi^\ast$ is an isomorphism), the family $(Q_v)_{v \in S} \in \prod_{v \in S} X_1(K_v)$ is orthogonal to $\tau_1(\Sha^2_S(K,M'))$. Then $(\psi(Q_v))_{v \in S} \in \prod_{v \in S} X_2(K_v)$ is orthogonal to $\tau_2(\Sha^2_S(K,M'))$. By Serre's generalized version of the implicit function theorem \cite[Part II, Chapter III, \S 10.2]{SerreLie}, we find for each $v \in S$ a $v$-adic neighborhood $\Vscr_v \subseteq X_2(K_v)$ of $\psi(Q_v)$ whose $K_v$-points can be lifted to $K_v$-points in $\phi^{-1}(\Uscr_v) \subseteq X_1(K_v)$.  Applying Theorem \ref{thm:MultiplicativeType} to $X_2 = M \backslash \SL_n$ yields a $K$-point $R \in X_2(K)$ belonging to $\prod_{v \in S} \Vscr_v$. Then the fibre $\psi^{-1}(R)$ contains a family $(Q'_v)_{v \in S} \in \prod_{v \in S}\phi^{-1}(\Uscr_v)$. By Lemma \ref{lm:Unipotent}, $\psi^{-1}(R)$ contains a $K$-point $Q \in \prod_{v \in S}\phi^{-1}(\Uscr_v)$. Then $\phi(Q) \in X(K)$ belongs to $\prod_{v \in S} \Uscr_v$.
    
    Since $X(K) \neq \varnothing$, any family $(P_v)$ is orthogonal to $\tau(\Sha^2(K,M'))$ (a subgroup consisting of everywhere locally constant classes). It follows that $X$ has weak approximation in $S$ whenever $\Sha^2_S(K,M') = \Sha^2(K,M')$. Conversely, suppose that $X$ has weak approximation in $S$. We show that it is also the case for $X_2$. Indeed, let $(R_v)_{v \in S} \in \prod_{v \in S} X_2(K_v)$ and let $\Vscr_v \subseteq X_2(K_v)$ be a $v$-adic neighborhood of $R_v$ for each $v \in S$. By Lemma \ref{lm:Unipotent}, each fibre $\psi^{-1}(R_v)$ contains a $K_v$-point $Q_v$. By Serre's generalized version of the implicit function theorem \cite[Part II, Chapter III, \S 10.2]{SerreLie}, we find for each $v \in S$ a $v$-adic neighborhood $\Uscr_v \subseteq X(K_v)$ of $\phi(Q_v)$ whose $K_v$-points can be lifted to $K_v$-points in $\psi^{-1}(\Vscr_v) \subseteq X_1(K_v)$. By our assumption on $X$, there is a $K$-point $P \in X(K)$ belonging to $\prod_{v \in S} \Uscr_v$. Then the fibre $\phi^{-1}(P)$ contains a family $(Q'_v)_{v \in S} \in \prod_{v \in S}\psi^{-1}(\Vscr_v)$. Since $G$ is special and has weak approximation, the fibre $\phi^{-1}(P)$ has a $K$-point $Q \in \prod_{v \in S}\psi^{-1}(\Vscr_v)$. Then $\psi(Q) \in X_2(K)$ belongs to $\prod_{v \in S} \Vscr_v$. Hence $X_2$ has weak approximation in $S$. By Theorem \ref{thm:MultiplicativeType}, one has $\Sha^2_S(K,M') = \Sha^2(K,M')$.
\end{proof}

\subsection{Alternative proofs} \label{subsec:Modified}

The idea of the ``fibration method'' at the end of the proofs of Theorems \ref{thm:Hasse} and \ref{thm:Weak} can be applied in an alternative way. They can be used to show that any homogeneous space of a special, $K$-rational algebraic group is $K$-stably birational to a $K$-torsor under a torus. If such a homogeneous space has a $K$-rational point, it is $K$-stably birational to a torus. This observation allows us to obtain Theorems \ref{thm:HasseModified} and \ref{thm:WeakModified} in this paragraph, which are variants of Theorems \ref{thm:Hasse} and \ref{thm:Weak} respectively (of course, they also imply Theorems \ref{customthm:Hasse} and \ref{customthm:Weak} respectively).

\begin{thm} [Theorem \ref{customthm:Hasse}] \label{thm:HasseModified}
    Let $K$ be the function field of a smooth projective geometrically integral curve $\Omega$ over a $p$-adic field $k$, and $X$ a homogeneous space of a special, $K$-rational algebraic group $G$, with geometric stabilizers $\ol{H}$ of type $\umult$. Let $M$ denote the natural $K$-form of the part of multiplicative type $\ol{M}$ of $\ol{H}$ and $M' = \wh{M} \otimes^{\Lbb} \Zbb(1)$. Then $X$ is $K$-stably birational to a $K$-torsor under a torus. Moreover, there exists a map
        \begin{equation} \label{eq:HasseModified1}
            \tau_1:\Sha^2_\omega(K,M') \to \tfrac{\H^3_{\nr}(K(X)/K,\Qbb/\Zbb(2))}{\Img \H^3(K,\Qbb/\Zbb(2))}
        \end{equation}
    with the following property. If there exists a family $(P_v)_{v \in \Omega^{(1)}}  \in \prod_{v \in \Omega^{(1)}} X(K_v)$ orthogonal to $\tau_1(\Sha^2(K,M'))$ relative to the pairing \eqref{eq:UnramifiedPairing}, then $X(K) \neq \varnothing$. In particular, the unramified first obstruction \eqref{eq:UnramifiedReciprocity} to the local-global principle for $X$ is the only one.
\end{thm}
\begin{proof}
    Consider the Springer $K$-lien $L_X$ and the Springer class $\eta_X \in \H^2(K,L_X)$. As in the proof of Theorem \ref{thm:Hasse}, we have a map $\H^2(K,L_X) \to \H^2(K,M)$. Let $\eta \in \H^2(K, M)$ denote the image of $\eta_X$ by this map. The first step is to find an embedding $\iota: M \hookrightarrow Q$ into a quasi-split torus such that $\iota_\ast \eta = 0 \in \H^2(K,Q)$. This can be done as follows. Let $L/K$ be a finite extension such that $X(L) \neq \varnothing$, then the restriction of $\eta_X$ to $\H^2(L,L_X)$ is neutral. Consequently, the restriction of $\eta$ to $\H^2(L,M)$ is $0$. Choose an embedding $i: M_L \hookrightarrow Q_1$ into a quasi-split $L$-torus, and let $\can: M \hookrightarrow \Res_{L/K}M_L$ be the canonical inclusion (where $\Res_{L/K}$ denotes the restriction of scalars {\em \`a la} Weil). The diagram
        \begin{equation*}
            \xymatrix{
                \H^2(K,M) \ar[r]^-{\can_\ast} \ar[rd]^{\res_{L/K}} & \H^2(K,\Res_{L/K}M_L) \ar[d]^{\cong} \ar[rr]^{(\Res_{L/K} i)_\ast} && \H^2(K,\Res_{L/K}Q_1) \ar[d]^{\cong} \\ & \H^2(L,M) \ar[rr]^{i_\ast} && \H^2(L,Q_1),
            }
        \end{equation*}
    where the vertical arrows are the isomorphisms from Shapiro's lemma, commutes. Indeed, its triangle commutes by \cite[Proposition 1.6.5]{NSW}. Let $Q = \Res_{L/K}Q_1$ (which is a quasi-split $K$-torus) and $\iota = (\Res_{L/K} i) \circ \can: M \hookrightarrow Q$, then $\iota_\ast \eta = 0$ as desired.

    Let $T$ be the cokernel of $\iota$ (it is a $K$-torus). By the long exact sequence associated with
        \begin{equation} \label{eq:HasseModified2}
            1 \to M \xrightarrow{\iota} Q \to T \to 1,
        \end{equation}
    we know that $\eta$ comes from $\H^1(K,T)$. If $\{m_{\sigma,\tau}\}_{\sigma,\tau}$ is a Galois $2$-cocycle representing $\eta$, then there is a $1$-cochain $\{q_{\sigma}\}_{\sigma}$ with coefficients in $Q(\ol{K})$ such that $\iota(m_{\sigma,\tau}) = q_\sigma \tensor[^\sigma]{q}{_\tau} q_{\sigma\tau}^{-1}$ for all $\sigma,\tau \in \Gamma_K$. Its image $\{t_\sigma\}_{\sigma}$ in $T$ is a $1$-cocycle, whose class $[t] \in \H^1(K,T)$ maps to $\eta \in \H^2(K,M)$. Let $Z = \tensor[_{-t}]{}{}T$ be the $K$-torsor under $T$ corresponding to the cocycle $-t$, that is, $Z(\ol{K}) = T(\ol{K})$ and the twisted Galois action on $Z(\ol{K})$ is given by
        \begin{equation} \label{eq:HasseModified3}
            \cdot_t: \Gamma_K \times Z(\ol{K}) \to Z(\ol{K}), \quad (\sigma, z)  \mapsto 
 \sigma \cdot_t z :=\tensor[^{\sigma}]{z}{} t_\sigma.
        \end{equation}
    The action of $T$ on $Z$ (by multiplication in $T(\ol{K}) = Z(\ol{K})$) makes $Z$ a homogeneous space of $Q$ with geometric stabilizers $\ol{M}$. We shall show that the Springer lien $L_Z$ is isomorphic to $\lien(M)$, and the Springer class $\eta_Z \in \H^2(K,M)$ is precisely $\eta$. To this end, we invoke the description of $L_Z$ and $\eta_Z$ in terms of cocycles (see \cite[\S 2.2.2]{DLA} or \cite[\S 5]{FSS}). Indeed, for each $\sigma \in \Gamma_K$, its action on $Q(\ol{K})$ restricts to the usual Galois action on $M(\ol{K})$, so that $L_Z = \lien(M)$. Next, fix the point $1 \in T(\ol{K}) = Z(\ol{K})$. Then \eqref{eq:HasseModified3} yields
        \begin{equation*}
            \sigma \cdot_t 1 = 1 t_\sigma = 1 \cdot q_\sigma, 
        \end{equation*}
    where $\cdot: Z \times Q \to Z$ denotes the action of $Q$ on $Z$ induced by that of $T$. We have $q_{\sigma}\tensor[^\sigma]{q}{_\tau} q_{\sigma\tau}^{-1} = \iota(m_{\sigma,\tau})$ for all $\sigma,\tau \in \Gamma_K$, hence $\eta_Z$ is represented by the $2$-cocycle $\{m_{\sigma,\tau}\}_{\sigma,\tau}$, {\em i.e.} $\eta_Z = \eta$. 
    
    To conclude, the map $\H^2(K,L_X) \to \H^2(K,M)$ sends $\eta_X$ to $\eta_Z$. Applying \cite[Theorem 3.4]{DLA}, we obtain a diagram
        \begin{equation} \label{eq:HasseModified4}
            \xymatrix{
                & X_1 \ar[ld]_{\phi} \ar[rd]^{\psi} \\
                X && Z,
            }
        \end{equation}
    where
    \begin{itemize}
        \item $X_1$ is a homogeneous space of $G \times_K Q$ with Springer lien $L_{X_1} \cong L_X$ and Springer class $\eta_{X_1} = \eta_X$,
        \item $\phi$ is a torsor under $Q$,
        \item the fibres of $\psi$ are homogeneous spaces of $G$ with geometric stabilizers $\Ker(\ol{H} \to \ol{M}) = :\ol{U}$.
    \end{itemize}
    Since $G$ is special, $\phi$ has a $K$-rational section. Since $G$ is $K$-rational, the extension $K(X_1)/K(X)$ is purely transcendental. Thus $X$ is $K$-stably birational to $X_1$. Next, by Lemma \ref{lm:Unipotent}, the generic fibre of $\psi$ is isomorphic to $U \backslash G_{K(Z)}$, where $U \subseteq G_{K(Z)}$ is a unipotent Zariski closed subgroup. Since both $U$ and $G_{K(Z)}$ are $K(Z)$-rational (for $U$, this is because its exponential map is a biregular isomorphism onto an affine space), the field extension $K(G)/K(X_1)$ and $K(G)/K(Z)$ are purely transcendental. It follows that $X_1$ (hence also $X$) is $K$-stably birational to $Z$. In particular, $\phi$ and $\psi$ induce isomorphisms between $\H^3_{\nr}(K(X)/K,\Qbb/\Zbb(2))$, $\H^3_{\nr}(K(X_1)/K,\Qbb/\Zbb(2))$, and $\H^3_{\nr}(K(Z)/K,\Qbb/\Zbb(2))$.

    In the course of establishing their theorem on the local-global principle for torsors under $K$-tori \cite[Theorem 5.1]{HShasse}, Harari and Szamuely constructed the first arrows in the composite
        \begin{equation} \label{eq:HasseModified5}
            \tau_2: \H^2(K,T') \to \tfrac{\H^3(Z,\Qbb/\Zbb(2))}{\Img \H^3(K,\Qbb/\Zbb(2))} \to \tfrac{\H^3(K(Z),\Qbb/\Zbb(2))}{\Img \H^3(K,\Qbb/\Zbb(2))},
        \end{equation}
    and Tian showed in his thesis that $\tau_2(\Sha^2_\omega(K,T'))$ lies in the subgroup $\tfrac{\H^3_{\nr}(K(Z)/K,\Qbb/\Zbb(2))}{\Img \H^3(K,\Qbb/\Zbb(2))}$ \cite[Corollary 1.4.5]{TianThese}. The map $\tau_2$ enjoys the following property. If there is a family of local points of $Z$ orthogonal to $\tau_2(\Sha^2(K,T'))$ relative to the pairing \eqref{eq:UnramifiedPairing}, then $Z(K) \neq \varnothing$. 
    
    Applying the functor $-\otimes^{\Lbb} \Zbb(1)$ to the exact sequence
        \begin{equation*} 
            0 \to \wh{T} \to \wh{Q} \to \wh{M} \to 0
        \end{equation*}
    dual to \eqref{eq:HasseModified2}, one obtains a distinguished triangle
        \begin{equation} \label{eq:HasseModified6}
            M' \to T' \to Q' \to M'[1].
        \end{equation}
    Since $\H^1(L,Q') = 0$ for any overfield $L/K$ and since $\Sha^2_\omega(K,Q') = 0$ by Lemma \ref{lm:QuasiTrivialTorusPAdic}, the long exact sequence associated with \eqref{eq:HasseModified6} yields an isomorphism 
        \begin{equation} \label{eq:HasseModified7}
            \Sha^2_\omega(K,M') \cong \Sha^2_\omega(K,T').
        \end{equation}
    More generally, one has $\Sha^2_S(K,M') \cong \Sha^2_S(K,T')$ for any finite set $S \subseteq \Omega^{(1)}$ (in particular, $\Sha^2(K,M') \cong \Sha^2(K,T')$). If we define the map $\tau_1$ in \eqref{eq:HasseModified1} as the composite 
        \begin{equation*}
            \Sha^2_\omega(K,M') \cong \Sha^2_\omega(K,T') \xrightarrow{\tau_2} \tfrac{\H^3_{\nr}(K(Z)/K,\Qbb/\Zbb(2))}{\Img \H^3(K,\Qbb/\Zbb(2))} \cong \tfrac{\H^3_{\nr}(K(X_1)/K,\Qbb/\Zbb(2))}{\Img \H^3(K,\Qbb/\Zbb(2))} \cong \tfrac{\H^3_{\nr}(K(X)/K,\Qbb/\Zbb(2))}{\Img \H^3(K,\Qbb/\Zbb(2))},
        \end{equation*}
    Then $\tau_1$ has the stated property. Indeed, suppose that there is a family $(P_v)_{v \in \Omega^{(1)}} \in \prod_{v \in \Omega^{(1)}}X(K_v)$ orthogonal to $\tau_1(\Sha^2(K,M')))$. Since $G$ is special, we may lift $(P_v)_{v \in \Omega^{(1)}}$ to a family $(Q_v)_{v \in \Omega^{(1)}} \in \prod_{v \in \Omega^{(1)}} X_1(K_v)$. Then $(\psi(Q_v))_{v \in \Omega^{(1)}} \in \prod_{v \in \Omega^{(1)}} Z(K_v)$ is orthogonal to $\tau_2(\Sha^2(K,T')))$, hence $Z(K) \neq \varnothing$. By Lemma \ref{lm:Unipotent}, one has $X_1(K) \neq \varnothing$, hence $X(K) \neq \varnothing$.
\end{proof}

For the problem of weak approximation, the above proof is actually simpler, since the involved homogeneous spaces already have $K$-rational points.

\begin{thm} [Theorem \ref{customthm:Weak}] \label{thm:WeakModified}
    Let $K$ be the function field of a smooth projective geometrically integral curve $\Omega$ over a $p$-adic field $k$, $H$ a $K$-group of type $\umult$, and $X = H \backslash G$ for some $K$-embedding $H \hookrightarrow G$ into a special, $K$-rational algebraic group $G$. Let $M$ denote the part of multiplicative type of $H$, and $M' = \wh{M} \otimes^{\Lbb} \Zbb(1)$.  Then $X$ is $K$-stably birational to a torus. Furthermore, let $\tau_1$ be the map $\Sha^2_\omega(K,M') \to \tfrac{\H^3_{\nr}(K(X)/K,\Qbb/\Zbb(2))}{\Img \H^3(K,\Qbb/\Zbb(2))}$ constructed as in \eqref{eq:HasseModified1}, and $S \subseteq \Omega^{(1)}$ a finite set of closed points. Then any family $(P_v)_{v \in S} \in \prod_{v \in S} X(K_v)$ orthogonal to $\tau_1(\Sha^2_S(K,M'))$ relative to the pairing \eqref{eq:UnramifiedPairingS} lies the closure of the diagonal image of $X(K)$. In particular, the reciprocity obstruction to weak approximation for $X$ is the only one. Moreover, $X$ has weak approximation in $S$ if and only if $\Sha^2_S(K,M') = \Sha^2(K,M')$.
\end{thm}
\begin{proof} [Proof (after J.-L. Colliot-Th\'el\`ene)]
    Following the proof of Theorem \ref{thm:HasseModified}, let $M \hookrightarrow Q$ be an embedding into a quasi-split torus. The embedding $H \hookrightarrow G$ and the composite $H \to M \hookrightarrow Q$ yield a diagonal embedding $H \hookrightarrow G \times_K Q$. This gives us a diagram
        \begin{equation*}
            \xymatrix{
                & X_1 \ar[ld]_{\phi} \ar[rd]^{\psi} \\
                X && T,
            }
        \end{equation*}
    where $X_1 = H \backslash (G \times_K Q)$, $T = Q/M$ (it is $K$-torus), $\phi$ is a torsor under $Q$, and the fibres of $\psi$ are homogeneous spaces of $G$ with geometric stabilizers $\Ker(\ol{H} \to \ol{M}) =: \ol{U}$. Note that $X$ is $K$-stably birational to $T$.
    
    In his thesis, Tian \cite[Proposition 1.3.1]{TianThese} showed that the restriction
        \begin{equation*}
            \Sha^2_\omega(K,T') \to \tfrac{\H^3_{\nr}(K(T)/K,\Qbb/\Zbb(2))}{\Img \H^3(K,\Qbb/\Zbb(2))}
        \end{equation*}
    of the map \eqref{eq:HasseModified5} coincides with the construction using ``flasque resolution'' by Harari, Scheiderer and Szamuely, which serves in the proof of their theorem on weak approximation for $K$-tori \cite[Theorem 5.2]{HSS}. Actually, combining this with Theorem 4.3(a) in {\em loc. cit.} gives us a more precise statement, that any family $(P_v)_{v \in S} \in \prod_{v \in S} T(K_v)$ orthogonal to $\tau_2(\Sha^2_S(K,T'))$ lies in the closure of the diagonal image of $T(K)$. Moreover, $T$ has weak approximation in $S$ if and only if $\Sha^2_S(K,T') = \Sha^2(K,T')$. By repeating the fibration argument as in the proof of Theorem \ref{thm:Weak}, one sees that the map $\tau_1$ defined in \eqref{eq:HasseModified1} has the stated property. Moreover, $X$ has weak approximation in $S$ if and only if $\Sha^2_S(K,M') = \Sha^2(K,M')$ (since $\Sha^2_S(K,T') \cong \Sha^2_S(K,M')$ and $\Sha^2(K,T') \cong \Sha^2(K,M')$).
\end{proof}

\subsection{Examples} \label{subsec:Example}

Let us discuss some corollaries to Theorem \ref{thm:Weak}. As always, let $K$ be the function field of a smooth projective geometrically integral curve $\Omega$ over a $p$-adic field $k$. Recall that when $M$ is a $K$-group of multiplicative type, $M' = M \otimes^{\Lbb} \Zbb(1)$ is quasi-isomorphic to a {\em $2$-term complex of tori} fitting into the distinguished triangle \eqref{eq:MQT}. If $M = F$ is finite abelian, then $M'$ is quasi-isomorphic to $F' = \iHom_K(F,\Qbb/\Zbb(2))$ (see Remark \ref{rmk:FiniteCase}). In what follows, if $H$ is a $K$-group of type $\umult$, then $M$ denotes its part of multiplicative type.

\begin{cor} \label{cor:Weak}
    Let $H$ be a $K$-group of type $\umult$ and $X = H \backslash \SL_n$ for some embedding $H \hookrightarrow \SL_n$. Then $X$ has weak approximation if and only if $\Sha^2_\omega(K,M') = \Sha^2(K,M')$.
\end{cor}

For example, by Lemma \ref{lm:QuasiTrivialTorusPAdic}, one has $\Sha^2_\omega(K,\Gbb_m) = 0$, hence $\Sha^2_\omega(K,\mu_n) = \Sha^2(K,\mu_n) 
= 0$ since $\Sha^2_\omega(K,\mu_n) \hookrightarrow \Sha^2_\omega(K,\Gbb_m)$ by the Kummer sequence and Hilbert's Theorem 90. Thus, the variety $\mu_n \backslash \SL_m$ (for any embedding $\mu_n \hookrightarrow \SL_m$) has weak approximation. By Lemma \ref{lm:WAForClassifyingSpaces}, the map $\H^1(K,\mu_n) \to \prod_{v \in S}\H^1(K_v,\mu_n)$ is surjective for all finite sets $S \subseteq \Omega^{(1)}$. Nonetheless, this follows easily from Kummer theory and the Artin--Whaples approximation theorem. 

\begin{prop} \label{prop:WAInInfinitePlaces}
    Let $H$ be a $K$-group of type $\umult$ and $X = H \backslash \SL_n$ for some embedding $H \hookrightarrow \SL_n$. There is an infinite set $S_0 \subseteq \Omega^{(1)}$ in which $X$ has weak approximation.
\end{prop}
\begin{proof}
    Let $L/K$ be a finite extension splitting the torus $T'$ in the triangle \eqref{eq:MQT}. It corresponds to a branched cover $f: \Omega' \to \Omega$ of smooth projective geometrically integral curves over $k$. For any non-empty open subset $U \subseteq \Omega$, a result of Poonen \cite[Corollary 2]{Poonen} assures the existence of a closed point $w \in f^{-1}(U)$ such that $k(w) = k(f(w))$. Hence, there are infinitely many points $w \in \Omega'$ having this property. If moreover $f$ is unramified over $f(w)$, then $L_w = K_{f(w)}$. Thus, the set $S_0$ of closed points $v \in \Omega$ for which there exists a point $w \in f^{-1}(v)$ with $L_w = K_v$ is infinite.

    We claim that $\Sha^2_S(K,M') = \Sha^2(K,M')$ for any finite set $S \subseteq S_0$ (which would conclude the proof by virtue of Theorem \ref{thm:Weak}). Indeed, let $\alpha \in \Sha^2_S(K,M')$, that is, $\loc_v(\alpha) = 0$ for any $v \notin S$. For $v \in S$, we have by definition of $S_0$ a closed point $w \in \Omega'$ lying over $v$ such that $L_w = K_v$.  We deduce from \eqref{eq:MapTau2} and Hilbert's Theorem 90 that $\Sha^2_\omega(L,M') = \Sha^1_\omega(L,T') \subseteq \H^1(L,T') = 0$, thus the restriction of $\alpha$ to $\H^2(L,M')$ is $0$. But then $\loc_v(\alpha) = 0$ since $\H^2(K_v,M') = \H^2(L_w,M')$. It follows that $\alpha \in \Sha^2(K,M')$.
\end{proof}

A theorem of Harbater \cite{Harbater} says that every finite group is a Galois group over $\Qbb_p(t)$. His original proof involves the technique of patching. There are several other proofs by Liu \cite{Liu}, Colliot-Th\'el\`ene \cite{CTgalois}, and Koll\'ar \cite{KollarFundamental,KollarRational}. As remarked by Colliot-Th\'el\`ene himself, the inverse Galois problem for a group $G$ (viewed as a finite constant group) over number fields can be reduced to the question of weak weak approximation (see below) on $G \backslash \SL_n$. Using this idea, we show

\begin{cor} \label{cor:InverseGalois}
    Any finite abelian group is a Galois group over $K$.
\end{cor}
\begin{proof}
    Let $F$ be a finite abelian group, which can be seen as a finite constant $K$-group scheme. By Proposition \ref{prop:WAInInfinitePlaces}, the variety $F \backslash \SL_n$ (for some embedding $F \hookrightarrow \SL_n$) has weak approximation in some infinite set $S_0 \subseteq \Omega^{(1)}$. By Lemma \ref{lm:WAForClassifyingSpaces}, this means the localization map $\H^1(K,F) \to \prod_{v \in S}\H^1(K_v,F)$ is surjective for every finite subset $S \subseteq S_0$. 

    Let $v: F \to S_0$ be any injective map. For each $x \in F$, choose a continuous homomorphism $c_x: \Gamma_{K_{v(x)}} \to F$ whose image contains $x$; this is possible since we have surjections $\Gamma_{K_{v(x)}} \twoheadrightarrow \Gamma_{k(v(x))} \twoheadrightarrow \wh{\Zbb}$ (the residue field of $k(v(x))$ being a finite field, its absolute Galois group is $\wh{\Zbb}$), then it is enough to inflate the continuous homomorphism $\wh{\Zbb} \to F$ mapping $1$ to $x$. Each homomorphism $c_x$ is an element of $\H^1(K_{v(x)},F)$. Then, there exists $c \in \H^1(K,F)$ such that $\loc_{v(x)}(c) = c_x$ for all $x \in F$. Thus $c$ is a continuous homomorphism $\Gamma_K \to F$ whose image contains every element of $F$, {\em i.e.} it is surjective. Its kernel is $\Gamma_L$ for some finite Galois extension $L$ of $K$, and $\Gal(L/K) = \Gamma_K / \Gamma_L \cong F$.
\end{proof}

Nevertheless, the abelian case of the {\em regular} inverse Galois problem over $\Qbb$ (hence over any field of characteristic $0$ by a base change argument) is known long before Harbater's theorem (see, {\em e.g.}, \cite[\S 4.2]{SerreTopic})\footnote{The author would like to thank Olivier Wittenberg for this remark.}.

We say that a smooth $K$-variety $X$ (with $X(K) \neq \varnothing$) satisfies the {\em weak weak approximation property} (resp. {\em countable weak weak approximation property}) if it has weak approximation away from a finite (resp. countable\footnote{The term ``countable'' means ``finite or countably infinite''.}) set $S_0 \subseteq \Omega^{(1)}$. This means $X$ has weak approximation in every finite set $S \subseteq \Omega^{(1)}$ with $S \cap S_0 = \varnothing$.

\begin{cor} \label{cor:WeakWeak}
    Let $H$ be a $K$-group of type $\umult$ and $X = H \backslash \SL_n$ for some embedding $H \hookrightarrow \SL_n$. Then $X$ satisfies the weak weak approximation property (resp. countable weak weak approximation property) if and only if $\Sha^2_\omega(K,M')$ is finite (resp. countable).
\end{cor}
\begin{proof}
    Suppose that $\Sha^2_\omega(K,M')$ is finite (resp. countable). Let $S_0$ be the set of closed points $v \in \Omega^{(1)}$ such that there exists $\alpha \in \Sha^2_\omega(K,M')$ with $\loc_v(\alpha) \neq 0$. Then $S_0$ is finite (resp. countable). For any finite set $S \subseteq \Omega^{(1)}$ with $S \cap S_0 = \varnothing$, one has $\Sha^2_S(K,M') = \Sha^2(K,M')$, thus $X$ has weak approximation in $S$ by Theorem \ref{thm:Weak}.

    Conversely, suppose that $X$ has weak approximation away from a finite (resp. countable) set $S_1 \subseteq \Omega^{(1)}$. By Theorem \ref{thm:Weak}, we have $\Sha^2_S(K,M') = \Sha^2(K,M')$ for any finite set $S \subseteq \Omega^{(1)}$ disjoint from $S_1$. We deduce from this the exactness of the bottom row of the diagram (induced by the distinguished triangle \eqref{eq:MQT})
        $$\xymatrix{
            0 \ar[r] & \Sha^1(K,T') \ar[r] \ar[d] & \Sha^1_\omega(K,T') \ar[r] \ar[d] & \bigoplus\limits_{v \in S_1} \H^1(K_v,T') \ar[d] \\
            0 \ar[r] & \Sha^2(K,M') \ar[r] & \Sha^2_\omega(K,M') \ar[r] & \bigoplus\limits_{v \in S_1} \H^2(K_v,M').
        }$$
   Indeed, any element $\alpha \in \Sha^2_\omega(K,M')$ lies in $\Sha^2_S(K,M')$ for some finite set $S \subseteq \Omega^{(1)}$. If the image of $\alpha$ in $\bigoplus_{v \in S_1}\H^2(K_v,M')$ is $0$, then $\alpha \in \Sha^2_{S \setminus S_1}(K,M') = \Sha^2(K,M')$. Since the left and the middle vertical arrows of the above diagram are isomorphisms by \eqref{eq:Sha2M}, a diagram chasing shows that the top row is also exact. Parts of local duality \eqref{eq:LocalDualityTorus} and global duality \eqref{eq:PoitouTateTorus} say that the groups $\Sha^1(K,T')$ and $\H^1(K_v,T')$ are finite, hence $\Sha^2_\omega(K,M') \cong \Sha^1_\omega(K,T')$ is finite (resp. countable).
\end{proof}

Note that in the case where $k$ is a number field and $F$ is a finite $\Gamma_k$-module, the defect of weak approximation on the quotient $F \backslash \SL_n$ (for some embedding $F \hookrightarrow \SL_{n,k}$) is explained by the group $\Sha^1_\omega(k,\wh{F})$. This group is always finite by an application of Chebotarev's theorem and the inflation-restriction sequence, hence weak weak approximation always holds. However, showing the finiteness (or the countability) of $\Sha^2_\omega(K,F')$ for a $p$-adic function field $K$ seems to be a difficult task in general. In fact, we have the following

\begin{prop} \label{prop:FailWeakWeak}
    There exists a finite Galois module over $K = \Qbb_p(t)$ such that $\Sha^2_\omega(K,F)$ is uncountable.
\end{prop}
\begin{proof}
    First, we find a $K$-torus $Q$ such that $\Sha^1_\omega(K,Q)$ is uncountable. It suffices to follow the proof of \cite[Proposition 4.5]{HSS}. Indeed, with the notations therein, one constructs a constant torus $Q$ with $\Sha^1_\omega(K,Q) = \H^1(K,Q)$, and, for each $b \in \Qbb_p = \Abb^1_{\Qbb_p}(\Qbb_p)$, an element $A_b \in \H^1(K,Q)$ whose ``residue'' at $b$ (see Remark 4.6 in {\em loc. cit.}) is non-zero. On the other hand, every element of $\H^1(K,Q)$ has non-zero residues at only finitely many points $b \in \Qbb_p$ (because it comes from some open subset $U \subseteq \Pbb^1_{\Qbb_p}$). Since $\Qbb_p$ is uncountable, if $\H^1(K,Q)$ was countable, there would exist $b \in \Qbb_p$ at which every element of $\H^1(K,Q)$ has vanishing residue, which is a contradiction. Thus $\Sha^1_\omega(K,Q) = \H^1(K,Q)$ must be uncountable.

    By Ono's Lemma \cite[Theorem 1.5.1]{Ono}, there are an integer $n \ge 1$, quasi-split $K$-tori $R, S$, and an isogeny $R \to Q^n \times_K S$. Let $F$ be the kernel of this isogeny. Since $\H^1(L,R) = \H^1(L,S) = 0$ for any overfield $L/K$ and $\Sha^2_\omega(K,R) = 0$ by Lemma \ref{lm:QuasiTrivialTorusPAdic}, the long exact sequence associated with $0 \to F \to R \to Q^n \times_K S \to 1$ gives $\Sha^2_\omega(K,F) \cong \Sha^1_\omega(K, Q^n \times_K S) \cong \Sha^1_\omega(K, Q)^n$, which is uncountable.
\end{proof}

Let $F$ be as in Proposition \ref{prop:FailWeakWeak}, $F' = \iHom_K(F,\Qbb/\Zbb(2))$, and $X = F' \backslash \SL_n$ for some embedding $F' \hookrightarrow \SL_n$. By Corollary \ref{cor:WeakWeak}, $X$ fails the countable weak weak approximation property, {\em a fortiori} it fails the weak weak approximation property and the weak approximation property. It would be interesting to see if weak weak approximation is strictly weaker than weak approximation, and if countable weak weak approximation is strictly weaker than weak weak approximation. By Corollaries \ref{cor:Weak} and \ref{cor:WeakWeak}, this is equivalent to the following questions.

\begin{question}
    Let $K$ be a $p$-adic function field. Does there exist a finite $\Gamma_K$-module $F$ such that $\Sha^2_\omega(K,F)$ is finite but $\Sha^2(K,F) \subsetneq \Sha_\omega^2(K,F)$?
\end{question}

\begin{question}
    Let $K$ be a $p$-adic function field. Does there exist a finite $\Gamma_K$-module $F$ such that $\Sha^2_\omega(K,F)$ is countably infinite?
\end{question}

Another closely related property is the {\em hyperweak approximation} \cite[\S 4]{HarariQuelques}. Let $F$ be a finite (not necessarily commutative) $K$-group, which extends to a finite \'etale group scheme $\Fcal \to U$ over a non-empty open subset $U \subseteq \Omega$. Then we say that $F$ satisfies the hyperweak approximation property if there is a non-empty open subset $V \subseteq U$ such that for every finite set $S \subseteq V^{(1)}$, the image of the localization $\H^1(K, F) \to \prod_{v \in S} \H^1(K_v, F)$ contains $\prod_{v \in S}\H^1(\Ocal_v,\Fcal)$. By Lemma \ref{lm:WAForClassifyingSpaces}, if the variety $F \backslash \SL_n$ (for some $K$-embedding $F \hookrightarrow \SL_n$) has weak weak approximation, then $F$ satisfies the hyperweak approximation property. It turns out that the converse holds when $F$ is commutative. 

\begin{lm}
    If $F$ is a finite $\Gamma_K$-module satisfying the hyperweak approximation property, then $F \backslash \SL_n$ (for some embedding $F \hookrightarrow \SL_n$)  satisfies the weak weak approximation property.
\end{lm}
\begin{proof}
    Indeed, if $F \backslash \SL_n$ fails the weak weak approximation property, then $\Sha^2_\omega(K, F')$ is infinite by Corollary \ref{cor:WeakWeak}. The same argument from its proof also gives an exact sequence
        \begin{equation*}
            0 \to \varinjlim_{S \subseteq V^{(1)}} \Sha^2_S(K,F') \to \Sha^2_\omega(K,F') \to \bigoplus_{v \notin V} \H^2(K_v,F')
        \end{equation*}
    for any non-empty open subset $V \subseteq U$. Since $\H^2(K_v,F')$ is finite for each $v \notin V$ and since $\Omega \setminus V$ is finite, $\varinjlim\limits_{S \subseteq V^{(1)}} \Sha^2_S(K,F')$ is infinite. On the other hand, there is an exact sequence
        \begin{equation*}
            \H^2(V,\Fcal') \to \prod_{v \notin V} \H^2(K_v,F') \times \prod_{v \in V^{(1)}}\H^2(\Ocal_v,\Fcal') \xrightarrow{\theta} \H^1(K,F)^D,
        \end{equation*}
    where the map $\theta$ is defined by $\theta((\alpha_v)_{v \in \Omega^{(1)}})(f) = \sum_{v \in \Omega^{(1)}} (\loc_v(f) \cup \alpha_v)$ (see the proof of \cite[Proposition 2.6]{IzquierdoDualite}). By the generalized Weil reciprocity law \eqref{eq:WeilReciprocity}, an element $\alpha \in \H^2(K,F')$ satisfies $\loc_v(\alpha) \in \H^2(\Ocal_v,\Fcal')$ for all $v \in V^{(1)}$ precisely when it comes from $\H^2(V,\Fcal')$. Since the group $\H^2(V,\Fcal')$ is finite, there exists a finite set $S \subseteq V^{(1)}$, a closed point $v_0 \in V^{(1)}$, and an element $\alpha \in \Sha^2_S(K,F')$ such that $\loc_{v_0}(\alpha) \notin \H^2(\Ocal_{v_0},\Fcal')$ (in particular, $v_0 \in S$). Since $\H^1(\Ocal_{v_0},\Fcal)$ and $\H^2(\Ocal_{v_0},\Fcal')$ are exact annihilators of each other under the cup-product pairing 
        \begin{equation*}
            \H^1(K_{v_0},F) \times \H^2(K_{v_0},F') \to \H^3(K_{v_0},\Qbb/\Zbb(2)) \cong \Qbb/\Zbb,
        \end{equation*}
    we find an element $f_{v_0} \in \H^1(\Ocal_{v_0},\Fcal)$ such that $f_{v_0} \cup \loc_{v_0}(\alpha) \neq 0$. Let $f_v = 0$ for $v \in S \setminus \{v_0\}$. Then the family $(f_v)_{v \in S} \in \prod_{v \in S}\H^1(\Ocal_v,\Fcal)$ is not orthogonal to $\Sha^2_S(K,F')$. In view of the exact sequence \eqref{eq:ExactSequenceFinite}, this family does not come from $\H^1(K, F)$, hence $F$ fails the hyperweak approximation property.
\end{proof}

The following vanishing result was communicated to the author by Jean-Louis Colliot-Th\'el\`ene. Recall that a $K$-group of multiplicative type $M$ is said to be {\em split} over a finite extension $L/K$ if $\Gamma_L$ acts trivially on $\wh{M}$. A finite group is said to be {\em metacyclic} if all of its Sylow subgroups are cyclic. 

\begin{prop} \label{prop:VanishingShaOmega}
    If $F$ is a finite $\Gamma_K$-module split by a metacyclic extension, then $\Sha^2_\omega(K,F) = 0$.
\end{prop}
\begin{proof}
    Let $L/K$ be a finite Galois extension splitting $F$, with metacyclic Galois group $G = \Gal(L/K)$. The ``coflasque resolution'' \cite[Lemma 0.6]{CTSflasque} provides an exact sequence
        \begin{equation*}
            0 \to \wh{T} \to \wh{Q} \to \wh{F} \to 0
        \end{equation*}
    of finitely generated $G$-modules, where $\wh{T}$ and $\wh{Q}$ are free as abelian groups, where $\wh{Q}$ is permutation, and where $\wh{T}$ is {\em coflasque} (that is, $\H^1(H,\wh{T}) = 0$ for all subgroups $H \subseteq G$). Since $G$ is metacyclic, a theorem of Endo--Miyata \cite[Theorem 1.5]{EM} says that $\wh{T}$ is a direct factor of a permutation module. By dualizing, we obtain an exact sequence
        \begin{equation*}
            1 \to F \to Q \to T \to 1,
        \end{equation*}
    where $Q$ is a quasi-split torus and $T$ is a direct factor of a quasi-split torus. Since $\H^1(K',T) = 0$ for any overfield $K'/K$ and $\Sha^2_\omega(K,Q) = 0$ by Lemma \ref{lm:QuasiTrivialTorusPAdic}, one has $\Sha^2_\omega(K,F) = 0$.
\end{proof}

\begin{cor}
    Let $F$ be a finite constant abelian group of exponent not divisible by 8 and $F' = \iHom_K(F,\Qbb/\Zbb(2))$. Then $\Sha^2_\omega(K,F') = 0$.
\end{cor}
\begin{proof}
    We may assume that $F = \Zbb/n$, where $n$ is either $2$, $4$, or a power of an odd prime. In all cases, the extension $K(\mu_n)/K$ is cyclic and splits $F' = \mu_n^{\otimes 2}$ (since $\wh{F'} = \mu_n^{\otimes (-1)}$), hence Proposition \ref{prop:VanishingShaOmega} implies $\Sha^2_\omega(K,\mu_n^{\otimes 2}) = 0$. 
\end{proof}

As for the group $\Zbb/2^m$, where $m \ge 3$, we cannot apply Proposition \ref{prop:VanishingShaOmega} unless the extension $K(\mu_{2^m})/K$ is cyclic. Nevertheless, we always have $\Sha^2_\omega(K(\sqrt{-1}),\mu_{2^m}^{\otimes 2}) = 0$ since the extension $K(\mu_{2^m})/K(\sqrt{-1})$ is cyclic. A restriction-corestriction argument yields

\begin{cor}
    The group $\Sha^2_\omega(K,\mu_{2^m}^{\otimes 2})$ is 2-torsion. It is trivial if the extension $K(\mu_{2^m})/K$ is cyclic (for example if $\sqrt{-1} \in K$).
\end{cor}

Thus, we are interested in the following question, to which a negative answer is expected.

\begin{question}
    Let $K$ be a $p$-adic function field. Is the group $\Sha^2_\omega(K,\mu_{2^m}^{\otimes 2})$ trivial for $m \ge 3$?
\end{question}
\section{Curves over higher-dimensional local fields} \label{sec:Higher}

In this section, we work with homogeneous spaces over function fields over higher-dimensional local fields. The main results here are Theorems \ref{customthm:HigherHasse} and \ref{customthm:HigherWeak}. Our approach is similar to that in paragraphs \ref{subsec:Hasse} and \ref{subsec:Weak}. The price we have to pay is that the obtained results are much coarser than those in the case of $p$-adic function fields. First, the geometric stabilizers are supposed to be finite abelian; the case of toric stabilizers is not treated because we do not have the corresponding duality theorems for tori at our disposal. Second, the constructed obstruction is not shown to be unramified.

\subsection{Construction of the obstruction} \label{subsec:HigherConstruction}

A single construction suffices for both local-global and weak approximation problems. Let $d \ge 0$ and let $\pi: X \to \Spec K$ be a smooth geometrically integral variety over a field $K$ of cohomological dimension $\cd(K) \le d+2$ and characteristic $0$. We construct a map
    \begin{equation} \label{eq:MapR}
        r: \H^{d+1}(K,\H^1(\ol{X},\Qbb/\Zbb(d+1))) \to \tfrac{\H^{d+2}(X,\Qbb/\Zbb(d+1))}{\Img\H^{d+2}(K,\Qbb/\Zbb(d+1))}
    \end{equation}
as follows. Consider the Hochschild--Serre spectral sequence 
    \begin{equation*}
        E_2^{p,q} = \H^p(K,\H^q(\ol{X},\Qbb/\Zbb(d+1))) \Rightarrow \H^{p+q}(X,\Qbb/\Zbb(d+1)).
    \end{equation*}
Under the assumption $\cd(K) \le d+2$, the outgoing differentials from $E_i^{d+1,1}$ vanish for all $i \ge 2$, hence there is a map $E_2^{d+1,1} \twoheadrightarrow E_{\infty}^{d+1,1} = \frac{F^{d+1} H^{d+2}}{F^{d+2} H^{d+2}} \subseteq \frac{H^{d+2}}{F^{d+2} H^{d+2}}$, where 
    \begin{equation*}
        H^{d+2} = F^0 H^{d+2} \supseteq F^1 H^{d+2} \supseteq \cdots \supseteq F^{d+2}H^{d+2}
    \end{equation*}
is a filtration of $H^{d+2} := \H^{d+2}(X,\Qbb/\Zbb(d+1))$. On the other hand, there is a surjection $E^{d+2,0}_2 \twoheadrightarrow E^{d+2,0}_\infty$, hence $F^{d+2}H^{d+2} = E^{d+2,0}_\infty = \Img(E^{d+2,0}_2 \to H^{d+2})$. Since $\ol{X}$ is integral, one has $\H^0(\ol{X},\Qbb/\Zbb(d+1)) = \Qbb/\Zbb(d+1)$, hence $E^{d+2,0} = \H^{d+2}(K,\Qbb/\Zbb(d+1))$, and the edge map $E^{d+2,0} \to \H^{d+2}$ is the natural pullback induced by $\pi$. Thus we obtain a map
    \begin{equation*}
        r: E^{d+1,1}_2 \to \tfrac{H^{d+2}}{\Img E^{d+2,0}_2}
    \end{equation*}
as in \eqref{eq:MapR}. It fits into the commutative diagram
    \begin{equation} \label{eq:MapR1}
        \xymatrix@C-1pc{
            \H^{d+2}(K,\Qbb/\Zbb(d+1)) \ar@{=}[d] \ar[r] & \H^{d+2}(K,\tau_{\le 1} \Rbb \pi_\ast \Qbb/\Zbb(d+1)) \ar[d] \ar[r] & \H^{d+1}(K,\H^1(\ol{X},\Qbb/\Zbb(d+1))) \ar[d]^{r} \\
            \H^{d+2}(K,\Qbb/\Zbb(d+1)) \ar[r] & \H^{d+2}(X,\Qbb/\Zbb(d+1)) \ar[r] & \frac{\H^{d+2}(X,\Qbb/\Zbb(d+1))}{\Img\H^{d+2}(K,\Qbb/\Zbb(d+1))}
        }
    \end{equation}
with exact rows, where the top row is associated with the distinguished triangle
    \begin{equation} \label{eq:MapR2}
        \Qbb/\Zbb(d+1) \to \tau_{\le 1} \Rbb \pi_\ast \Qbb/\Zbb(d+1) \to \H^1(\ol{X},\Qbb/\Zbb(d+1))[-1] \to \Qbb/\Zbb(d+1)[1]  
    \end{equation}
in $\Dcal^+(K)$. The middle vertical arrow in \eqref{eq:MapR1} is the composite 
    \begin{equation*}
        \H^{d+2}(K,\tau_{\le 1} \Rbb \pi_\ast \Qbb/\Zbb(d+1)) \to \H^{d+2}(K,\Rbb \pi_\ast \Qbb/\Zbb(d+1)) = \H^{d+2}(X,\Qbb/\Zbb(d+1)),
    \end{equation*}
where the last identification is due to the fact that $\Hbb(K,-) \circ \Rbb \pi_\ast = \Hbb(X,-)$ \cite[Corollary 10.8.3]{Weibel}.

To prove Theorem \ref{customthm:HigherWeak}, we shall need the following higher-dimensional generalization of Proposition \ref{prop:Descent2} (which is actually an analogue of \cite[Lemme 3.5.2]{CTSdescent} and \cite[Lemma 3]{SkorobogatovBeyond}).

\begin{prop} \label{prop:HigherWeak}
    Keep the above notations and assume in addition that $\ol{K}[X]^\times = \ol{K}^\times$. Let $F$ be a finite $\Gamma_K$-module and $\lambda: \wh{F} \to \Pic\ol{X}$ a $\Gamma_K$-equivariant homomorphism.
    Then the following claims hold.
    \begin{enumerate}
        \item \label{prop:HigherWeak1} $\lambda$ factors through a unique $\Gamma_K$-equivariant homomorphism $\lambda^{(0)}: \wh{F} \to \H^1(\ol{X},\Qbb/\Zbb(1))$. 

        \item \label{prop:HigherWeak2} Let $F' = \wh{F} \otimes \Qbb/\Zbb(d)  = \iHom_K(F,\Qbb/\Zbb(d+1))$ and let $\lambda^{(d)}: F' \to \H^1(\ol{X},\Qbb/\Zbb(d+1))$ denote the twist by $\Qbb/\Zbb(d)$ of the map $\lambda^{(0)}$ from \ref{prop:HigherWeak1}. Then, for any torsor $Y \to X$ of type $\lambda$ (see the discussion at the beginning of paragraph \ref{subsec:DescentRemark}), the diagram
            \begin{equation*}
                \xymatrix{
                    \H^{d+1}(K,F') \ar[r]^-{[Y] \cup \pi^\ast(-)} \ar[d]^{\lambda^{(d)}_\ast} & \H^{d+2}(X,\Qbb/\Zbb(d+1)) \ar[d] \\
                    \H^{d+1}(K,\H^1(\ol{X},\Qbb/\Zbb(d+1))) \ar[r]^-{r} & \tfrac{\H^{d+2}(X,\Qbb/\Zbb(d+1))}{\Img \H^{d+2}(K,\Qbb/\Zbb(d+1))}
                }
            \end{equation*}
        commutes. Here the cup-product $\H^{1}(X,F) \times \H^{d+1}(X,F') \to \H^{d+2}(X,\Qbb/\Zbb(d+1))$ is induced by the pairing $F \otimes F' \to \Qbb/\Zbb(d+1)$, and the map $r$ was constructed in \eqref{eq:MapR}.
    \end{enumerate}
\end{prop}
\begin{proof}
    We prove \ref{prop:HigherWeak1}. Since $\ol{K}[X]^\times = \ol{K}^\times$ is divisible, for each $n \ge 1$, the Kummer sequence yields an identification $\H^1(K,\mu_n) = \tensor[_n]{\Pic}{}\ol{X}$. Furthermore, one checks that for $n|m$, the map $\H^1(\ol{X},\mu_n) \to \H^1(\ol{X},\mu_m)$ (induced by the inclusion $\mu_n \hookrightarrow \mu_m$) is precisely the inclusion $\tensor[_n]{\Pic}{}\ol{X} \hookrightarrow \tensor[_m]{\Pic}{}\ol{X}$. Taking limit yields an identification $\H^1(\ol{X},\Qbb/\Zbb(1)) = (\Pic\ol{X})_{\tors}$. Since $\wh{F}$ is finite, $\lambda$ factors through a unique morphism $\lambda^{(0)}: \wh{F} \to \H^1(\ol{X},\Qbb/\Zbb(1))$. 

    Let us now show \ref{prop:HigherWeak2}. Let $\alpha \in \H^{d+1}(K,F') = \Ext^{d+1}_{K}(\Zbb,F')$ and let $\Zbb \to F'[d+1]$ be a morphism in $\Dcal^+(K)$ representing it. Such a morphism gives rise to the vertical arrows in the following commutative diagram in $\Dcal^+(\Ab)$:
        \begin{equation} \label{eq:HigherWeak1}
            \xymatrix{
                \Rbb\Hom_K(F', \tau_{\le 1} \Rbb \pi_\ast \Qbb/\Zbb(d+1))[1] \ar[r] \ar[d] & \Rbb\Hom_K(F', \H^1(\ol{X},\Qbb/\Zbb(d+1))) \ar[d] \\
                \Hbb(K,\tau_{\le 1} \Rbb \pi_\ast \Qbb/\Zbb(d+1))[d+2] \ar[r] & \Hbb(K,\H^1(\ol{X},\Qbb/\Zbb(d+1)))[d+1].
            }
        \end{equation}
    The horizontal arrows in \eqref{eq:HigherWeak1} are induced by triangle \eqref{eq:MapR2}. We claim that there is a commutative diagram
        \begin{equation} \label{eq:HigherWeak2}
            \xymatrix@C-1pc{
                \H^1(X,F) \ar@{=}[r] \ar[d]^{\cong} & \H^1(X,F) \ar[r]^-{\type} \ar[d]^{\cong} & \Hom_K(\wh{F},\Pic\ol{X}) \ar@{=}[d] \\
                \Ext^1_X(\wh{F},\Gbb_m) & \Ext^1_K(\wh{F},\tau_{\le 1} \Rbb \pi_\ast \Gbb_{m,X}) \ar[r] \ar[l]_-{\cong} & \Hom_K(\wh{F},\Pic\ol{X}) \\
                \Ext^1_X(\wh{F},\Qbb/\Zbb(1)) \ar[u]_{\gamma_1} \ar[d]^{- \otimes \Qbb/\Zbb(d)} & \Ext^1_K(\wh{F},\tau_{\le 1} \Rbb \pi_\ast \Qbb/\Zbb(1)) \ar[l]_-{\cong} \ar[r] \ar[u] \ar[d]^{- \otimes \Qbb/\Zbb(d)} & \Hom_K(\wh{F}, \H^1(\ol{X},\Qbb/\Zbb(1))) \ar[u]_{\cong} \ar[d]^{- \otimes \Qbb/\Zbb(d)} \\
                \Ext^1_X(F',\Qbb/\Zbb(d+1)) \ar[d]^{\sqdot \, \pi^\ast \alpha} & \Ext^1_K(F',\tau_{\le 1} \Rbb \pi_\ast \Qbb/\Zbb(d+1)) \ar[l]_-{\cong} \ar[r] \ar[d]^{\sqdot \, \alpha} & \Hom_K(F',\H^1(\ol{X},\Qbb/\Zbb(d+1))) \ar[d]^{\sqdot \, \alpha} \\ \H^{d+2}(X,\Qbb/\Zbb(d+1)) \ar@{=}[rd] & \H^{d+2}(K, \tau_{\le 1} \Rbb \pi_\ast \Qbb/\Zbb(d+1)) \ar[l] \ar[r] \ar[d] & \H^{d+1}(K, \H^1(\ol{X},\Qbb/\Zbb(d+1))) \ar[d]^r \\
                & \H^{d+2}(X,\Qbb/\Zbb(d+1)) \ar[r] & \tfrac{\H^{d+2}(X,\Qbb/\Zbb(d+1))}{\Img \H^{d+2}(K,\Qbb/\Zbb(d+1))},
            }
        \end{equation}
        where $\sqdot$ means the Yoneda product. For further use, we recall that for any $\Gamma_K$-module $M$ and any sheaf $\Fcal$ on $X_{\et}$, the distinguished triangle 
            \begin{equation*}
                \tau_{\le 1} \Rbb \pi_\ast \Fcal \to \Rbb \pi_\ast \Fcal \to \tau_{\ge 2} \Rbb \pi_\ast \Fcal \to (\tau_{\le 1} \Rbb \pi_\ast \Fcal)[1]
            \end{equation*}
        in $\Dcal^+(K)$ yields an isomorphism $\Ext^1_K(M,\tau_{\le 1} \Rbb \pi_\ast \Fcal) \xrightarrow{\cong} \Ext^1_K(M,\Rbb \pi_\ast \Fcal)$, since $\tau_{\ge 2} \Rbb \pi_\ast \Fcal$ is acyclic in degrees $0$ and $1$. But $\Ext^1_K(M,\Rbb \pi_\ast \Fcal) = \Ext^1_X(\pi^\ast M,\Fcal)$ since we have $\Rbb\Hom_K(M,-) \circ \Rbb \pi_\ast = \Rbb\Hom_X(\pi^\ast M,-)$ \cite[Corollary 10.8.3]{Weibel}. Thus, we obtain an isomorphism
            \begin{equation} \label{eq:HigherWeak3}
                \Ext^1_K(M, \tau_{\le 1} \Rbb \pi_\ast \Fcal) \xrightarrow{\cong} \Ext^1_X(\pi^\ast M, \Fcal).
            \end{equation}
        To see why \eqref{eq:HigherWeak2} commutes, let us consider the four rectangles of \eqref{eq:HigherWeak2} from the top to the bottom. As for the first rectangle, the right bottom horizontal arrow is induced by the map from the distinguished triangle
            \begin{equation*}
                \Gbb_m \to \tau_{\le 1} \Rbb \pi_\ast \Gbb_{m,X} \to \Pic\ol{X}[-1] \to \Gbb_m[1].
            \end{equation*}
        in $\Dcal^+(K)$ (here we use the fact that $\ol{K}[X]^\times = \ol{K}^\times$). The left vertical arrow is the isomorphism \eqref{eq:ExtOfGroupOfMultiplicativeType}. The left bottom horizontal arrow is the isomorphism \eqref{eq:HigherWeak3} for $M = \wh{F}$ and $\Fcal = \Gbb_m$. The middle vertical arrow is the obvious isomorphism that makes the left square commute. For the commutativity of the right square, see \cite[Appendix B]{HSdescent}. 
        
        As for the second rectangle, the left bottom horizontal arrow is the isomorphism \eqref{eq:HigherWeak3} for $M = \wh{F}$ and $\Fcal = \Qbb/\Zbb(1)$, the right bottom horizontal arrow is induced by the map from the distinguished triangle
            \begin{equation*}
                \Qbb/\Zbb(1) \to \tau_{\le 1} \Rbb \pi_\ast \Qbb/\Zbb(1) \to \H^1(\ol{X},\Qbb/\Zbb(1))[-1] \to \Qbb/\Zbb[1]
            \end{equation*}
        in $\Dcal^+(K)$ (here we use the fact that $\H^0(\ol{X},\Qbb/\Zbb(1)) = \Qbb/\Zbb(1))$ since $\ol{X}$ is integral). This rectangle obviously commutes. The right vertical arrow is an isomorphism since $\wh{F}$ is finite and since $\H^1(\ol{X},\Qbb/\Zbb(1)) = (\Pic\ol{X})_{\tors}$.

        The third rectangle obviously commutes. Its left bottom horizontal arrow is the isomorphism \eqref{eq:HigherWeak3} for $M = F'$ and $\Fcal = \Qbb/\Zbb(d+1)$ and its right bottom horizontal arrow is induced by the map from triangle \eqref{eq:MapR2}. As for the fourth rectangle, the left square obviously commutes (bearing in mind that $\H^{d+2}(X,\Qbb/\Zbb(d+1)) = \H^{d+2}(K,\Rbb\pi_\ast \Qbb/\Zbb(d+1))$ since $\Hbb(K,-) \circ \Rbb \pi_\ast = \Hbb(X,-)$ \cite[Corollary 10.8.3]{Weibel}), and the right square is obtained by taking cohomology of \eqref{eq:HigherWeak1}. Finally, the bottom square of \eqref{eq:HigherWeak2} is part of \eqref{eq:MapR1}; the triangle obviously commutes.

        By a similar argument as in the proof of Proposition \ref{prop:Descent2}, we have a commutative diagram
        \begin{equation} \label{eq:HigherWeak4}
            \xymatrix@C-1pc{
                & \H^1(X,\iHom_X(\wh{F},\Gbb_m)) \ar[r]^-{\cong} & \Ext^1_X(\wh{F},\Gbb_m)  \\
                & \H^1(X,\iHom_X(\wh{F},\Qbb/\Zbb(1))) \ar[u]_{\cong} \ar[r] & \Ext^1_X(\wh{F},\Qbb/\Zbb(1)) \ar@{=}[rd] \ar[u]^{\gamma_1} \\
                \H^1(X,F) \ar[r] \ar@/^2pc/[ruu]^-{\cong} \ar[ru]^-{\cong} \ar[rd] & \H^1(X,\iHom_X(\wh{F},F \otimes \wh{F})) \ar[u] \ar[d] \ar[r] & \Ext^1_X(\wh{F},F \otimes \wh{F}) \ar[d]^{-\otimes \Qbb/\Zbb(d)} \ar[u]^{\gamma_2} \ar[r]^-{\gamma_2} & \Ext^1_X(\wh{F},\Qbb/\Zbb(1)) \ar[d]^{-\otimes \Qbb/\Zbb(d)} \\
                & \H^1(X,\iHom_X(F',F \otimes F')) \ar[r] & \Ext^1_X(F', F \otimes F') \ar[r]^-{\gamma_3} & \Ext^1_X(F',\Qbb/\Zbb(d+1)).
            }
        \end{equation}
        
         Let $Y \to X$ be a torsor under $F$ of type $\lambda \in \Hom_K(\wh{F},\Pic\ol{X})$, and let $\varepsilon,\varepsilon'$ denote the respective images in $\Ext^1_X(\wh{F},F \otimes \wh{F})$ and $\Ext^1_X(F',F \otimes F')$ of $[Y] \in \H^1(X,F)$ by \eqref{eq:HigherWeak4}. The image of $[Y]$ in $\Ext^1_X(\wh{F},\Gbb_m)$ by \eqref{eq:HigherWeak2} is precisely $\gamma_1(\gamma_2(\varepsilon))$. The class $\gamma_2(\varepsilon) \in \Ext^1_X(\wh{F},\Qbb/\Zbb(1))$ corresponds to an element of $\Ext^1_K(\wh{F},\tau_{\le 1} \Rbb\pi_\ast \Qbb/\Zbb(1))$, whose image in $\Hom_K(\wh{F},\H^1(\ol{X},\Qbb/\Zbb(1)))$ is the homomorphism $\lambda^{(0)}$ from \ref{prop:HigherWeak1}. Moreover, the image of $\gamma_2(\varepsilon)$ in $\Ext^1(F',\Qbb/\Zbb(d+1))$ by \eqref{eq:HigherWeak2} is precisely $\gamma_3(\varepsilon')$. By exploiting the commutativity of \eqref{eq:HigherWeak2}, we see that the image of $\gamma_3(\varepsilon') \sqdot \pi^\ast \alpha \in \H^{d+2}(X,\Qbb/\Zbb(d+1))$ in $\tfrac{\H^{d+2}(X,\Qbb/\Zbb(d+1))}{\Img \H^{d+2}(K,\Qbb/\Zbb(d+1))}$ is precisely $r(\lambda_\ast^{(d)} \alpha)$. Finally, we have a commutative diagram of pairings
        \begin{equation*}
            \xymatrix@C-2pc{
                \H^1(X,F) \ar[d]  & \times & \H^{d+1}(X,F') \ar@{=}[d] \ar[rrrr]^-{\cup} &&&& \H^{d+2}(X,F \otimes F') \ar@{=}[d]  \\
                \Ext^1_X(F',F \otimes F') \ar[d]^{\gamma_3} & \times & \H^{d+1}(X,F') \ar@{=}[d] \ar[rrrr]^-{\sqdot} &&&& \H^{d+2}(X,F \otimes F') \ar[d] \\
                \Ext^1_X(F', \Qbb/\Zbb(d+1)) & \times & \H^{d+1}(X,F') \ar[rrrr]^-{\sqdot} &&&& \H^{d+2}(X,\Qbb/\Zbb(d+1)),
            }
        \end{equation*}    
        where $\sqdot$ means the Yoneda product, and where the top square commutes by \cite[Chapter V, Proposition 1.20]{MilneEtale}. This yields the identity $\gamma_3(\varepsilon') \sqdot \pi^\ast \alpha = [Y] \cup \pi^\ast \alpha \in \H^{d+2}(X,\Qbb/\Zbb(d+1))$, which proves \ref{prop:HigherWeak2}. 
\end{proof}

\subsection{The main theorems} \label{subsec:HigherProof}

We prove Theorems \ref{customthm:HigherHasse} and \ref{customthm:HigherWeak} in this paragraph. Let $k$ be a $d$-dimensional local field satisfying the condition \eqref{eq:Star} from page \pageref{eq:Star}, $\Omega$ a smooth projective geometrically integral curve over $k$, and $K = k(\Omega)$ its function field. Let $X$ be a homogeneous space of a simply connected semisimple linear algebraic group $G$ over $K$, with finite abelian geometric stabilizers $\ol{F}$. By Lemma \ref{lm:PicardGroupOfHomogeneousSpace} (and the discussion preceding it), $\ol{F}$ has a natural $K$-form $F$, and $\wh{F} = \Pic\ol{X}$ as $\Gamma_K$-modules. Recall that $\ol{K}[X]^\times = \ol{K}[G]^\times = \ol{K}^\times$ by Rosenlicht's lemma \cite[Proposition 3]{Rosenlicht}. Hence we may apply Proposition \ref{prop:HigherWeak}\ref{prop:HigherWeak2}, which says that $\id: \wh{F} \to \Pic\ol{X}$ yields an isomorphism $\id^{(d)}: F' \xrightarrow{\cong} \H^1(\ol{X},\Qbb/\Zbb(d+1))$, where $F' = \iHom_K(F,\Qbb/\Zbb(d+1))$. Let
    \begin{equation} \label{eq:MapRd}
        r^{(d)} = r \circ \id_{\ast}^{(d)}: \H^{d+1}(K,F') \to \tfrac{\H^{d+2}(X,\Qbb/\Zbb(d+1))}{\Img \H^{d+2}(K,\Qbb/\Zbb(d+1))},
    \end{equation}
where $r$ is the map from \eqref{eq:MapR}. We can now state the main results of this section.

\begin{thm} [Theorem \ref{customthm:HigherHasse}] \label{thm:HigherHasse}
    Let $K$ be the function field of a smooth projective geometrically integral curve $\Omega$ over a $d$-dimensional local field $k$ satisfying the condition \eqref{eq:Star} from page \pageref{eq:Star}, $G$ a special, simply connected semisimple algebraic group $G$ over $K$, and $X$ a homogeneous space of $G$ with finite abelian geometric stabilizers $\ol{F}$. Let $F$ be the natural $K$-form of $\ol{F}$, $F' = \iHom_K(F,\Qbb/\Zbb(d+1))$, and $r^{(d)}$ the map defined in \eqref{eq:MapRd}. If there is a point of $X(\Abb_K)$ orthogonal to $r^{(d)}(\Sha^{d+1}(K,F'))$ relative to the pairing \eqref{eq:AdelicPairing}, then $X(K) \neq\varnothing$. In particular, the first adelic obstruction \eqref{eq:AdelicReciprocity} to the local-global principle for $X$ is the only one.
\end{thm}

\begin{thm} [Theorem \ref{customthm:HigherWeak}] \label{thm:HigherWeak}
    Let $K$ be the function field of a smooth projective geometrically integral curve $\Omega$ over a $d$-dimensional local field $k$ satisfying the condition \eqref{eq:Star} from page \pageref{eq:Star}, $F$ a finite $\Gamma_K$-module, equipped with an embedding $F \hookrightarrow G$ into a special, simply connected semisimple algebraic group over $K$ that has weak approximation, and $X = F \backslash G$. Let $F' = \iHom_K(F,\Qbb/\Zbb(d+1))$, $r^{(d)}$ the map defined in \eqref{eq:MapRd}, and $S \subseteq \Omega^{(1)}$ a finite set. Then any family $(P_v)_{v \in S} \in \prod_{v \in S} X(K_v)$ orthogonal to $r^{(d)}(\Sha^{d+1}_S(K,F'))$ 
    relative to the pairing \eqref{eq:AdelicPairingS} lies in the closure of the diagonal image of $X(K)$. In particular, the generalized Brauer--Manin obstruction to weak approximation for $X$ is the only one. Moreover, $X$ has weak approximation in $S$ if and only if $\Sha^{d+1}_S(K,F') = \Sha^{d+1}(K,F')$.
\end{thm}

Theorem \ref{thm:HigherHasse} would follow from the following higher-dimensional generalization of Proposition \ref{prop:Descent1} (which is actually an analogue of \cite[Lemme 3.3.3]{CTSdescent} and \cite[(6.4)]{SkorobogatovTorsors}).

\begin{prop} \label{prop:HigherHasse}
    Keep the notations from Theorem \ref{thm:HigherHasse}. Let $\eta_X \in \H^2(K,F)$ be the Springer class of $X$ (see the discussion at the beginning of paragraph of \ref{subsec:Hasse}). Assume in addition that $X(\Abb_K) \neq \varnothing$ (hence $\eta_X \in \Sha^2(K,F)$). For any $\alpha \in \Sha^2(K,F')$, there holds
        \begin{equation*}
            \rho_X(r^{(d)}(\alpha)) = -\pair{\eta_X,\alpha}_{\PT}.
        \end{equation*}
    Here, the map $\rho_X$ was defined in \eqref{eq:AdelicReciprocity}, and  $\pair{-,-}_{\PT}$ is the pairing \eqref{eq:PoitouTateFinite}.
\end{prop}
\begin{proof}
    Let $e_X \in \Ext^2_K(\wh{F},\Gbb_m)$ denote the elementary obstruction of $X$ (see the discussion at the beginning of paragraph \ref{subsec:DescentRemark}). By Lemma \ref{lm:Hasse}, the isomorphism \eqref{eq:ExtOfGroupOfMultiplicativeType} sends $\eta_X$ to $e_X$. By the same argument as at the beginning of the proof of Proposition \ref{prop:Descent1}, $-e_X$ is represented by a morphism $\wh{F} \to \Gbb_m[2]$ associated with the distinguished triangle
        \begin{equation} \label{eq:HigherHasse1}
            \Gbb_m \to \tau_{\le 1} \Rbb\pi_\ast \Gbb_{m,X} \to \wh{F}[-1] \to \Gbb_m[1]
        \end{equation}
    in $\Dcal^+(K)$. Let $\pi^U: \Xcal \to U$ be an integral model of $X$ over some non-empty open subset $U \subseteq \Omega$. We may assume that $F$ extends to a locally constant finite \'etale $U$-group scheme $\Fcal$. Then $\wh{F}$ (resp. $F'$) extends to the locally constant finite \'etale $U$-group scheme $\wh{\Fcal} = \iHom_U(\Fcal,\Gbb_m)$ (resp. $\Fcal' = \wh{\Fcal} \otimes \Qbb/\Zbb(d) = \iHom(\Fcal,\Qbb/\Zbb(d+1))$). We have $\pi^U_\ast \Gbb_{m,\Xcal} = \Gbb_m$ (since the canonical map $\Gbb_m \to \pi^U_\ast \Gbb_{m,\Xcal}$ induces isomorphisms on the stalks of geometric points, by Rosenlicht's lemma \cite[Proposition 3]{Rosenlicht}). Similarly, $\wh{\Fcal} = \Rbb^1\pi_\ast^U \Gbb_{m,\Xcal}$ (by applying Lemma \ref{lm:PicardGroupOfHomogeneousSpace} to the stalks of geometric points). Thus, \eqref{eq:HigherHasse1} extends to a distinguished triangle 
        \begin{equation} \label{eq:HigherHasse2}
            \Gbb_m \to \tau_{\le 1} \Rbb\pi^U_\ast \Gbb_{m,X} \to \wh{\Fcal}[-1] \to \Gbb_m[1]
        \end{equation}
    in $\Dcal^+(U)$. The inverse class $\lambda^\ast e_U \in \Ext^2_U(\wh{\Fcal},\Gbb_m)$ of the morphism $\wh{\Fcal} \to \Gbb_m[2]$ associated with \eqref{eq:HigherHasse2} is a lifting of $\lambda^\ast e_X$. The spectral sequence $\H^p(U,\iExt^q_U(\wh{\Fcal},\Gbb_m)) \Rightarrow \Ext^{p+q}_U(\wh{\Fcal},\Gbb_m)$ provides an edge map
        \begin{equation*}
            \H^2(U,\Fcal) \cong \H^2(U,\iHom_U(\wh{\Fcal},\Gbb_m)) \to \Ext^2_U(\wh{\Fcal},\Gbb_m),
        \end{equation*}
    which is an isomorphism by \cite[Lemma 2.3.7]{SkorobogatovTorsors}. It maps a lifting $\eta_U$ of $\eta_X \in \H^2(K,F)$ to $\lambda^\ast e_U$. 
    
    Note that $\pi_\ast^U \Qbb/\Zbb(1) = \Qbb/\Zbb(1)$ since the fibres of $\pi^U$ are geometrically integral. Furthermore, by the smooth base change theorem \cite[Chapter VI, \S 4]{MilneEtale}, $\id^{(0)}: 
    \wh{F} \xrightarrow{\cong} \H^1(\ol{X},\Qbb/\Zbb(1))$ extends to an isomorphism $\id^{(0)}: \wh{\Fcal} \xrightarrow{\cong} \Rbb^1\pi_\ast^U \Qbb/\Zbb(1)$. Thus, we have a distinguished triangle
        \begin{equation} \label{eq:HigherHasse3}
            \Qbb/\Zbb(1) \to \tau_{\le 1} \Rbb\pi^U_\ast \Qbb/\Zbb(1) \to \wh{\Fcal}[-1] \to \Qbb/\Zbb(1)[1]
        \end{equation}
    in $\Dcal^+(U)$, extending the distinguished triangle
        \begin{equation*}
            \Qbb/\Zbb(1) \to \tau_{\le 1} \Rbb\pi_\ast \Qbb/\Zbb(1) \to \wh{F}[-1] \to \Qbb/\Zbb(1)[1]
        \end{equation*}
    in $\Dcal^+(K)$. To replace $\Gbb_m$ by $\Qbb/\Zbb(1)$, we need the following technical

\begin{lm} \label{lm:HigherHasse}
    The edge map $\H^2(U,\Fcal) \cong \H^2(U,\iHom_U(\wh{\Fcal},\Qbb/\Zbb(1))) \to \Ext^2_U(\wh{\Fcal},\Qbb/\Zbb(1))$,
    induced by the spectral sequence $\H^p(U,\iExt^q_U(\wh{\Fcal},\Qbb/\Zbb(1))) \Rightarrow \Ext^{p+q}_U(\wh{\Fcal},\Qbb/\Zbb(1))$, sends $\eta_U$ to the inverse class of the morphism $\wh{\Fcal} \to \Qbb/\Zbb(1)[2]$ associated with triangle \eqref{eq:HigherHasse3}.
\end{lm}
\begin{proof}
    Let $n \ge 1$ such that $n\Fcal = 0$, then $\wh{\Fcal} = \iHom_{U,\Zbb/n}(\Fcal, \mu_n)$, where the subscript ${}_{U,\Zbb/n}$ means we are taking $\iHom$ in the category of {\em $n$-torsion} sheaves over $U_{\et}$. In the spectral sequence $\H^p(U,\iExt^p_{U,\Zbb/n}(\wh{\Fcal},\mu_n)) \Rightarrow \Ext^{p+q}_{U,\Zbb/n}(\wh{\Fcal},\mu_n)$, one has $\iExt^p_{U,\Zbb/n}(\wh{\Fcal},\mu_n) = 0$ since we have $\mu_n \cong \Zbb/n$ (which is an injective $\Zbb/n$-module) locally for the \'etale topology on $U$. It follows that the induced edge map  $\H^2(U,\Fcal) \cong \H^2(U,\iHom_{U,\Zbb/n}(\wh{\Fcal},\mu_n)) \to \Ext^2_{U,\Zbb/n}(\wh{\Fcal},\mu_n)$ is an isomorphism. It follows that in the commutative diagram
        \begin{equation} \label{eq:eqHigherHasse4}
            \xymatrix{
                && \H^2(U,\Fcal) \ar[lld]_{\cong} \ar[ld] \ar[d] \ar[rd]^{\cong} \\
                \Ext^2_{U,\Zbb/n}(\wh{\Fcal},\mu_n) \ar[r] & \Ext^2_{U}(\wh{\Fcal},\mu_n) \ar[r] & \Ext^2_U(\wh{\Fcal},\Qbb/\Zbb(1)) \ar[r] & \Ext^2_U(\wh{\Fcal},\Gbb_m),
            }
        \end{equation}
    the composite of the bottom row is an isomorphism. On the other hand, the isomorphism $\wh{F} \xrightarrow{\cong} \H^1(\ol{X},\Qbb/\Zbb(1))$ factors through an isomorphism $\wh{F} \cong \H^1(\ol{X},\mu_n)$, which extends to an isomorphism $\wh{\Fcal} \xrightarrow{\cong} \Rbb^1 \pi^U_\ast \mu_n$ by the smooth base change theorem \cite[Chapter VI, \S 4]{MilneEtale}. Thus, we have a commutative diagram
        \begin{equation*}
            \xymatrix{
                \mu_n \ar[r] \ar[d] & \tau_{\le 1} \Rbb \pi^U_\ast \mu_n \ar[r] \ar[d] & \wh{\Fcal}[-1] \ar[r] \ar@{=}[d] & \mu_n[1] \ar[d] \\
                \Qbb/\Zbb(1) \ar[r] \ar[d] & \tau_{\le 1} \Rbb \pi^U_\ast \Qbb/\Zbb(1) \ar[d] \ar[r] & \wh{\Fcal}[-1] \ar[r] \ar@{=}[d] & \Qbb/\Zbb(1)[1] \ar[d] \\
                \Gbb_m \ar[r]  & \tau_{\le 1} \Rbb \pi^U_\ast \Gbb_m(1) \ar[r] & \wh{\Fcal}[-1] \ar[r]  & \Gbb_m[1]
            }
        \end{equation*}
    in $\Dcal^+(U)$, with distinguished rows (the middle row being \eqref{eq:HigherHasse3}). The existence of this diagram implies that the map $\H^2(U,\Fcal) \to \Ext^2_U(\wh{\Fcal},\Qbb/\Zbb(1))$ from \eqref{eq:eqHigherHasse4} effectively sends $\eta_U$ to the inverse class of the morphism $\wh{\Fcal} \to \Qbb/\Zbb(1)[2]$ associated with triangle \eqref{eq:HigherHasse3}.
\end{proof}

    Return to the proof of Proposition \ref{prop:HigherHasse}. Just like in the proof of Proposition \ref{prop:HigherWeak}, we have a commutative diagram
        \begin{equation} \label{eq:HigherHasse5}
            \xymatrix@C-1pc{
                & \H^2(U,\iHom_U(\wh{\Fcal},\Gbb_m)) \ar[r]^-{\cong} & \Ext^2_U(\wh{\Fcal},\Gbb_m)  \\
                & \H^2(U,\iHom_U(\wh{\Fcal},\Qbb/\Zbb(1))) \ar[u]_{\cong} \ar[r] & \Ext^2_U(\wh{\Fcal},\Qbb/\Zbb(1)) \ar@{=}[rd] \ar[u]^{\gamma_1} \\
                \H^2(U,\Fcal) \ar[r] \ar@/^2pc/[ruu]^-{\cong} \ar[ru] \ar[rd] & \H^2(U,\iHom_U(\wh{\Fcal},\Fcal \otimes \wh{\Fcal})) \ar[u] \ar[d] \ar[r] & \Ext^2_U(\wh{\Fcal},\Fcal \otimes \wh{\Fcal}) \ar[d]^{-\otimes \Qbb/\Zbb(d)} \ar[u]^{\gamma_2} \ar[r]^-{\gamma_2} & \Ext^2_U(\wh{\Fcal},\Qbb/\Zbb(1)) \ar[d]^{-\otimes \Qbb/\Zbb(d)} \\
                & \H^2(U,\iHom_U(\Fcal',\Fcal \otimes \Fcal')) \ar[r] & \Ext^2_U(\Fcal', \Fcal \otimes \Fcal') \ar[r]^-{\gamma_3} & \Ext^2_U(\Fcal',\Qbb/\Zbb(d+1)),
            }
        \end{equation}
    Let $\varepsilon_U,\varepsilon_U'$ denote the respective images of $\eta_U$ in $\Ext^2_U(\wh{\Fcal},\Fcal \otimes \wh{\Fcal})$ and $\Ext^2_U(\Fcal',\Fcal \otimes \Fcal')$ by \eqref{eq:HigherHasse5}. By Lemma \ref{lm:HigherHasse}, $-\gamma_2(\varepsilon_U)$ is the class of the the morphism $\Fcal \to \Qbb/\Zbb(1)[2]$ associated with triangle \eqref{eq:HigherHasse3}. The image in $\Ext^2(\Fcal',\Qbb/\Zbb(d+1))$ by \eqref{eq:HigherHasse5} of $\gamma_2(\varepsilon_U)$ is $\gamma_3(\varepsilon'_U)$, {\em i.e.  $-\gamma_3(\varepsilon'_U)$} is the class of the morphism $\Fcal' \to \Qbb/\Zbb(d+1)[2]$ from the distinguished triangle
        \begin{equation} \label{eq:HigherHasse6}
            \Qbb/\Zbb(d+1) \to \tau_{\le 1} \Rbb\pi^U_\ast \Qbb/\Zbb(d+1) \to \Fcal'[-1] \to \Qbb/\Zbb(d+1)[1]
        \end{equation}
    extending \eqref{eq:MapR2} (recall that $F' = \H^1(\ol{X},\Qbb/\Zbb(d+1))$). 
    
    Now, we have a commutative diagram of pairings
        \begin{equation} \label{eq:HigherHasse7}
            \xymatrix@C-2pc{
                \H^2(U,\Fcal) \ar[d]  & \times & \H^{d+1}_c(U,\Fcal') \ar@{=}[d] \ar[rrrr] &&&& \H^{d+3}_c(U,\Fcal \otimes \Fcal') \ar@{=}[d] \\
                \Ext^2_U(\Fcal',\Fcal \otimes \Fcal') \ar[d]^{\gamma_3} & \times & \H^{d+1}_c(U,\Fcal') \ar@{=}[d] \ar[rrrr]^-{\sqdot} &&&& \H^{d+3}_c(U,\Fcal \otimes \Fcal') \ar[d] \\
                \Ext^2_U(\Fcal', \Qbb/\Zbb(d+1)) & \times & \H^{d+1}_c(U,\Fcal') \ar[rrrr]^-{\sqdot} &&&& \H^{d+3}_c(U,\Qbb/\Zbb(d+1)) \cong \Qbb/\Zbb,
            }
        \end{equation}
     where $\sqdot$ means the Yoneda product, and where the top square commutes thanks to the construction of the cup-product (Artin--Verdier) pairing for cohomology with compact support. Let $\alpha \in \Sha^2(K,F')$, which lifts to an element $\alpha_U \in \H^{d+1}(U,\Fcal')$ (after shrinking $U$ if necessary). By the localization sequence \eqref{eq:LocalizationFinite}, $\alpha_U$ comes from an element $\alpha_U^c \in \H^{d+1}_c(U,\Fcal')$. The commutativity of \eqref{eq:HigherHasse7} implies $\gamma_3(\varepsilon'_U) \sqdot \alpha_U^c = \pair{\eta_U,\alpha_U^c}_{\AV}$. In other words, the map $\H^{d+1}_c(U,\Fcal') \to \H^{d+3}_c(U,\Qbb/\Zbb(d+1)) \cong \Qbb/\Zbb$, induced by triangle \eqref{eq:HigherHasse6}, sends $\alpha_U^c$ to $-\pair{\eta_U,\alpha_U^c}_{\AV}$, which is $-\pair{\eta_X,\alpha}_{\PT}$ by the construction of $\pair{-,-}_{\PT}$. Finally, by repeating the argument with the snake lemma construction at the end of the proof of Proposition \ref{prop:Descent2}, we see that the same map sends $\alpha_U^c$ to an element $\beta_U^c$, which would become $\rho_X(r^{(d)}(\alpha))$ after taking limit over $U$.
\end{proof}

\begin{proof}[Proof of Theorem \ref{thm:HigherHasse}]
    If there is a point of $X(\Abb_K)$ orthogonal to $r^{(d)}(\Sha^{d+1}(K,F'))$, then Proposition \ref{prop:HigherHasse} implies that $\pair{\eta_X,\alpha}_{\PT} = 0$ for all $\alpha \in \Sha^{d+1}(K,F')$. Since the Poitou--Tate pairing is non-degenerate, we have $\eta_X = 0$. Since $G$ is special, this is equivalent to $X(K) \neq \varnothing$.
\end{proof}

\begin{proof}[Proof of Theorem \ref{thm:HigherWeak}]
    We recall that the map $G \to X$ is a torsor under $F$ of type $\id: \wh{F} \to \Pic\ol{X}$. Let $(P_v)_{v \in S} \in \prod_{v \in S} X(K_v)$ be a family orthogonal to $r^{(d)}(\Sha^{d+1}_S(K,F'))$ relative to the pairing \eqref{eq:AdelicPairingS}. By Proposition \ref{prop:HigherWeak},  this is equivalent to
        \begin{equation*}
            \sum_{v \in S} [G](P_v) \cup \loc_v(\alpha) = \sum_{v \in S} ([G] \cup \pi^\ast \alpha)(P_v) = 0
        \end{equation*}
    for all $\alpha \in \Sha^{d+1}_S(K,M')$, where $\pi: X \to \Spec K$ is the structure morphism. Using the exact sequence \eqref{eq:ExactSequenceFinite}, we see that the family $([G](P_v))_{v \in S} \in \prod_{v \in S} \H^1(K_v,F)$ comes from $\H^1(K,F)$. By Lemma \ref{lm:WAForClassifyingSpaces}, $(P_v)_{v \in S}$ lies in the closure of the diagonal image of $X(K)$. The claim that $X$ has weak approximation in $S$ if and only if $\Sha^{d+1}_S(K,F') = \Sha^{d+1}(K,F')$ can be proved using the same argument as in the proof of Theorem \ref{thm:MultiplicativeType} (and by using \eqref{eq:ExactSequenceFinite} instead of \eqref{eq:ExactSequenceComplexTori}).
\end{proof}

We conclude by noting that both Theorems \ref{thm:HigherHasse} and \ref{thm:HigherWeak} fail for toric stabilizers. 

\begin{ex} \label{ex:CounterExampleToricStabilizer}
     Let $k = \Cbb\ps{t}$ (so that $d = 0$) and $K = k(\Omega)$ for some smooth projective geometrically integral $k$-curve. For a torus $T$ over $K$, (local and global) duality theorems between $T$ and $\wh{T}$ are established in \cite{CTH}. It is possible to construct such a torus $T$ and a homogeneous $X$ of $\SL_{n,K}$ with Springer lien $\lien(T)$ for which the Brauer--Manin obstruction attached to the group $\Sha^1(K,\wh{T})$ (resp. $\Sha^1_\omega(K,\wh{T})$) is insufficient to explain the failure of the local-global principle (resp. weak approximation). Indeed, an example in the case of the local-global principle is given at the beginning of \cite[\S 4]{ILA}. As for weak approximation, see Zhang's recent work \cite[Proposition 4.1]{Zhang}.
\end{ex}
	
\bibliographystyle{alpha}
\bibliography{ref}

\begin{thebibliography}{CTPS16}

\bibitem[BDH13]{BDH}
Mikhail~V. Borovoi, Cyril Demarche, and David Harari.
\newblock {Complexes de groupes de type multiplicatif et groupe de Brauer non
  ramifi\'e des espaces homog\`enes}.
\newblock {\em Annales scientifiques de l'\'Ecole Normale Sup\'erieure}, 4e
  s{\'e}rie, 46(4):651--692, 2013.

\bibitem[Blo86]{Bloch}
Spencer Bloch.
\newblock {Algebraic cycles and higher $K$-theory}.
\newblock {\em Advances in Mathematics}, 61(3):267--304, 1986.

\bibitem[BLR90]{BLR}
Siegfried Bosch, Werner L{\"u}tkebohmert, and Michel Raynaud.
\newblock {\em {N\'eron Models}}, volume~21 of {\em Ergebnisse der Mathematik
  und ihrer Grenzgebiete. 3. Folge / A Series of Modern Surveys in
  Mathematics}.
\newblock Springer Berlin Heidelberg, 1990.

\bibitem[BO74]{BlochOgus}
Spencer Bloch and Arthur Ogus.
\newblock {Gersten's conjecture and the homology of schemes}.
\newblock {\em Annales scientifiques de l'\'Ecole Normale Sup\'erieure}, 4e
  s{\'e}rie, 7(2):181--201, 1974.

\bibitem[Bor93]{BorovoiSecond}
Mikhail~V. Borovoi.
\newblock {Abelianization of the second nonabelian Galois cohomology}.
\newblock {\em Duke Mathematical Journal}, 72(1):217--239, 1993.

\bibitem[Bor96]{BorovoiFirst}
Mikhail~V. Borovoi.
\newblock {Abelianization of the first Galois cohomology of reductive groups}.
\newblock {\em International Mathematics Research Notices}, 1996(8):401--407,
  1996.

\bibitem[Bor98]{BorovoiReductive}
Mikhail~V. Borovoi.
\newblock {\em {Abelian Galois Cohomology of Reductive Groups}}, volume 132 of
  {\em Memoirs of the American Mathematical Society}.
\newblock American Mathematical Society, 1998.

\bibitem[BvH09]{BvHcrelle}
Mikhail~V. Borovoi and Joost van Hamel.
\newblock {Extended Picard complexes and linear algebraic groups}.
\newblock {\em Journal f\"ur die reine und angewandte Mathematik}, 627:53--82,
  2009.

\bibitem[BvH12]{BvHtrans}
Mikhail~V. Borovoi and Joost van Hamel.
\newblock {Extended equivariant Picard complexes and homogeneous spaces}.
\newblock {\em Transformation Groups}, 17:51--86, 2012.

\bibitem[Che54]{Chevalley}
Claude Chevalley.
\newblock {On algebraic group varieties}.
\newblock {\em Journal of the Mathematical Society of Japan}, 6(3-4):303--324,
  1954.

\bibitem[Che95]{Chernousov}
Vladimir~I. Chernousov.
\newblock {Galois cohomology and a weak approximation property for factor
  varieties $\mathbb{A}^n/G$}.
\newblock {\em Trudy Matematicheskogo Instituta imeni V.A. Steklova},
  208:335--349, 1995.
\newblock (in Russian).

\bibitem[Con12]{Conrad}
Brian Conrad.
\newblock {Weil and Grothendieck approaches to adelic points}.
\newblock {\em L'Enseignement Math\'ematique}, 58(1-2):61--97, 2012.

\bibitem[CT95]{CTgersten}
Jean-Louis Colliot-Th\'el\`ene.
\newblock {Birational invariants, purity and the Gersten conjecture}.
\newblock {\em Proc. Sympos. Pure Math.}, 58(1):1--64, 1995.

\bibitem[CT00]{CTgalois}
Jean-Louis Colliot-Th\'el\`ene.
\newblock {Rational connectedness and Galois covers of the projective line}.
\newblock {\em Annals of Mathematics}, 151(1):359--373, 2000.

\bibitem[CTH15]{CTH}
Jean-Louis Colliot-Th\'el\`ene and David Harari.
\newblock {Dualité et principe local-global pour les tores sur une courbe
  au-dessus de $\mathbb{C}(\!(t)\!)$}.
\newblock {\em Proceedings of the London Mathematical Society},
  110(6):1475--1516, 2015.

\bibitem[CTPS16]{CTPS}
Jean-Louis Colliot-Th\'el\`ene, Raman Parimala, and Venapally Suresh.
\newblock {Lois de r\'eciprocit\'e sup\'erieures et points rationnels}.
\newblock {\em Transactions of the American Mathematical Society},
  368(6):4219--4255, 2016.

\bibitem[CTS87a]{CTSdescent}
Jean-Louis Colliot-Th\'el\`ene and Jean-Jacques Sansuc.
\newblock {La descente sur les vari\'et\'es rationnelles, II}.
\newblock {\em Duke Mathematical Journal}, 54(1):375–492, 1987.

\bibitem[CTS87b]{CTSflasque}
Jean-Louis Colliot-Th\'el\`ene and Jean-Jacques Sansuc.
\newblock {Principal homogeneous spaces under flasque tori: Applications}.
\newblock {\em Journal of Algebra}, 106(1):148--205, 1987.

\bibitem[CTS07]{CTSrational}
Jean-Louis Colliot-Th\'el\`ene and Jean-Jacques Sansuc.
\newblock {The rationality problem for fields of invariants under linear
  algebraic groups (with special regards to the rationality problem)}.
\newblock In {\em Proceedings of the International Colloquium on Algebraic
  groups and Homogeneous Spaces (Mumbai 2004)}, pages 113--186. TIFR Mumbai,
  Narosa Publishing House, 2007.

\bibitem[Deb01]{Debarre}
Olivier Debarre.
\newblock {\em {Higher-Dimensional Algebraic Geometry}}.
\newblock Universitext. Springer New York, 2001.

\bibitem[DH19]{DH}
Cyril Demarche and David Harari.
\newblock {Artin--Mazur--Milne duality for fppf cohomology}.
\newblock {\em Algebra and Number Theory}, 13(10):2323–2357, 2019.

\bibitem[DLA19]{DLA}
Cyril Demarche and Giancarlo Lucchini~Arteche.
\newblock {Le principe de Hasse pour les espaces homog\`enes : r\'eduction au
  cas des stabilisateurs finis}.
\newblock {\em Compositio Mathematica}, 158(8):1568--1593, 2019.

\bibitem[Dou76]{Douai}
Jean-Claude Douai.
\newblock {\em {$2$-cohomologie galoisienne des groupes semi-simples}}.
\newblock PhD thesis, Universit\'e de Lille I, 1976.

\bibitem[EM75]{EM}
Shizuo Endo and Takehiko Miyata.
\newblock {On a classification of the function fields of algebraic tori}.
\newblock {\em Nagoya Mathematical Journal}, 56:85–104, 1975.

\bibitem[FSS98]{FSS}
Yuval Flicker, Claus Scheiderer, and Ramdorai Sujatha.
\newblock {Grothendieck's theorem on non-abelian $H^2$ and local-global
  principles}.
\newblock {\em Journal of the American Mathematical Society}, 11(3):731–750,
  1998.

\bibitem[Fu11]{Fu}
Lei Fu.
\newblock {\em {\'Etale Cohomology Theory}}, volume~13 of {\em Nankai Tracts in
  Mathematics}.
\newblock World Scientific Publishing Co. Pte. Ltd., Hackensack, NJ, 2011.

\bibitem[Gir71]{Giraud}
Jean Giraud.
\newblock {\em {Cohomologie non ab\'elienne}}, volume 179 of {\em Grundlehren
  der mathematischen Wissenschaften}.
\newblock Springer Berlin Heidelberg, 1971.

\bibitem[GS17]{GS}
Philippe Gille and Tam\'as Szamuely.
\newblock {\em {Central Simple Algebras and Galois Cohomology}}, volume 165 of
  {\em Cambridge University Press}.
\newblock World Scientific Publishing Co. Pte. Ltd., Hackensack, NJ, second
  edition, 2017.

\bibitem[Har87]{Harbater}
David Harbater.
\newblock {Galois coverings of the arithmetic line}.
\newblock In {\em Number Theory}, volume 1240 of {\em Lecture Notes in
  Mathematics}, pages 165--195. Springer Berlin Heidelberg, 1987.

\bibitem[Har07]{HarariQuelques}
David Harari.
\newblock {Quelques propri\'et\'es d'approximation reli\'ees \`a la cohomologie
  galoisienne d'un groupe alg\'ebrique fini}.
\newblock {\em Bulletin de la Soci\'et\'e' Math\'ematique de France},
  135(4):549--564, 2007.

\bibitem[Har17]{HarariCohomologie}
David Harari.
\newblock {\em {Cohomologie galoisienne et théorie du corps de classes}}.
\newblock Savoirs Actuels. EDP Sciences, 2017.

\bibitem[HS08]{HSmotive}
David Harari and Tam\'as Szamuely.
\newblock {Local-global principles for $1$-motives}.
\newblock {\em Duke Mathematical Journal}, 143(3):531–557, 2008.

\bibitem[HS13]{HSdescent}
David Harari and Alexei~N. Skorobogatov.
\newblock {Descent theory for open varieties}.
\newblock In {\em Torsors, \'Etale Homotopy and Applications to Rational
  Points}, volume 405 of {\em London Mathematical Society Lecture Note Series},
  pages 250--279. Cambridge University Press, 2013.

\bibitem[HS16]{HShasse}
David Harari and Tam\'as Szamuely.
\newblock {Local-global questions for tori over $p$-adic function fields}.
\newblock {\em Journal of Algebraic Geometry}, 25(3):571–605, 2016.

\bibitem[HSS15]{HSS}
David Harari, Claus Scheiderer, and Tam\'as Szamuely.
\newblock {Weak approximation for tori over $p$-adic function fields}.
\newblock {\em International Mathematics Research Notices},
  2015(10):2751--2783, 2015.

\bibitem[ILA21]{ILA}
Diego Izquierdo and Giancarlo Lucchini~Arteche.
\newblock {Local-global principles for homogeneous spaces over some
  two-dimensional geometric global fields}.
\newblock {\em Journal f\"ur die reine und angewandte Mathematik},
  781:165--186, 2021.

\bibitem[Izq14]{IzquierdoArxiv}
Diego Izquierdo.
\newblock {Th\'eor\`emes de dualit\'e pour les corps de fonctions sur des corps
  locaux sup\'erieurs et applications arithm\'etiques}, 2014.
\newblock preprint, arXiv:1405.2056v2 [math.AG],
  \url{https://arxiv.org/abs/1405.2056}.

\bibitem[Izq15]{IzquierdoI}
Diego Izquierdo.
\newblock {Principe local-global pour les corps de fonctions sur des corps
  locaux sup\'erieurs I}.
\newblock {\em Journal of Number Theory}, 157:250–270, 2015.

\bibitem[Izq16]{IzquierdoDualite}
Diego Izquierdo.
\newblock {Th\'eor\`emes de dualit\'e pour les corps de fonctions sur des corps
  locaux sup\'erieurs}.
\newblock {\em Mathematische Zeitschrift}, 284(1-2):615–642, 2016.

\bibitem[Izq17]{IzquierdoII}
Diego Izquierdo.
\newblock {Principe local-global pour les corps de fonctions sur des corps
  locaux sup\'erieurs II}.
\newblock {\em Bulletin de la Soci\'et\'e' Math\'ematique de France},
  145(2):267–293, 2017.

\bibitem[Kah12]{Kahn}
Bruno Kahn.
\newblock {Classes de cycles motiviques \'etales}.
\newblock {\em Algebra and Number Theory}, 6(7):1369--1407, 2012.

\bibitem[Kat86]{Kato}
Kazuya Kato.
\newblock {A Hasse principle for two dimensional global fields}.
\newblock {\em Journal f\"ur die reine und angewandte Mathematik},
  366:142--180, 1986.

\bibitem[Kle71]{Kleiman}
Steven Kleiman.
\newblock {Les theoremes de finitude pour le foncteur de Picard}.
\newblock In {\em Th\'eorie des Intersections et Th\'eor\`eme de Riemann--Roch
  : S\'eminaire de G\'eom\'etrie Alg\'ebrique du Bois Marie 1966 /67 (SGA 6)
  S\'eminaire Bourbaki}, volume 255 of {\em Lecture Notes in Mathematics},
  pages 616--666. Springer Berlin Heidelberg, 1971.

\bibitem[Kol00]{KollarFundamental}
J\'anos Koll\'ar.
\newblock {Fundamental groups of rationally connected varieties}.
\newblock {\em The Michigan Mathematical Journal}, 48(1):359--368, 2000.

\bibitem[Kol03]{KollarRational}
J\'anos Koll\'ar.
\newblock {Rationally connected varieties and fundamental groups}.
\newblock In {\em Higher Dimensional Varieties and Rational Points}, volume~12
  of {\em Bolyai Society Mathematical Studies}, pages 69--92. Springer Berlin
  Heidelberg, 2003.

\bibitem[LA14]{LA}
Giancarlo Lucchini~Arteche.
\newblock {Approximation faible et principe de Hasse pour des espaces
  homog\`enes \`a stabilisateur fini r\'esoluble}.
\newblock {\em Mathematische Annalen}, 360:1021--1039, 2014.

\bibitem[Lic87]{LichtenbaumConstruction}
Stephen Lichtenbaum.
\newblock {The construction of weight-two arithmetic cohomology}.
\newblock {\em Inventiones mathematicae}, 88:183--216, 1987.

\bibitem[Lic90]{LichtenbaumNew}
Stephen Lichtenbaum.
\newblock {New results on weight-two motivic cohomology}.
\newblock In {\em The Grothendieck Festschrift, Volume III: A Collection of
  Articles Written in Honor of the 60th Birthday of Alexander Grothendieck},
  volume~88 of {\em Modern Birkh\"auser Classics}, pages 35--55. Birkh\"auser
  Boston, 1990.

\bibitem[Liu95]{Liu}
Qing Liu.
\newblock {Tout groupe fini est un groupe de Galois sur $\mathbb{Q}_p(T)$,
  d'apr\`es Harbater}.
\newblock In {\em Recent Developments in the Inverse Galois Problem: A Joint
  Summer Research Conference on Recent Developments in the Inverse Galois
  Problem (July 17-23, 1993, University of Washington, Seattle)}, volume 186 of
  {\em Contemporary Mathematics}, pages 261--265. American Mathematical
  Society, 1995.

\bibitem[Man71]{Manin}
Yuri~I. Manin.
\newblock {Le groupe de Brauer--Grothendieck en g\'eom\'etrie diophantienne}.
\newblock In {\em Actes du Congr{\`e}s international des Math{\'e}maticiens
  (Nice, 1970)}, pages 401--411. Gauthier-Villars, Paris, 1971.

\bibitem[Mil80]{MilneEtale}
James~S. Milne.
\newblock {\em {Etale Cohomology}}, volume~33 of {\em Princeton Mathematical
  Series}.
\newblock Princeton University Press, 1980.

\bibitem[Mil06]{MilneDuality}
James~S. Milne.
\newblock {\em {Arithmetic Duality Theorems}}.
\newblock BookSurge, LLC, Charleston, SC, second edition, 2006.

\bibitem[Nag63]{Nagata}
Masayoshi Nagata.
\newblock {A generalization of the imbedding problem of an abstract variety in
  a complete variety}.
\newblock {\em Journal of Mathematics of Kyoto University}, 3(1):89--102, 1963.

\bibitem[Nis55]{Nishimura}
Hajime Nishimura.
\newblock {Some remark on rational points}.
\newblock {\em Memoirs of the College of Science, University Kyoto, Series A.
  Mathematics}, 29(2):189--192, 1955.

\bibitem[NSW08]{NSW}
J\"urgen Neukirch, Alexander Schmidt, and Kay Wingberg.
\newblock {\em Cohomology of Number Fields}, volume 323 of {\em Grundlehren der
  mathematischen Wissenschaften}.
\newblock Springer Berlin Heidelberg, 2008.

\bibitem[Ono61]{Ono}
Takashi Ono.
\newblock {Arithmetic of algebraic tori}.
\newblock {\em Annals of Mathematics}, 74(1):101--139, 1961.

\bibitem[Poo01]{Poonen}
Bjorn Poonen.
\newblock {Points having the same residue field as their image under a
  morphism}.
\newblock {\em Journal of Algebra}, 243(1):224--227, 2001.

\bibitem[Pop74]{Popov}
Vladimir~L. Popov.
\newblock {Picard group of homogeneous spaces of linear algebraic groups and
  one-dimensional vector bundles}.
\newblock {\em Mathematics of the USSR-Izvestiya}, 8(2):301--327, 1974.

\bibitem[Ros57]{Rosenlicht}
Maxwell Rosenlicht.
\newblock {Some rationality questions on algebraic groups}.
\newblock {\em Annali di Matematica Pura ed Applicata}, 43:25--50, 1957.

\bibitem[Ros96]{Rost}
Markus Rost.
\newblock {Chow groups with coefficients}.
\newblock {\em Documenta Mathematica}, 1:319--393, 1996.

\bibitem[San81]{Sansuc}
Jean-Jacques Sansuc.
\newblock {Groupe de Brauer et arithm\'etique des groupes alg\'ebriques
  lin\'eaires sur un corps de nombres}.
\newblock {\em Journal f\"ur die reine und angewandte Mathematik}, 327:12--80,
  1981.

\bibitem[Ser58]{SerreChow}
Jean-Pierre Serre.
\newblock {Espaces fibr\'es alg\'ebriques}.
\newblock In {\em Anneaux de Chow et Applications}, volume~3 of {\em
  S\'eminaire Claude Chevalley: 2e ann\'ee: 1958}. Secr\'etariat
  math\'ematique, Paris, 1958.
\newblock Expos\'e no. 1, p. 1-37.

\bibitem[Ser92]{SerreLie}
Jean-Pierre Serre.
\newblock {\em {Lie Algebras and Lie Groups: 1964 Lectures Given at Harvard
  University}}, volume 1500 of {\em Lecture Notes in Mathematics}.
\newblock Springer Berlin Heidelberg, 1992.

\bibitem[Ser94]{SerreGalois}
Jean-Pierre Serre.
\newblock {\em {Cohomologie Galoisienne: Cinqui\`eme \'edition, r\'evis\'ee et
  compl\'et\'ee}}, volume~5 of {\em Lecture Notes in Mathematics}.
\newblock Springer-Verlag Berlin Heidelberg, 1994.

\bibitem[Ser07]{SerreTopic}
Jean-Pierre Serre.
\newblock {\em {Topics in Galois Theory}}, volume~1 of {\em Research Notes in
  Mathematics}.
\newblock A K Peters/CRC Press, second edition, 2007.

\bibitem[Sko99]{SkorobogatovBeyond}
Alexei~N. Skorobogatov.
\newblock {Beyond the Manin obstruction}.
\newblock {\em Inventiones mathematicae}, 135:399--424, 1999.

\bibitem[Sko01]{SkorobogatovTorsors}
Alexei~N. Skorobogatov.
\newblock {\em {Torsors and Rational Points}}, volume 144 of {\em Cambridge
  Tracts in Mathematics}.
\newblock Cambridge University Press, 2001.

\bibitem[Tia20]{TianThese}
Yisheng Tian.
\newblock {\em {Arithm\'etique des groupes alg\'ebriques au-dessus du corps des
  fonctions d'une courbe sur un corps $p$-adique}}.
\newblock PhD thesis, Universit\'e Paris-Saclay, 2020.

\bibitem[Tia21]{TianWeak}
Yisheng Tian.
\newblock {Obstructions to weak approximation for reductive groups over
  $p$-adic function fields}.
\newblock {\em Journal of Number Theory}, 220:128–162, 2021.

\bibitem[Vos98]{Voskresenskii}
Valentin~E. Voskresenskii.
\newblock {\em {Algebraic Groups and Their Birational Invariants}}, volume 179
  of {\em Translations of Mathematical Monographs}.
\newblock American Mathematical Society, 1998.

\bibitem[Wei94]{Weibel}
Charles~A. Weibel.
\newblock {\em {An Introduction to Homological Algebra}}, volume~38 of {\em
  Cambridge Studies in Advanced Mathematics}.
\newblock Cambridge University Press, 1994.

\bibitem[Wit22]{Wittenberg}
Olivier Wittenberg.
\newblock Some aspects of rational points and rational curves.
\newblock In {\em Proceedings of the International Congress of Mathematicians,
  2022}, 2022.

\bibitem[Zha23]{Zhang}
Haowen Zhang.
\newblock {Weak approximation for homogeneous spaces over some two-dimensional
  geometric global fields}.
\newblock {\em Bulletin of the London Mathematical Society}, 55(1):308--320,
  2023.

\end{thebibliography}
\end{document}